\documentclass[10pt]{article}
\usepackage{amsfonts,color}
\usepackage{amssymb}
\usepackage{footnote}
\usepackage{mathtools}
\usepackage{amsthm,wasysym}
\usepackage[english]{babel}
\usepackage{bbm}
\usepackage{colonequals}
\usepackage{setspace}
\usepackage{titlesec}
\usepackage{float}

\usepackage{graphicx}
\usepackage{amsfonts}
\usepackage{hyperref}
\usepackage{subcaption}
\usepackage{amsmath}
\usepackage{nicefrac}
\usepackage{gensymb}
\usepackage{enumerate}
\usepackage[nameinlink,capitalize]{cleveref}
\usepackage{chngcntr}

\usepackage{cite}
\usepackage[nottoc,numbib]{tocbibind}

\usepackage{pgf,tikz,pgfplots}
\pgfplotsset{compat=1.14}
\usepackage{mathrsfs}
\usetikzlibrary{arrows}
\usetikzlibrary{arrows.meta}

\DeclareMathOperator{\at}{\bigg\vert}
\newcommand{\vb}[1]{\mathbf{#1}}
\newcommand{\bm}[1]{\boldsymbol{#1}}

\DeclareMathOperator{\sym}{\mathrm{sym}}
\DeclareMathOperator{\spa}{\mathrm{span}}

\DeclareMathOperator{\skw}{\mathrm{skew}}

\DeclareMathOperator{\tr}{\mathrm{tr}}

\newcommand{\one}{\bm{\mathbbm{1}}}

\newcommand{\dd}{\mathrm{d}}
\newcommand{\D}{\mathrm{D}}
\DeclareMathOperator{\di}{\mathrm{div}}
\DeclareMathOperator{\Di}{\mathrm{Div}}
\newcommand{\rot}[1]{\mathrm{curl}_\mathrm{2D}{#1}}
\newcommand{\srot}{\mathrm{div}\bm{R}}
\newcommand{\rog}{\bm{R} \nabla}
\DeclareMathOperator{\curl}{\mathrm{curl}}
\DeclareMathOperator{\Curl}{\mathrm{Curl}}

\newcommand{\so}{\mathfrak{so}}

\newcommand{\Ned}{\mathcal{N}_{I}}

\newcommand{\Nedtwo}{\mathcal{N}_{II}}
\newcommand{\Lag}{\mathcal{L}}
\newcommand{\Ber}{\mathcal{B}}
\newcommand{\Po}{\mathit{P}}
\newcommand{\Le}{\mathit{L}^2}

\newcommand{\Hone}{\mathit{H}^1}

\newcommand{\Hr}[1]{\mathit{H}(\srot{#1})}

\newcommand{\Hc}[1]{\mathit{H}(\mathrm{curl}{#1})}

\newcommand{\Hcz}[1]{\mathit{H}_0(\mathrm{curl}{#1})}

\newcommand{\body}{V}
\newcommand{\surf}{A}
\newcommand{\curv}{s}

\newcommand{\R}{\mathbb{R}}

\newcommand{\X}{\mathit{X}}
\newcommand{\C}{\mathit{C}}

\newcommand{\tem}{\mathcal{T}}
\newcommand{\ver}{\mathcal{V}}
\newcommand{\edge}{\mathcal{E}}
\newcommand{\face}{\mathcal{F}}
\newcommand{\cell}{\mathcal{C}}

\newcommand{\lame}{\lambda_{\mathrm{e}}}
\newcommand{\lammi}{\lambda_{\mathrm{micro}}}
\newcommand{\lamma}{\lambda_{\mathrm{macro}}}
\newcommand{\mue}{\mu_{\mathrm{e}}}
\newcommand{\muc}{\mu_{\mathrm{c}}}
\newcommand{\mumi}{\mu_{\mathrm{micro}}}
\newcommand{\muma}{\mu_{\mathrm{macro}}}
\newcommand{\Lc}{L_\mathrm{c}}
\newcommand{\Ce}{\mathbb{C}_{\mathrm{e}}}
\newcommand{\Cc}{\mathbb{C}_{\mathrm{c}}}
\newcommand{\Cmic}{\mathbb{C}_{\mathrm{micro}}}
\newcommand{\Cmac}{\mathbb{C}_{\mathrm{macro}}}
\newcommand{\J}{\mathbb{J}}

\newcommand{\wrt}{\text{w.r.t.}}

\newcommand{\Pm}{\bm{P}}
\newcommand{\ud}{\vb{u}}

\allowdisplaybreaks[1]

\makeindex

\newtheoremstyle{break}
{\topsep}{\topsep}%
{\itshape}{}%
{\bfseries}{}%
{\newline}{}%
\theoremstyle{break}

\newtheorem{theorem}{Theorem}

\newtheorem{remark}{Remark}
\newtheorem{observation}{Observation}
\newtheorem{definition}{Definition}

\counterwithin{theorem}{section}
\counterwithin{lemma}{section}
\counterwithin{remark}{section}
\counterwithin{observation}{section}
\counterwithin{definition}{section}

\counterwithin{figure}{section}
\counterwithin{equation}{section}

\makeatletter
\let\@fnsymbol\@arabic
\makeatother

\crefname{Problem}{Problem.}{Problem.}

\title{Higher order Bernstein-B\'ezier and N\'ed\'elec finite elements for the relaxed micromorphic model}
\author{\normalsize{Adam Sky}\thanks{Corresponding author: Adam Sky, Institute of Structural Mechanics, Statics and Dynamics, Technische Universit\"at Dortmund, August-Schmidt-Str. 8, 44227 Dortmund, Germany, email: adam.sky@tu-dortmund.de}
	, \quad
	\normalsize{Ingo Muench}\thanks{Ingo Muench, Institute of Structural Mechanics, Statics and Dynamics, Technische Universit\"at Dortmund, August-Schmidt-Str. 8, 44227 Dortmund, Germany, email: ingo.muench@tu-dortmund.de}
    , \quad
	\normalsize{Gianluca Rizzi}\thanks{Gianluca Rizzi, Institute of Structural Mechanics, Statics and Dynamics, Technische Universit\"at Dortmund, August-Schmidt-Str. 8, 44227 Dortmund, Germany, email: gianluca.rizzi@tu-dortmund.de}
	\quad and \quad
	\normalsize{Patrizio Neff}\thanks{Patrizio Neff,  \ \ Chair for Nonlinear 
		Analysis and Modelling, Faculty of Mathematics, Universit\"{a}t Duisburg-Essen,
		Thea-Leymann Str. 9, 45127 Essen, Germany, email: patrizio.neff@uni-due.de}
}

\begin{document}

\maketitle

\begin{abstract}
The relaxed micromorphic model is a generalized continuum model that is well-posed in the space $\X = [\Hone]^3 \times [\Hc{}]^3$.
Consequently, finite element formulations of the model rely on $\Hone$-conforming subspaces and N\'ed\'elec elements for discrete solutions of the corresponding variational problem. This work applies the recently introduced polytopal template methodology for the construction of N\'ed\'elec elements.
This is done in conjunction with Bernstein-B\'ezier polynomials and dual numbers in order to compute hp-FEM solutions of the model.
Bernstein-B\'ezier polynomials allow for optimal complexity in the assembly procedure due to their natural factorization into univariate Bernstein base functions.
In this work, this characteristic is further augmented by the use of dual numbers in order to compute their values and their derivatives simultaneously.
The application of the polytopal template methodology for the construction of the N\'ed\'elec base functions allows them to directly inherit the optimal complexity of the underlying Bernstein-B\'ezier basis.
We introduce the Bernstein-B\'ezier basis along with its factorization to univariate Bernstein base functions, the principle of automatic differentiation via dual numbers and a detailed construction of N\'ed\'elec elements based on Bernstein-B\'ezier polynomials with the polytopal template methodology.
This is complemented with a corresponding technique to embed  Dirichlet boundary conditions, with emphasis on the consistent coupling condition.
The performance of the elements is shown in examples of the relaxed micromorphic model.
\\
\vspace*{0.25cm}
\\
\end{abstract}

\section{Introduction}
One challenge that arises in the computation of materials with a pronounced micro-structure is the necessity of modelling the complex geometry of the domain as a whole, in order to correctly capture its intricate kinematics.
In other words, unit-cell geometries in metamaterials or various hole-shapes in porous media have to be accounted for in order to assert the viability of the model.
Naturally, this correlates with the resolution of the discretization in finite element simulations, resulting in longer computation times.   

The relaxed micromorphic model \cite{Neff2014} offers an alternative approach by introducing a continuum model with enriched kinematics, accounting for the independent distortion arising from the micro-structure.
As such, for each material point, the model introduces the microdistortion field $\Pm$ in addition to the standard displacement  field $\ud$. 
Consequently, each material point is endowed with twelve degrees of freedom, effectively turning into an affine-deformable micro-body with its own orientation.
In contrast to the classical micromorphic model \cite{Forest_full} by Eringen \cite{Eringen1999} and Mindlin \cite{Mindlin1964}, the relaxed micromorphic model does not employ the full gradient of the microdistortion $\D \Pm$ in its energy functional but rather its skew-symmetric part $\Curl \Pm$, designated as the micro-dislocation. Therefore, the micro-dislocation $\Curl \Pm$ remains a second-order tensor, whereas $\D \Pm$ is a third-order tensor.
Further, the model allows the transition between materials with a pronounced micro-structure and  homogeneous materials using the characteristic length scale parameter $\Lc$, which governs the influence of the micro-structure. 
In highly homogeneous materials the characteristic length scale parameter approaches zero $\Lc \to 0$, and for materials with a pronounced micro-structure its value is related to the size of the underlying unit-cell geometry.
Recent works demonstrate the effectiveness of the model in the simulation of band-gap metamaterials \cite{Madeo2016,Madeo2018,Agostino2020Band,BARBAGALLO2019148,DEMORE2022104995} and shielding against elastic waves \cite{Rizzi_shield,Leo, Rizzi2021Wave,ALBERDI2021104540}.
Furthermore, analytical solutions are already available for bending \cite{Rizzi_bending}, torsion \cite{Rizzi_torsion}, shear \cite{Rizzi_shear}, and extension \cite{Rizzi_extension} kinematics.

We note that the usage of the curl operator in the free energy functional directly influences the appropriate Hilbert spaces for existence and uniqueness of the related variational problem.
Namely, the relaxed micromorphic model is well-posed in $\{\ud, \Pm \} \in \X = [\Hone]^3 \times [\Hc{}]^3$ \cite{GNMPR15,Neff2015}, although the regularity of the microdistortion can be improved to $\Pm \in [\Hone]^{3\times 3}$ for certain smoothness of the data \cite{Knees,Reg}.
As shown in \cite{SKY2022115298}, the $\X$-space asserts well-posedness according to the Lax-Milgram theorem, such that $\Hone$-conforming subspaces and N\'ed\'elec elements \cite{Ned2,Nedelec1980,BERGOT20101937} inherit the well-posedness property as well. 

In this work we apply the polytopal template methodology introduced in \cite{skypoly} in order to construct higher order N\'ed\'elec elements based on Bernstein polynomials \cite{Ming} and apply the formulation to the relaxed micromorphic model.
Bernstein polynomials are chosen due to their optimal complexity property in the assembly procedure \cite{AinsworthOpt}.
We further enhance this feature by employing dual numbers \cite{Fike} in order to compute the values of the base functions and their derivatives simultaneously.
The polytopal template methodology allows to extend this property to the assembly of the N\'ed\'elec base functions, resulting in fast computations.
Alternatively, the formulation of higher order elements on the basis of Legendre polynomials can be found in \cite{Zaglmayr2006,Joachim2005,Solin}.
The construction of low order N\'ed\'elec elements can be found in \cite{Anjam2015,SkyOn} and specifically in the context of the 
the relaxed micromorphic model in \cite{Sky2021,Schroder2022,SKY2022115298,Sarhil2}.

This paper is structured as follows.
First, we introduce the relaxed micromorphic model and its limit cases with respect to the characteristic length scale parameter $\Lc$, after which we reduce it to a model of antiplane shear \cite{Voss2020}.
Next, we shortly discuss Bernstein polynomials and dual numbers for automatic differentiation.
The B\'ezier polynomial basis for triangles and tetrahedra is introduced, along with its factorization, highlighting its compatibility with dual numbers.
We consider a numerical example in antiplane shear for two-dimensional elements, a three-dimensional example for convergence of cylindrical bending, and a benchmark for the behaviour of the model with respect to the characteristic length scale parameter $\Lc$.
Lastly, we present our conclusions and outlook.

\hfill \break
The following definitions are employed throughout this work:
\begin{itemize}
    \item vectors are indicated by bold letters. Non-bold letters represent scalars;
    \item in general, formulas are defined using the Cartesian basis, where the base vectors are denoted by $\vb{e}_1$, $\vb{e}_2$ and $\vb{e}_3$;
    \item three-dimensional domains in the physical space are denoted with $\body \subset \R^3$. The corresponding reference domain is given by $\Omega$;
    \item analogously, in two dimensions we employ $\surf \subset \R^2$ for the physical domain and $\Gamma$ for the reference domain;
    \item curves on the physical domain are denoted by $\curv$, whereas curves in the reference domain by $\mu$;
    \item the tangent and normal vectors in the physical domain are given by $\vb{t}$ and $\vb{n}$, respectively.
    Their counterparts in the reference domain are $\bm{\tau}$ for tangent vectors and $\bm{\nu}$ for normal vectors.
\end{itemize}

\section{The relaxed micromorphic model}
The relaxed micromorphic model \cite{Neff2014} is governed by a free energy functional, incorporating the gradient of the displacement field $\D \ud$, the microdistortion $\Pm$ and the Curl of the microdistortion
	\begin{align}
		I(\ud, \Pm) = \dfrac{1}{2} \int_\body &\langle \sym(\D \ud - \Pm) , \, \Ce \sym(\D \ud - \Pm) \rangle + \langle \sym\Pm , \, \Cmic \sym \Pm \rangle \notag \\
		& + \langle \skw(\D \ud - \Pm) , \, \Cc \skw(\D \ud - \Pm) \rangle + \muma\Lc^2 \langle \Curl \Pm , \, \mathbb{L} \Curl \Pm \rangle \, \dd \body \notag \\ 
		& \qquad - \int_\body \langle \ud , \, \vb{f} \rangle + \langle \Pm , \, \bm{M} \rangle \, \dd \body \to  \min \quad \wrt \quad \{\ud, \Pm\}
		\, ,
	\end{align}
    where the Curl operator for second order tensors is defined row-wise as
    \begin{align}
    	\Curl \Pm &= \begin{bmatrix}
    		\curl(\begin{bmatrix}
    			P_{11} & P_{12} & P_{13}
    		\end{bmatrix}) \\
    	     \curl(\begin{bmatrix}
    	     	P_{21} & P_{22} & P_{23}
    	     \end{bmatrix}) \\
             \curl(\begin{bmatrix}
             	P_{31} & P_{32} & P_{33}
             \end{bmatrix}) 
    	\end{bmatrix} = \begin{bmatrix}
    	P_{13,y} - P_{12,z} & P_{11,z} - P_{13,x} & P_{12,x} - P_{11,y} \\
    	P_{23,y} - P_{22,z} & P_{21,z} - P_{23,x} & P_{22,x} - P_{21,y} \\
    	P_{33,y} - P_{32,z} & P_{31,z} - P_{33,x} & P_{32,x} - P_{31,y} 
    \end{bmatrix} \, , \notag \\  \curl \vb{p} &= \nabla \times \vb{p} \, , \qquad \vb{p}: \overline{\body} \subset \R^3 \to \R^3 \, ,
    \end{align}
    and $\curl(\cdot)$ is the vectorial curl operator.
	The displacement field and the microdistortion field are functions of the reference domain
	\begin{align}
		&\ud : \overline{\body} \subset \R^3 \to \R^3 \, , && \Pm : \overline{\body} \subset \R^3 \to \R^{3 \times 3} \, .
	\end{align}
	The tensors $\Ce, \Cmic,\mathbb{L} \in \R^{3\times 3 \times 3 \times 3}$ are standard positive definite fourth order elasticity tensors. For isotropic materials they take the form
	\begin{align}
		&\Ce = \lame \one \otimes \one + 2 \mue \, \J \, , &&
		\Cmic = \lammi \one \otimes \one + 2 \mumi \, \J \, .
	\end{align}
    where $\one$ is the second order identity tensor and $\J$ is the fourth order identity tensor.
	The fourth order tensor $\Cc \in \R^{3 \times 3 \times 3 \times 3}$ is a positive semi-definite material tensor related to Cosserat micro-polar continua and accounts for infinitesimal rotations $\Cc: \so(3) \to \so(3)$, where $\so(3)$ is the space of skew-symmetric matrices.
 
    For isotropic materials there holds $\Cc = 2\muc \, \J$, where $\muc \geq 0$ is called the Cosserat couple modulus. Further, for simplicity, we assume $\mathbb{L}=\J$ in the following.
	The macroscopic shear modulus is denoted by $\muma$ and $\Lc$ represents the characteristic length scale motivated by the geometry of the microstructure.
	The forces and micro-moments are given by $\vb{f}$ and $\bm{M}$, respectively.
	
\begin{figure}
		\centering
		\definecolor{asl}{rgb}{0.4980392156862745,0.,1.}
		\definecolor{asb}{rgb}{0.,0.4,0.6}
		\begin{tikzpicture}[line cap=round,line join=round,>=triangle 45,x=1.0cm,y=1.0cm]
			\clip(-0.5,-0.5) rectangle (16,4.5);
			
			\fill [asb, opacity=0.1] plot [smooth cycle] coordinates {(1,3) (3,4) (7, 2) (10,3) (12,1) (10,0) (5,1) (2,1)};
			
			\begin{scope}
				\clip(5,-0.5) rectangle (12.5,1.5);
				\draw [asl, dashed] plot [smooth cycle] coordinates {(1,3) (3,4) (7, 2) (10,3) (12,1) (10,0) (5,1) (2,1)};
			\end{scope}
		    \begin{scope}
		    	\clip(5,1.5) rectangle (12.5,4.5);
		    	\draw [asl, dashed] plot [smooth cycle] coordinates {(1,3) (3,4) (7, 2) (10,3) (12,1) (10,0) (5,1) (2,1)};
		    \end{scope}
			\begin{scope}
				\clip(0,-0.5) rectangle (5,4.5);
				\draw [asb] plot [smooth cycle] coordinates {(1,3) (3,4) (7, 2) (10,3) (12,1) (10,0) (5,1) (2,1)};
			\end{scope}
			
			\draw [-to,color=black,line width=1.pt] (0,0) -- (1,0);
			\draw [-to,color=black,line width=1.pt] (0,0) -- (0,1);
			\draw (1,0) node[color=black,anchor=west] {$x$};
			\draw (0,1) node[color=black,anchor=south] {$y$};
			
			\draw [-to,color=asb,line width=1.pt] (10.,1.2) -- (10.3,0.9);
			\draw [-to,color=asb,line width=1.pt] (9.6,1.2) -- (9.3,0.9);
			\draw [-to,color=asb,line width=1.pt] (9.8,1.2) -- (9.8,0.7);
			
			\draw [-to,color=asl,line width=1.pt] (11.35,2) -- (12,2.55);
			\draw (11.8,2.3) node[color=asl,anchor=north] {$\vb{n}$};
			\draw (9.8,1.2) node[color=asb,anchor=south] {$\vb{f}$};
			\draw (2.8,2.3) node[color=asb,anchor=south] {$\bm{M}$};
			\draw [-to,asb,domain=0:180,line width=1.pt] plot ({0.5*cos(\x-180)+2.8}, {0.5*sin(\x-180)+2.5});
			
			\draw (6.5,1.15) node[color=asb,anchor=south] {$\body$};
			\draw (4.3,3.5) node[color=asb,anchor=west] {$\surf_D=\surf^u_D=\surf^P_D$};
			\draw (10.6,2.95) node[color=asl,anchor=south] {$\surf_N =\surf_N^u= \surf_N^P$};
			
			\draw [black,domain=0:360,densely dotted] plot ({0.3*cos(\x)+11.3}, {0.3*sin(\x)+1});
			\draw [black,domain=0:360,densely dotted] plot ({0.78*cos(\x)+14}, {0.78*sin(\x)+1});
			\draw [color=black, densely dotted] (11.3,1.3) -- (14,1.78);
			\draw [color=black, densely dotted] (11.3,0.7) -- (14,0.22);
			
			\draw [color=black, line width=1] (13.9,1.3) -- (14.1,1.3) -- (14.1,1.1) -- (14.3,1.1) -- (14.3,0.9) -- (14.1,0.9) -- (14.1,0.7) -- (13.9,0.7) -- (13.9,0.9)
			-- (13.7,0.9) -- (13.7,1.1) -- (13.7,1.1) -- (13.9,1.1) -- (13.9,1.3);
			\draw [color=black, line width=1] (13.7,1.7) -- (13.9,1.7) -- (13.9,1.5) -- (14.1,1.5) -- (14.1,1.7) -- (14.3,1.7);
			\draw [color=black, line width=1] (14.7,1.3) -- (14.7,1.1) -- (14.5,1.1) -- (14.5,0.9) -- (14.7,0.9) -- (14.7,0.7);
			\draw [color=black, line width=1] (13.3,1.3) -- (13.3,1.1) -- (13.5,1.1) -- (13.5,0.9) -- (13.3,0.9) --  (13.3,0.7);
			\draw [color=black, line width=1] (13.7,0.3) -- (13.9,0.3) -- (13.9,0.5) -- (14.1,0.5) -- (14.1,0.3) -- (14.3,0.3);
			
			\fill [asb,domain=0:360, opacity=0.1] plot ({0.78*cos(\x)+14}, {0.78*sin(\x)+1});
			
			\fill[white] (13.9,1.7) -- (13.9,1.5) -- (14.1,1.5) -- (14.1,1.7) -- cycle;
			
			\fill[white] (14.7,1.1) -- (14.5,1.1) -- (14.5,0.9) -- (14.7,0.9)  -- cycle;
			
			\fill[white] (13.3,1.1) -- (13.5,1.1) -- (13.5,0.9) -- (13.3,0.9) -- cycle;
			
			\fill[white] (13.9,0.3) -- (13.9,0.5) -- (14.1,0.5) -- (14.1,0.3) -- cycle;
			
			\fill[white] (13.9,1.3) -- (14.1,1.3) -- (14.1,1.1) -- (14.3,1.1) -- (14.3,0.9) -- (14.1,0.9) -- (14.1,0.7) -- (13.9,0.7) -- (13.9,0.9)
			-- (13.7,0.9) -- (13.7,1.1) -- (13.7,1.1) -- (13.9,1.1) -- cycle;
			
			\begin{scope}
				\clip(13.3,1.3) rectangle (13.1,0.7);
				\fill [white,domain=0:360] plot ({0.78*cos(\x)+14}, {0.78*sin(\x)+1});
			\end{scope}
		    \begin{scope}
		    	\clip(13.7,1.7) rectangle (14.3,2.3);
		    	\fill [white,domain=0:360] plot ({0.78*cos(\x)+14}, {0.78*sin(\x)+1});
		    \end{scope}
	        \begin{scope}
	        	\clip(13.7,0.3) rectangle (14.3,-0.3);
	        	\fill [white,domain=0:360] plot ({0.78*cos(\x)+14}, {0.78*sin(\x)+1});
	        \end{scope}
            \begin{scope}
            	\clip(14.7,1.3) rectangle (14.9,0.7);
            	\fill [white,domain=0:360] plot ({0.78*cos(\x)+14}, {0.78*sin(\x)+1});
            \end{scope}
		\end{tikzpicture}
		\caption{The domain in the relaxed micromorphic model with Dirichlet and Neumann boundaries under internal forces and micro-moments. The Dirichlet boundary of the microdistortion is given by the consistent coupling condition. The model can capture the complex kinematics of an underlying micro-structure.}
		\label{fig:domain}
	\end{figure}
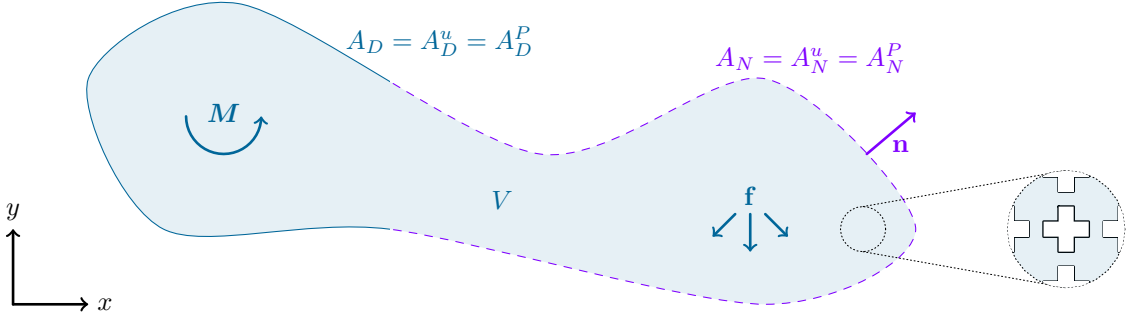

	Equilibrium is found at minima of the energy functional, which is strictly convex (also for $\Cc \equiv 0$). As such, we consider variations with respect to its parameters, namely the displacement and the microdistortion.
	Taking variations of the energy functional with respect to the displacement field $\ud$ yields
	\begin{align}
		\delta_u I = \int_\body \langle \sym \D \delta \ud , \, \Ce \sym(\D \ud - \Pm) \rangle + \langle \skw \D \delta \ud , \, \Cc \skw(\D \ud - \Pm) \rangle - \langle \delta \ud , \, \vb{f} \rangle \, \dd \body = 0 \, .
		\label{eq:weak_disp}
	\end{align}
	The variation with respect to the microdistortion $\Pm$ results in
	\begin{align}
		\delta_P I = \int_\body & \langle \sym \delta \Pm , \, \Ce \sym(\D \ud - \Pm) \rangle + \langle \skw \delta \Pm , \, \Cc \skw(\D \ud - \Pm) \rangle \notag \\ &- \langle \sym \delta \Pm , \, \Cmic \sym \Pm \rangle - \muma \Lc^2 \langle \Curl \delta \Pm , \, \Curl \Pm \rangle + \langle \delta \Pm , \, \bm{M} \rangle \, \dd \body = 0 \, .
		\label{eq:weak_p}
	\end{align}
	From the total variation we extract the bilinear form
	\begin{align}
		a(\{\delta \ud , \delta \Pm\},\{\ud, \Pm\}) = \int_\body & \langle \sym(\D \delta \ud - \delta \Pm) , \, \Ce \sym(\D \ud - \Pm) \rangle + \langle \sym \delta \Pm, \, \Cmic \sym \Pm \rangle \notag \\
		& + \langle \skw(\D \delta \ud - \delta \Pm) , \, \Cc \skw(\D \ud - \Pm) \rangle + \muma \Lc^2 \langle \Curl \delta \Pm , \, \Curl \Pm \rangle \, \dd \body \, ,
		\label{eq:bi_full}
	\end{align}
	and linear form of the loads
	\begin{align}
		l(\{\delta \ud , \delta \Pm\}) = \int_\body \langle \delta \ud , \, \vb{f} \rangle + \langle \delta \Pm , \, \bm{M} \rangle \, \dd \body \, .
		\label{eq:li_full}
	\end{align}
	Applying integration by parts to \cref{eq:weak_disp} yields
	\begin{align}
		\int_{\partial \body} &\langle \delta \vb{u} \, , [\Ce \sym (\D \ud - \Pm) + \Cc \skw(\D \ud - \Pm)] \,\vb{n} \rangle \, \dd \surf \notag \\ &-  \int_\body \langle \delta \ud \, , \Di[\Ce \sym (\D \ud - \Pm) + \Cc \skw(\D \ud - \Pm)] -  \vb{f} \rangle \, \dd \body = 0 \, .
		\label{eq:du}
	\end{align}
	Likewise, integration by parts of \cref{eq:weak_p} 
	results in
	\begin{align}
		\int_\body \langle &\delta \Pm , \, \Ce \sym(\D \ud - \Pm) + \Cc \skw(\D \ud - \Pm) - \Cmic \sym \Pm - \muma \Lc^2 \Curl \Curl \Pm + \bm{M} \rangle \, \dd \body \notag \\
		& - \muma \Lc^2 \int_{\partial \body} \langle \delta \Pm , \, \Curl \Pm \times \vb{n} \rangle \, \dd \surf = 0 \, .
		\label{eq:dp}
	\end{align}
	The strong form is extracted from \cref{eq:du} and \cref{eq:dp} by splitting the boundary
	\begin{align}
		&\surf = \surf_D \cup \surf_N \, , && \surf_D \cap \surf_N = \emptyset \, ,
	\end{align}
	into a Dirichlet boundary with embedded boundary conditions and a Neumann boundary with natural boundary conditions, such that no tractions are imposed on the Neumann boundary   
	\begin{subequations}
		\begin{align}
			-\Di[\Ce \sym (\D \vb{u} - \bm{P}) + \Cc \skw (\D \vb{u} - \bm{P})] &= \vb{f} && \text{in} \quad \body \, , \label{eq:strong_u} \\
			-\Ce  \sym (\D \vb{u} - \Pm) - \Cc  \skw(\D \vb{u} - \Pm) + \Cmic \sym \Pm + \muma \, \Lc ^ 2  \Curl\Curl\Pm &= \bm{M} && \text{in} \quad \body \, , \label{eq:strong_p}  \\
			\vb{u} &= \widetilde{\vb{u}} && \text{on} \quad \surf_D^u \, , \\
			\Pm \times \, \vb{n} &= \widetilde{\Pm} \times \, \vb{n} && \text{on} \quad \surf_D^P \, , \label{eq:pdir} \\
			[\Ce \sym (\D \vb{u}- \bm{P}) + \Cc \skw (\D \vb{u} - \bm{P})] \, \vb{n} &= 0 && \text{on} \quad \surf_N^u \, ,\\
			\Curl \Pm \times \, \vb{n}  &= 0  && \text{on} \quad \surf_N^P \, .
		\end{align}
	    \label[Problem]{eq:full_relaxed}
	\end{subequations}
 The force stress tensor $\widetilde{\bm{\sigma}}\coloneqq\Ce \sym (\D \vb{u} - \bm{P}) + \Cc \skw (\D \vb{u} - \bm{P})$ is symmetric if and only if $\Cc \equiv 0$, a case which is permitted. 
 \cref{eq:full_relaxed} represents a tensorial Maxwell-problem coupled to linear elasticity.
    We observe that the Dirichlet boundary condition for the microdistortion controls only its tangential components. 
    It is unclear, how to control the micro-movements of a material point without also affecting the displacement. Therefore, the relaxed micromorphic model introduces the so called consistent coupling condition \cite{dagostino2021consistent}
	\begin{align}
		\Pm \times \vb{n} = \D \widetilde{\ud} \times \vb{n} \quad \text{on} \quad \surf_D^P  \, , \label{eq:consistent_coupling}
	\end{align}
    where the prescribed displacement on the Dirichlet boundary $\widetilde{\ud}$ automatically dictates the tangential component of the microdistortion on that same boundary.
    Consequently, the consistent coupling condition enforces the definitions $\surf_D = \surf_D^u = \surf_D^P$ and $\surf_N = \surf_N^u = \surf_N^P$ (see \cref{fig:domain}).
    Further, the consistent coupling condition substitutes \cref{eq:pdir}.
    The set of equations in \cref{eq:full_relaxed} remains well-posed for $\Cc \equiv 0$ due to the generalized Korn inequality for incompatible tensor fields \cite{Lewintan2021K,LewintanInc,LewintanInc2,Neff2012}. The inequality relies on a non-vanishing Dirichlet boundary for the microdistortion field $\surf_D^P \neq \emptyset$, which the consistent coupling condition guarantees.
 
\subsection{Limits of the characteristic length scale parameter - a true two scale model} \label{sub:limits}
	In the relaxed micromorphic model the characteristic length $\Lc$ takes the role of a scaling parameter between the well-defined macro and the micro scales. 
    This property, unique to the relaxed micromorphic model, allows the theory to interpolate between materials with a pronounced micro-structure and homogeneous materials, thus relating the characteristic length scale parameter $\Lc$ to the size of the micro-structure in metamaterials. 
	In the lower limit $\Lc \to 0$ the continuum is treated as homogeneous and the solution of the classical Cauchy continuum theory is retrieved \cite{Neff2019,Aivaliotis2020}.
 This can be observed by reconsidering \cref{eq:strong_p} for $\Lc=0$,
	\begin{align}
		-\Ce  \sym (\D \vb{u} - \Pm) - \Cc  \skw(\D \vb{u} - \Pm) + \Cmic \sym \Pm = \bm{M} \, ,
	\end{align}
	which can now be used to express the microdistortion $\Pm$ algebraically
		\begin{align}
			\sym\Pm &= (\Ce + \Cmic)^{-1}(\sym \bm{M} + \Ce \sym \D \ud) \, , &&
			\skw\Pm = \Cc^{-1}\skw \bm{M} + \skw \D \ud \, .
		\end{align}
	Setting $\bm{M} = 0$ corresponds to Cauchy continua, where micro-moments are not accounted for. Thus, one finds
	\begin{align}
		&\Cc \skw(\D \ud - \Pm) = 0 \, , && \Ce \sym (\D \ud - \Pm) = \Cmic \sym \Pm \, , && \sym \Pm = (\Ce + \Cmic)^{-1}\Ce \sym \D \ud \, .
  \label{eq:lc_to_zero}
	\end{align}
	Applying the former results to \cref{eq:strong_u} yields
	\begin{align}
		-\Di[\Ce \sym (\D \vb{u} - \bm{P})]  = -\Di[\Cmic (\Ce + \Cmic)^{-1}\Ce \sym \D \ud ] = -\Di[\Cmac \sym \D \ud ] = \vb{f} \, , 
		\label{eq:cmacro}
	\end{align}
	where the definition 
         \begin{equation}
         \Cmac = \Cmic (\Ce + \Cmic)^{-1}\Ce    
         \end{equation}
         relates the meso- and micro-elasticity tensors to the classical macro-elasticity tensor of the Cauchy continuum. In fact, $\Cmac$ contains the material constants that arise from standard homogenization for large periodic structures \cite{Neff2019,Aivaliotis2020}.
	For isotropic materials one can directly express the macro parameters \cite{Neff2007forest}
	\begin{align}
		&\muma = \dfrac{\mue \, \mumi}{\mue + \mumi} \, ,
		&& 2 \muma + 3 \lamma = \dfrac{(2\mue + 3 \lame)(2 \mumi + 3 \lammi)}{(2\mue + 3 \lame) + (2 \mumi + 3 \lammi)}
		\label{eq:muma}
	\end{align}
	in terms of the parameters of the relaxed micromorphic model.
	
	In the upper limit $\Lc \to +\infty$, the stiffness of the micro-body becomes dominant. As the characteristic length $\Lc$ can be viewed as a zoom-factor into the microstructure, the state $\Lc \to + \infty$ can be interpreted as the entire domain being the micro-body itself. However, this is only theoretically possible as in practice, the limit is given by the size of one unit cell.
    Since the energy functional being minimized contains $\muma\Lc^2\| \Curl \Pm \|^2$, on contractible domains and bounded energy this implies the reduction of the microdistortion to a gradient field $\Pm \to \D \vb{v}$ due to the classical identity 
	\begin{align}
		\Curl \D \vb{v} = 0 \quad \forall \, \vb{v} \in [\C^\infty(\body)]^3 \, ,
	\end{align}
	thus asserting finite energies of the relaxed micromorphic model for arbitrarily large characteristic length values $\Lc$.
    The corresponding energy functional in terms of the reduced kinematics $\{ \ud, \vb{v} \} : V\to \mathbb{R}^{3}$ now reads
    \begin{align}
		I(\ud, \vb{v}) = \dfrac{1}{2} \int_\body &\langle \sym(\D \ud - \D \vb{v}) , \, \Ce \sym(\D \ud - \D \vb{v}) \rangle + \langle \sym\D \vb{v} , \, \Cmic \sym \D \vb{v} \rangle \notag \\
		& + \langle \skw(\D \ud - \D \vb{v}) , \, \Cc \skw(\D \ud - \D \vb{v}) \rangle \, \dd \body  - \int_\body \langle \ud , \, \vb{f} \rangle + \langle \D \vb{v} , \,  \bm{M} \rangle \, \dd \body \, , 
	\end{align}
    such that variation with respect to the two vector fields $\ud$ and $\vb{v}$ leads to
    \begin{subequations}
        \begin{align}
        \delta_u I &= \int_\body \langle \sym \D \delta \ud , \, \Ce \sym(\D \ud - \D \vb{v}) \rangle + \langle \skw \D \delta \ud , \, \Cc \skw(\D \ud - \D \vb{v}) \rangle - \langle \delta \ud , \, \vb{f} \rangle \, \dd \body = 0 \, , \label{eq:dv1} \\
        \delta_v I &= \int_\body \langle \sym \D \delta\vb{v} , \, \Ce \sym(\D \ud - \D \vb{v}) \rangle + \langle \skw \D \delta \vb{v} , \, \Cc \skw(\D \ud - \D \vb{v}) \rangle \notag \\ & \quad \quad - \langle \sym  \D \delta \vb{v} , \, \Cmic \sym \D \vb{v} \rangle + \langle \D \delta \vb{v} , \,  \bm{M} \rangle \, \dd \body = 0 \, . \label{eq:dv2}
    \end{align}
    \end{subequations}
    The resulting bilinear form is given by
    \begin{align}
		a(\{\delta \ud ,  \delta \vb{v}\},\{\ud,  \vb{v}\}) = \int_\body & \langle \sym(\D \delta \ud - \D \delta \vb{v}) , \, \Ce \sym(\D \ud - \D \vb{v}) \rangle + \langle \sym \D \delta \vb{v}, \, \Cmic \sym \D \vb{v} \rangle \notag \\
		& + \langle \skw(\D \delta \ud - \D \delta \vb{v}) , \, \Cc \skw(\D \ud - \D \vb{v}) \rangle \, \dd \body \, .
        \label{eq:new_bi}
	\end{align}
    By partial integration of \cref{eq:dv1} and \cref{eq:dv2} one finds the equilibrium equations
    \begin{subequations}
        \begin{align}
            -\Di[\Ce \sym (\D \vb{u} - \D \vb{v}) + \Cc \skw (\D \vb{u} - \D \vb{v})] &= \vb{f} && \text{in} \quad \body \, , 
            \label{eq:lc_strong_u}
            \\
			-\Di[\Ce  \sym (\D \vb{u} - \D \vb{v}) + \Cc  \skw(\D \vb{u} - \D \vb{v})] + \Di[\Cmic \sym \D \vb{v}]  &= \Di \bm{M} && \text{in} \quad \body \, .
            \label{eq:lc_strong_p}
        \end{align}
        \label{eq:lc_eq}
    \end{subequations}
    We can now substitute the right-hand side of \cref{eq:lc_strong_u} into \cref{eq:lc_strong_p} to find
    \begin{align}
		-\Di (\Cmic \sym \D \vb{v}) = \vb{f} - \Di \bm{M}  \, .
		\label{eq:upper}
	\end{align}
    Clearly, setting $\vb{v} = \ud$ satisfies both local equilibrium equations \cref{eq:lc_strong_u} and \cref{eq:lc_strong_p} for $\vb{f} = 0$. Further, the consistent coupling condition \cref{eq:consistent_coupling} is also automatically satisfied, asserting the equivalence of the tangential projections of both fields on the boundary of the domain. Since, as shown in \cite{SKY2022115298,Neff2019} using the extended Brezzi theorem, the case $\Lc \to +\infty$ is well-posed (including $\Cc \equiv 0$), the solution $\vb{v} = \ud$ is the unique solution to the bilinear form \cref{eq:new_bi} with the right-hand side 
    \begin{align}
        l(\{\delta \ud , \delta \vb{v}\}) =  \langle \D \delta \vb{v} , \,  \bm{M} \rangle \, \dd \body \, .
    \end{align}
    Effectively, equation~\cref{eq:upper} implies that the limit $\Lc \to +\infty$ defines a classical Cauchy continuum with a finite stiffness governed by $\Cmic$, representing the upper limit of the stiffness for the relaxed micromorphic continuum \cite{Neff2019}, where the corresponding forces read $\vb{m} = \Di \bm{M}$.
    We emphasize that this interpretation  of $\Cmic$ is impossible in the classical micromorphic model since there the limit $\Lc \to +\infty$ results in a constant microdistortion field $\Pm:\body \to \R^{3 \times 3}$ as its full gradient $\D \Pm$ is incorporated via $\muma\Lc^2 \| \D \Pm \|^2$ into the energy functional \cite{Barbagallo2017}. 
 
    %

\subsection{Antiplane shear}
\begin{figure}[H]
        \centering
        \begin{tikzpicture}[line cap=round,line join=round,>=triangle 45,x=1.0cm,y=1.0cm]
				\clip(-8,-9) rectangle (8,4.5);
				\draw (0,0) node[] {relaxed micromorphic};
                \draw [-Triangle,line width=.5pt] (0,-1) -- (0,-3);
				\draw [-Triangle,line width=.5pt] (2,0) -- (4,0);
                \draw [-Triangle,line width=.5pt] (-2,0) -- (-4,0);
                
                \draw (6,0) node[] {Cosserat elasticity};
                \draw [-Triangle,line width=.5pt] (6,-1) -- (6,-3);

                \draw (-6,0) node[] {\begin{tabular}{c}
                     linear elasticity  \\
                     with $\Cmac$
                \end{tabular}};
                \draw [-Triangle,line width=.5pt] (-6,-1) -- (-6,-3);

                \draw (0,-4) node[] {\begin{tabular}{c}
                     antiplane relaxed   \\
                     micromorphic
                \end{tabular}};
                \draw [-Triangle,line width=.5pt] (2,-4) -- (4,-4);
                \draw [line width=.5pt] (3.25-0.075,-3.85) -- (2.75-0.075,-4.15);
                \draw [line width=.5pt] (3.25+0.075,-3.85) -- (2.75+0.075,-4.15);
                \draw [-Triangle,line width=.5pt] (-2,-4) -- (-4,-4);

                \draw (6,-4) node[] {\begin{tabular}{c}
                     antiplane Cosserat  \\
                     elasticity 
                \end{tabular}};

                \draw (-6,-4) node[] {\begin{tabular}{c}
                     antiplane linear   \\
                     elasticity \\ 
                     with $\muma$
                \end{tabular}};

                \draw (-3,0) node[anchor=south] {$\Lc \to 0$};

                \draw (3,0) node[anchor=south] {\begin{tabular}{c}
                     $\Cmic \to + \infty \, ,$  \\
                    $\muc > 0$
                \end{tabular}};
                
                \draw (-3,-4) node[anchor=south] {\begin{tabular}{c}
                     $\Lc \to 0 \, ,$  \\
                    $\muc \equiv 0$
                \end{tabular}};

                \draw (3,-4) node[anchor=south] {\begin{tabular}{c}
                     $\mumi \to + \infty \, ,$  \\
                    $\muc > 0$ 
                \end{tabular}};

                \draw (3,-4.05) node[anchor=north] {(contradiction)};

                \draw (-6,-2) node[anchor=west] {\begin{tabular}{c}
                     antiplane \\
                     shear
                \end{tabular}};
                \draw (0,-2) node[anchor=west] {\begin{tabular}{c}
                     antiplane \\
                     shear
                \end{tabular}};
                \draw (6,-2) node[anchor=west] {\begin{tabular}{c}
                     antiplane \\
                     shear
                \end{tabular}};

                \draw [-Triangle,line width=.5pt] (0,1) -- (0,3);

                \draw [-Triangle,line width=.5pt] (0,-5) -- (0,-7);

                \draw (0,-8) node[] {\begin{tabular}{c}
                     antiplane linear   \\
                     elasticity \\ 
                     with $\mumi$
                \end{tabular}};

                \draw (0,4) node[] {\begin{tabular}{c}
                     linear elasticity  \\
                     with $\Cmic$
                \end{tabular}};

                \draw (0,2) node[anchor=east] {$\Lc \to +\infty$};

                \draw (0,-6) node[anchor=west] {$\Lc \to +\infty$};

                \draw [Triangle-,line width=.5pt] (-6,1-6) -- (-5,1.75-7.5);
                \draw [-Triangle,line width=.5pt] (-3,3.25-10.5) -- (-2,4-12);

                \draw [Triangle-,line width=.5pt] (-6,1) -- (-5,1.75);
                \draw [-Triangle,line width=.5pt] (-3,3.25) -- (-2,4);

                \node at (-4,2.5) [rectangle,draw] {\begin{tabular}{c}
                     two-scale  \\
                     model
                \end{tabular}};

                \node at (-4,-6.5) [rectangle,draw] {\begin{tabular}{c}
                     two-scale  \\
                     model
                \end{tabular}};

                \node at (3.75,-2) [rectangle,draw] {\begin{tabular}{c}
                     non-  \\
                     commutative
                \end{tabular}};
			\end{tikzpicture}
        \caption{Kinematic reduction of the relaxed micromorphic model to antiplane shear and consistency at limit cases according to \cref{re:muc} and \cref{sub:limits}. The two-scale nature of the relaxed micromorphic model can be clearly observed.}
        \label{fig:antiplane_flow}
    \end{figure}
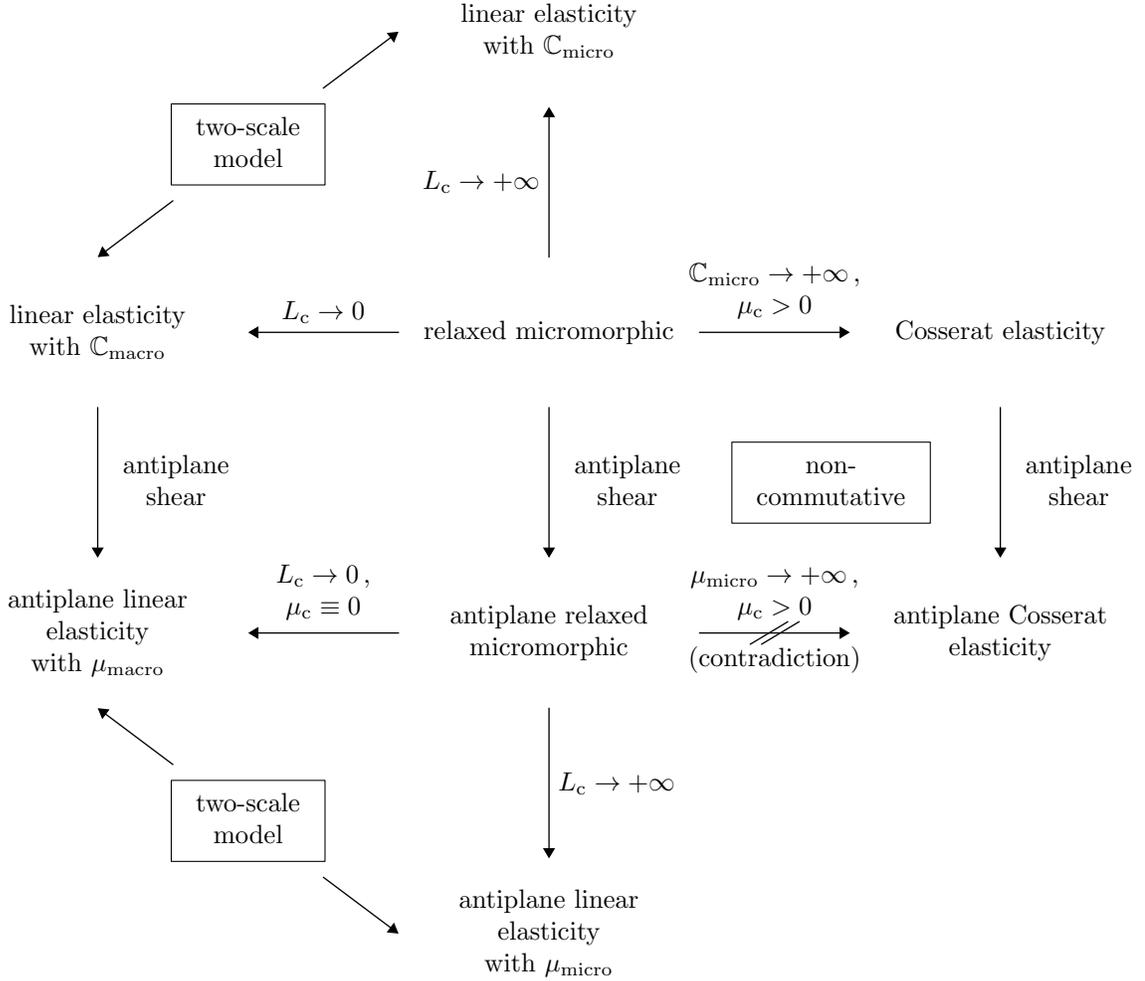

We introduce the relaxed micromorphic model of antiplane shear\footnote{Note that the antiplane shear model encompasses $1+2 = 3$ degrees of freedom and is the simplest non-trivial active version for the relaxed micromorphic model, as the one-dimensional elongation ansatz features only $1+1=2$ degrees of freedom and eliminates the curl operator
$$I(u,p) = \dfrac{1}{2} \int_\curv (\lame + 2 \mue) |u' - p|^2 + (\lammi + 2\mumi)|p|^2 \, \dd \curv - \int_\curv  u  \, f  +  p \, m  \, \dd \curv \to  \min \quad \wrt \quad \{u, p\} \, ,$$
since $\D \ud = u' \, \vb{e}_1 \otimes \vb{e}_1$ and $\Pm = p \, \vb{e}_1 \otimes \vb{e}_1$, such that $\skw (\D \ud - \Pm) = 0$ and $\Curl \Pm = 0$. 
This is not to be confused with uniaxial extension, which entails $1+3 = 4$ degrees of freedom \cite{Rizzi_extension}. 
} \cite{Voss2020} by reducing the displacement field to
	\begin{align}
		\ud = \begin{bmatrix}
			0, & 0, & u
		\end{bmatrix}^T \, ,
	\end{align}
    such that $\ud = \ud(x,y)$ is a function of the $x-y$-plane. 
    Consequently, its gradient reads
    \begin{align}
    	\D \ud = \begin{bmatrix}
    		0 & 0 & 0 \\
    		0 & 0 & 0 \\
    		u_{,x} & u_{,y} & 0 
    	\end{bmatrix} \, .
    \end{align}
    The structure of the microdistortion tensor is chosen accordingly 
    \begin{align}
    	\Pm = \begin{bmatrix}
    		0 & 0 & 0 \\
    		0 & 0 & 0 \\
    		p_{1} & p_{2} & 0 
    	\end{bmatrix} \, ,
     &&
    	\Curl \Pm = \begin{bmatrix}
    		0 & 0 & 0 \\
    		0 & 0 & 0 \\
    		  0 & 0 & p_{2,x}-p_{1,y} 
    	\end{bmatrix}  
     = \begin{bmatrix}
    		0 & 0 & 0 \\
    		0 & 0 & 0 \\
    		  0 & 0 & \rot{\vb{p}}
    	\end{bmatrix}\, .
    \end{align}
    Analogously to the displacement field $\ud$, the microdistortion $\Pm$ is also set to be a function of the $\{x,y\}$-variables $\Pm = \Pm(x,y)$.
    We observe the following sym-skew decompositions of the gradient and microdistortion tensors
    \begin{align}
    	 \sym \Pm &= \dfrac{1}{2} \begin{bmatrix}
    	 	0 & 0 & p_{1} \\
    	 	0 & 0 & p_{2} \\
    	 	p_{1} & p_{2} & 0 
    	 \end{bmatrix} \, , \qquad
    	 \sym(\D \ud - \Pm) = \dfrac{1}{2} \begin{bmatrix}
    	 	0 & 0 & u_{,x}-p_1 \\
    	 	0 & 0 & u_{,y}-p_2 \\
    	 	u_{,x}-p_1 & u_{,y}-p_2 & 0 
    	 \end{bmatrix} \, ,   \notag \\
    	 \skw(\D \ud - \Pm) &= \dfrac{1}{2} \begin{bmatrix}
    	 	0 & 0 & p_1-u_{,x} \\
    	 	0 & 0 & p_2-u_{,y} \\
    	 	u_{,x}-p_1 & u_{,y}-p_2 & 0 
    	 \end{bmatrix} \, .
    \end{align}
    Clearly, there holds
    \begin{align}
    	\tr [\sym \Pm] = \tr [\sym(\D \ud - \Pm)] = \tr [\skw(\D \ud - \Pm)] = 0 \, ,
    \end{align}
    such that the contraction with the material tensors reduces to
    \begin{align}
    	\Ce \sym(\D \ud - \Pm) &= 2 \mue \sym(\D \ud - \Pm) \, , & 
    	\Cmic \sym(\D \ud - \Pm) = 2 \mumi \sym\Pm \, , \notag \\ \Cc \skw(\D \ud - \Pm) &= 2 \muc \skw(\D \ud - \Pm) \, .
    \end{align}
	As such, the quadratic forms of the energy functional are given by
	\begin{subequations}
		\begin{align}
			\langle \sym(\D \ud - \Pm) , \, \Ce \sym(\D \ud - \Pm) \rangle &= \mue \| \nabla u - \vb{p} \|^2 \, , \\
			\langle \skw(\D \ud - \Pm) , \, \Cc \skw(\D \ud - \Pm) \rangle &= \muc \| \nabla u - \vb{p} \|^2 \, , \\
			\langle \sym\Pm , \, \Cmic \sym \Pm \rangle &= \mumi \| \vb{p} \|^2 \, ,
		\end{align}
	\end{subequations}
    with the definitions
    \begin{align}
    	&\nabla u = \begin{bmatrix} 
    		u_{,x} \\
    		u_{,y} 
    	\end{bmatrix} \, , &&
    	\vb{p} = \begin{bmatrix} 
    		p_1 \\
    		p_2 
    	\end{bmatrix} \, .
    \end{align}
    The resulting energy functional for antiplane shear reads therefore
    \begin{align}
    	I(u, \vb{p}) = \dfrac{1}{2} \int_\surf (\mue + \muc) \| \nabla u - \vb{p} \|^2  + \mumi \| \vb{p} \|^2 + \muma \Lc^2 \| \rot{\vb{p}} \|^2 \, \dd \surf - \int_\surf u \, f  + \langle \vb{p} , \, \vb{m} \rangle \, \dd \surf \, .
    \end{align} 
    In order to maintain consistency with the three-dimensional model we must choose $\muc = 0$. The reasoning for this choice is explained upon in \cref{re:muc} (see also \cref{fig:antiplane_flow}). Consequently, the energy functional is given by
	\begin{align}
		I(u, \vb{p}) = \dfrac{1}{2} \int_\surf& \mue \| \nabla u - \vb{p} \|^2 + \mumi \| \vb{p} \|^2 +  \muma \Lc^2 \| \rot{\vb{p}} \|^2 \, \dd \surf \notag \\ & - \int_\surf u\, f \,  + \langle \vb{p} , \, \vb{m} \rangle \, \dd \surf \to  \min \quad \wrt \quad \{ u, \vb{p} \} \, .
	\end{align}
    Note that on two-dimensional domains the differential operators are reduced to
\begin{align}
    &\nabla u = \begin{bmatrix} u_{,x} \\ u_{,y} \end{bmatrix} \, ,
    && \rog u = \begin{bmatrix} u_{,y} \\ -u_{,x}  \end{bmatrix} \, , 
    && \bm{R} = \begin{bmatrix}
    0 & 1 \\
    -1 & 0
    \end{bmatrix} \, , &&  \rot{\vb{p}} = \di (\bm{R} \vb{p}) = p_{2,x} - p_{1,y} \, ,
\end{align}
where we note that $\curl_{\text{2D}}$ is just a rotated divergence.
	Taking variations of the energy functional with respect to the displacement field results in
	\begin{align}
		\delta_u I = \int_\surf \mue \langle \nabla \delta u , \, \nabla u - \vb{p} \rangle - \delta u \, f \, \dd \surf = 0 \, ,
		\label{eq:var_antiplane_disp}
	\end{align}
	and variation with respect to the microdistortion yields
	\begin{align}
		\delta_p I = \int_\surf \mue \langle \delta \vb{p} , \, \nabla u - \vb{p} \rangle - \mumi \langle \delta \vb{p} , \, \vb{p} \rangle - \muma \Lc^2 (\rot{\delta \vb{p}}) \rot{\vb{p}} + \langle \delta \vb{p} , \, \vb{m} \rangle \, \dd \surf = 0 \, .
		\label{eq:var_p_anti}
	\end{align}
	Consequently, one finds the bilinear and linear forms
	\begin{subequations}
		\begin{align}
			a(\{\delta u, \delta \vb{p}\},\{u,\vb{p}\}) &= \int_\surf \mue \langle \nabla \delta u - \delta \vb{p} , \, \nabla u - \vb{p} \rangle + \mumi \langle \delta \vb{p} , \, \vb{p} \rangle + \muma \Lc^2 (\rot{\delta \vb{p}}) \rot{\vb{p}} \, \dd \surf \, , \\
			l(\{\delta u, \delta \vb{p}\}) &= \int_\surf \delta u \, f + \langle \delta \vb{p} , \, \vb{m} \rangle \, \dd \surf \, .
		\end{align}
		\label{eq:bili_anti}
	\end{subequations}
	Partial integration of \cref{eq:var_antiplane_disp} results in
	\begin{align}
		\int_{\partial \surf} \delta u \, \langle \mue( \nabla u - \vb{p} ) , \, \vb{n} \rangle \, \dd \curv - \int_\surf \delta u \, [ \mue \di(\nabla u - \vb{p}) + f ] \, \dd \surf =  0 \, , 
	\end{align}
	and analogously for \cref{eq:var_p_anti}, yielding
	\begin{align}
		\int_\surf \langle \delta \vb{p} , \, \mue(\nabla u - \vb{p}) - \mumi  \, \vb{p} - \muma \Lc^2 \rog \rot{\vb{p}}  + \vb{m} \rangle \, \dd \surf  
		- \int_{\partial \surf} \langle \delta \vb{p} , \,  \muma \Lc^2 (\rot{\vb{p}}) \, \vb{t} \rangle \, \dd \curv = 0 \, .
	\end{align}
	Consequently, the strong form reads
	\begin{subequations}
		\begin{align}
			-\mue \di(\nabla u - \vb{p}) &= f  && \text{in} \quad \surf \, , \label{eq:strong_anti_p}  \\
			-\mue(\nabla u - \vb{p}) + \mumi \, \vb{p} + \muma \Lc^2 \rog \rot{\vb{p}} &= \vb{m}  && \text{in} \quad \surf \, , \label{eq:strong_anti_u} \\
			u &= \widetilde{u} && \text{on} \quad \curv_D^u \, , \\
			\langle \vb{p} , \, \vb{t} \rangle &= \langle \widetilde{\vb{p}} , \, \vb{t} \rangle && \text{on} \quad \curv_D^P \, , \\
			\langle \nabla u , \, \vb{n} \rangle &= \langle \vb{p} , \, \vb{n} \rangle && \text{on} \quad \curv_N^u \, ,\\
			\rot{\vb{p}} &= 0  && \text{on} \quad \curv_N^P \, .
		\end{align}
        \label[Problem]{eq:ap_strong}
	\end{subequations}
The consistent coupling condition accordingly reduces to
\begin{align}
    \langle \vb{p} , \, \vb{t} \rangle = \langle \nabla \widetilde{u} , \, \vb{t} \rangle \qquad &\qquad \text{on} \quad \curv_D = \curv_D^P = \curv_D^u  \, .
\end{align}
\begin{remark} \label{re:muc}
		Note that without setting $\muc = 0$ in the antiplane shear model, the analogous result to \cref{eq:cmacro} in the limit $\Lc \to 0$ would read
		\begin{align}
			- \underbrace{\left ( \dfrac{\mumi \, [\mue + \muc]}{\mue + \muc + \mumi} \right )}_{\neq \muma} \Delta u = f \, ,
		\end{align} 
		where the relation to the macro parameter $\muma$ in \cref{eq:muma} is lost. Further, the limit defined in \cref{eq:lc_to_zero} with $\bm{M} = 0$ yields the contradiction
	    \begin{align}
	    	\sym\Pm = (\Ce + \Cmic)^{-1} \Ce \sym \D \ud \, , &&
	    	\Cc\skw\Pm = \Cc\skw \D \ud \, ,
	    \end{align}
        since the equations degenerate to
        \begin{align}
        	&\vb{p} = \dfrac{\mue}{\mue + \mumi} \nabla u \, , &&  \muc\vb{p} = \muc\nabla u \, ,
        	\label{eq:anti_cond}
        \end{align}
        due to the equivalent three-dimensional forms for antiplane shear.
        Choosing $\mumi = 0$ leads to a loss of structure in the strong form \cref{eq:ap_strong}, while satisfying \cref{eq:anti_cond}.
        As such, we must set the Cosserat couple modulus $\muc = 0$ to preserve the structure of the equations and satisfy both \cref{eq:muma} and \cref{eq:anti_cond}.
        
        Although the relaxed micromorphic model includes the Cosserat model as a singular limit for $\Cmic \to +\infty$ ($\mumi \to +\infty$), it is impossible to deduce the Cosserat model of antiplane shear as a limit of the antiplane relaxed micromorphic model, since one needs to satisfy \cref{eq:anti_cond} for $\muc > 0$ and $\mumi \to +\infty$, which is impossible.
	\end{remark}
The kinematic reduction of the relaxed micromorphic model to antiplane shear and its behaviour in the limit cases of its material parameters is depicted in \cref{fig:antiplane_flow}.

\section{Polynomial basis}
In this section we briefly introduce Bernstein polynomials and dual numbers. Bernstein polynomials are used to construct both the $\Hone$-conforming subspace and, in conjunction with the polytopal template methodology, the N\'ed\'elec elements. The computation of derivatives of the Bernstein base functions is achieved by employing dual numbers, thus enabling the calculation of the value and the derivative of a base function simultaneously.  

\subsection{Bernstein polynomials}
	Bernstein polynomials of order $p$ are given by the binomial expansion of the barycentric representation of the unit line
	\begin{align}
		1 = (\lambda_1 + \lambda_2)^p = ((1-\xi) + \xi)^p = \sum_{i=0}^p \begin{pmatrix} p \\ i \end{pmatrix} \xi^i (1 - \xi)^{p-i} = \sum_{i=0}^p \dfrac{p!}{i!(p-i)!} \xi^i (1 - \xi)^{p-i} \, ,
	\end{align}
	where $\xi \in [0, 1]$. The Bernstein polynomial reads
	\begin{align}
		b_i^p(\xi) = \begin{pmatrix} p \\ i \end{pmatrix} \xi^i (1 - \xi)^{p-i} \, .
	\end{align}
	A direct result of the binomial expansion is that Bernstein polynomials form a partition of unity, see also \cref{fig:bernstein}
 \begin{figure}
	\centering
	\definecolor{asl}{rgb}{0.4980392156862745,0.,1.}
	\definecolor{asb}{rgb}{0.,0.4,0.6}
	\begin{tikzpicture}[scale=5,line cap=round,line join=round,>=triangle 45,x=1.0cm,y=1.0cm]
		\clip(-0.2,-0.2) rectangle (1.2,1.2);
		
		\draw [asb,domain=0:1, samples=50] plot ({\x}, {(1-\x)^4});
		\draw (0.07,0.75) node[color=asb,anchor=west] {$b_0^4(\xi)$};
		
		\draw [violet,domain=0:1, samples=50] plot ({\x}, {-4*\x*(\x - 1)^3});
		
		\draw (0.28,0.41) node[color=violet,anchor=south] {$b_1^4(\xi)$};
		
		\draw [cyan,domain=0:1, samples=50] plot ({\x}, {6*\x^2*(\x - 1)^2});
		
		\draw (0.5,0.37) node[color=cyan,anchor=south] {$b_2^4(\xi)$};
		
		\draw [asl,domain=0:1, samples=50] plot ({\x}, {4*\x^3*(1 - \x)});
		
		\draw (1-0.28,0.41) node[color=asl,anchor=south] {$b_3^4(\xi)$};
		
		\draw [blue,domain=0:1, samples=50] plot ({\x}, {(\x)^4});
		
		\draw (1-0.07,0.75) node[color=blue,anchor=east] {$b_4^4(\xi)$};
		
		\draw [-to,color=black] (-0.1,0) -- (1.1,0);
		
		\draw (1.1,0) node[color=black,anchor=west] {$\xi$};
		
		\draw [color=black, densely dashed] (0,-0.1) -- (0,1.1);
		
		\draw [color=black, densely dashed] (1,-0.1) -- (1,1.1);
		\draw (1.1,1) node[color=black,anchor=west] {$1$};
		\draw (-0.1,1) node[color=black,anchor=east] {$1$};
		
		\draw [color=black, dotted] (-0.1,1) -- (1.1,1);
		\draw (0,-0.1) node[color=black,anchor=north] {$0$};
		\draw [color=black, dotted] (-0.1,1) -- (1.1,1);
		\draw (1,-0.1) node[color=black,anchor=north] {$1$};
		
		\draw (0.5,-0.1) node[color=black,anchor=north] {$1/2$};
		\draw (0.5,1.1) node[color=black,anchor=south] {$1/2$};
		
		\draw [color=black, dotted] (0.5,-0.1) -- (0.5,1.1);
	\end{tikzpicture}
	\caption{Bernstein base functions of degree $p=4$ on the unit domain. Their sum forms a partition of unity. The base functions are symmetric for $\xi=0.5$ with respect to their indices and always positive.}
	\label{fig:bernstein}
\end{figure}
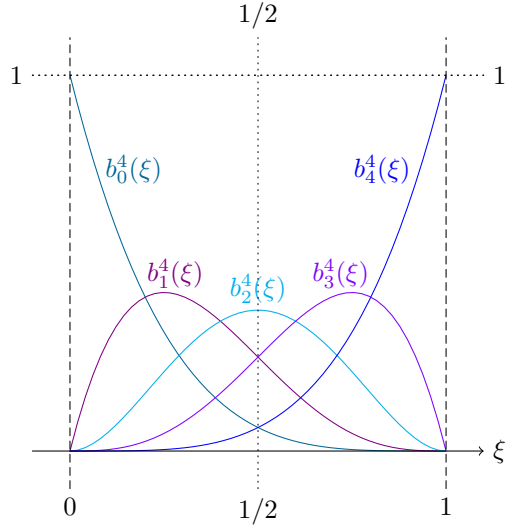
	\begin{align}
		\sum_{i = 0}^p b_i^p(\xi) = 1 \, . 
		\label{eq:unity}
	\end{align}
	Another consequence is that Bernstein polynomials are non-negative and less than or equal to 1
	\begin{align}
		&0\leq b_i^p(\xi) \leq 1 \, , && \xi \in [0,1] \, .
	\end{align}
	A necessary condition for the use of Bernstein polynomials in finite element approximations is for them to span the entire polynomial space.
	\begin{theorem} [Span of Bernstein polynomials]
		
		The span of Bernstein polynomials forms a basis of the one-dimensional polynomial space
		\begin{align}
			&\Po^p(\xi) = \spa \{ b_i^p \} \, , && \xi \subseteq \R \, .
		\end{align}
	\end{theorem}
	\begin{proof}
		First we observe 
		\begin{align}
			\dim (\spa \{ b_i^p \}) = \dim \Po^p(\xi) = p + 1 \, .
		\end{align}
		The proof of linear independence is achieved by contradiction. Let the set $\spa \{ b_i^{p} \}$ with $0<i\leq p$ be linearly dependent, then there exists some combination with at least one non-zero constant $c_i \neq 0$ such that 
		\begin{align}
			&\sum_{i=1}^{p} c_i b_i^p(\xi) = 0 \, , && \dfrac{\dd}{\dd \xi} \sum_{i=1}^{p} c_i b_i^p(\xi) = 0 \, .
		\end{align}
		However, by the partition of unity property \cref{eq:unity}, only the full combination ($0 \leq i \leq p$) generates a constant and by the exact sequence property the kernel of the differentiation operator is exactly the space of constants $\ker(\partial) = \R$. The linear independence of the full span also follows from the partition of unity property, since constants cannot be constructed otherwise.
	\end{proof}
	Bernstein polynomials can be evaluated efficiently using the recursive formula 
	\begin{align}
		& b_0^p(\xi) = (1-\xi)^p \, , &&b_{i+1}^p(\xi) = \dfrac{(p-i)\xi}{(p+1)(1-\xi)} b_i^p(\xi) \, , && i \in \{0,1,...,p-1\} \, , 
		\label{eq:rec}
	\end{align}
	which allows for fast evaluation of the base functions.
	\begin{remark}
		Note that the formula \cref{eq:rec} implies $\lim_{\xi \to 1} b_{i+1}^p(\xi) = \infty$. As such, evaluations using the formula are required to use $\xi < 1$ preferably with additional tolerance. The limit case $\xi = 1$ is zero for all Bernstein base functions aside from the last function belonging to the vertex, which simply returns one
		\begin{align}
		    b_i^p(1) = 0 \quad \forall \, i \neq p \, , && b_p^p(1) = 1 \, .
		\end{align}
	\end{remark}

\subsection{Dual numbers}
    Dual numbers \cite{Fike} can be used to define define an augmented algebra, where the derivative of a function can be computed simultaneously with the evaluation of the function.
	This enhancement is also commonly used in forward automatic differentiation \cite{Neidinger,baydin2018automatic}, not to be confused with numerical differentiation, since unlike in numerical differentiation, automatic differentiation is no approximation and yields the exact derivative.
    The latter represents an alternative method to finding the derivatives of base functions, as opposed to explicit formulas or approximations.
	Dual numbers augment the classical numbers by adding a non-zero number $\varepsilon$ with a zero square $\varepsilon^2 = 0$.
	\begin{definition} [Dual number]
		
		The dual number is defined by
		\begin{align}
			&x + x' \varepsilon \, , && \varepsilon \ll 1 \, ,
		\end{align}
		where $x'$ is the derivative (only in automatic differentiation), $\varepsilon$ is an abstract number (infinitesimal) and formally $\varepsilon^2 = 0$.
	\end{definition}
	The augmented algebra results automatically from the definition of the dual number.
	\begin{definition} [Augmented dual algebra]
		
		The standard algebraic operations take the following form for dual numbers
		\begin{enumerate}
			\item Addition and subtraction 
			\begin{align}
				(x + x'\varepsilon) \pm (y + y'\varepsilon) = x \pm y + (x' \pm y')\varepsilon \, .
			\end{align}
			\item Multiplication
			\begin{align}
				(x + x'\varepsilon)  (y + y'\varepsilon) = xy + (xy' + x'y)\varepsilon \, ,
			\end{align}
			since formally $\varepsilon^2 = 0$.
			\item Division is achieved by first defining the inverse element
			\begin{align}
				(x+ x' \varepsilon)(y + y'\varepsilon) = 1 \quad \iff \qquad y = \dfrac{1}{x}, \quad y' = -\dfrac{x'}{x^2} \, ,
			\end{align}
			such that
			\begin{align}
				(x + x'\varepsilon) / (y + y'\varepsilon) = x/y + (x'/y - xy'/y^2 )\varepsilon \, .
			\end{align}
		\end{enumerate}
	\end{definition}
	Application of the above definitions to polynomials
	\begin{align}
		p(x + \varepsilon) = \sum_{i = 0}^\infty c_i (x + \varepsilon)^i = \sum_{i=0}^\infty \sum_{j = 0}^1 c_i \begin{pmatrix}
		i \\ j
	\end{pmatrix} x^{i-j} \varepsilon^j = \sum_{i=0}^\infty  c_i  x^{i} + \varepsilon \sum_{i = 1}^\infty i \, c_i x^{i-1} = p(x) + p'(x) \varepsilon \, ,
	\end{align}
	allows the extension to various types of analytical functions with a power-series representation (such as trigonometric or hyperbolic).
	\begin{definition} [General dual numbers function] 
		
		A function of a dual number is defined in general by
		\begin{align}
			f(x+\varepsilon) = f(x) + f'(x) \varepsilon \, ,
		\end{align}
	    being a fundamental formula for forward automatic differentiation.
	\end{definition}
    The definition of dual numbers makes them directly applicable to the general rules of differentiation, such as the chain rule or product rule, in which case the derivative is simply the composition of previous computations with $\varepsilon$.
	The logic of dual numbers can be understood intuitively by the directional derivative 
	\begin{align}
		\dfrac{\dd }{ \dd x } f(x) = \partial_{x'}f(x) = \dfrac{\dd}{\dd \varepsilon} f(x + x'\varepsilon) \at_{\varepsilon = 0} = \lim_{\varepsilon \to 0} \dfrac{f(x + x'\varepsilon) - f(x)}{\varepsilon} \, ,
	\end{align}
	where dividing by $\varepsilon$ and setting $\varepsilon = 0$ are deferred to the last step of the computation, being the extraction of the derivative and equivalent to the operation $f(x + \varepsilon) - f(x)$ with the augmented algebra of dual numbers. 
	
In this work we apply dual numbers for the computation of Bernstein polynomials using the recursive formula \cref{eq:rec}, thus allowing to iteratively compute each base function simultaneously with its derivative.

\section{Triangular elements}
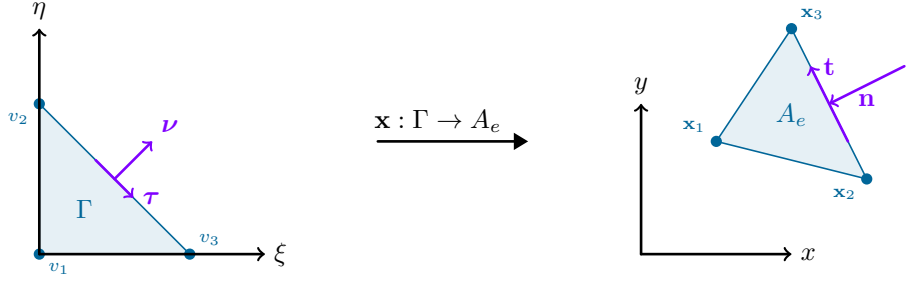
\begin{figure}
		\centering
		\definecolor{asl}{rgb}{0.4980392156862745,0.,1.}
		\definecolor{asb}{rgb}{0.,0.4,0.6}
		\begin{tikzpicture}[line cap=round,line join=round,>=triangle 45,x=1.0cm,y=1.0cm]
			\clip(-1,-0.5) rectangle (12,3.5);
			\draw (0,0) node[circle,fill=asb,inner sep=1.5pt] {};
			\draw (0,2) node[circle,fill=asb,inner sep=1.5pt] {};
			\draw (2,0) node[circle,fill=asb,inner sep=1.5pt] {};
			\draw [color=asb,line width=.6pt] (0,0) -- (0,2) -- (2,0) -- (0,0);
			\fill[opacity=0.1, asb] (0,0) -- (0,2) -- (2,0) -- cycle;
			\draw (0,0) node[color=asb,anchor=north west] {$_{v_1}$};
			\draw (2,0) node[color=asb,anchor=south west] {$_{v_3}$};
			\draw (0,2) node[color=asb,anchor=north east] {$_{v_2}$};
			\draw (0.6,0.6) node[color=asb] {$\Gamma$};
			
			\draw [-to,color=asl,line width=1.pt] (0.75,1.25) -- (1.25,0.75);
			\draw [-to,color=asl,line width=1.pt] (1,1) -- (1.5,1.5);
			\draw (1.5,1.5) node[color=asl,anchor=south west] {$\bm{\nu}$};
			\draw (1.25,0.75) node[color=asl,anchor=west] {$\bm{\tau}$};
			
			\draw [-to,color=black,line width=1.pt] (0,0) -- (3,0);
			\draw [-to,color=black,line width=1.pt] (0,0) -- (0,3);
			\draw (3,0) node[color=black,anchor=west] {$\xi$};
			\draw (0,3) node[color=black,anchor=south] {$\eta$};
			
			\draw (9,1.5) node[circle,fill=asb,inner sep=1.5pt] {};
			\draw (10,3) node[circle,fill=asb,inner sep=1.5pt] {};
			\draw (11,1) node[circle,fill=asb,inner sep=1.5pt] {};
			\draw (9,1.5) node[color=asb,anchor=south east] {$_{\vb{x}_1}$};
			\draw (10,3) node[color=asb,anchor=south west] {$_{\vb{x}_3}$};
			\draw (11,1) node[color=asb,anchor=north east] {$_{\vb{x}_2}$};
			\draw [color=asb,line width=.6pt] (9,1.5) -- (10,3) -- (11,1) -- (9,1.5);
			\fill[opacity=0.1, asb] (9,1.5) -- (10,3) -- (11,1) -- cycle;
			\draw (10,1.85) node[color=asb] {$\surf_e$};
			\draw [-to,color=asl,line width=1.pt] (10.75,1.5) -- (10.25,2.5);
			\draw (10.3,2.5) node[color=asl,anchor=west] {$\vb{t}$};
			\draw [to-,color=asl,line width=1.pt] (10.5,2) -- (11.5,2.5);
			\draw (11,2.25) node[color=asl,anchor=north] {$\vb{n}$};
			
			\draw [-to,color=black,line width=1.pt] (8,0) -- (10,0);
			\draw [-to,color=black,line width=1.pt] (8,0) -- (8,2);
			\draw (10,0) node[color=black,anchor=west] {$x$};
			\draw (8,2) node[color=black,anchor=south] {$y$};
			
			\draw (5.3,1.5) node[color=black,anchor=south] {$\vb{x}:\Gamma \to \surf_e$};
			\draw [-Triangle,color=black,line width=1.pt] (4.5,1.5) -- (6.5,1.5);
		\end{tikzpicture}
		\caption{Barycentric mapping of the reference triangle to an element in the physical domain.}
		\label{fig:trimap}
	\end{figure}

The triangle elements are mapped from the reference element $\Gamma$ to the physical domain $\surf_e$ via barycentric coordinates
	\begin{align}
		& \vb{x}(\xi, \eta) = (1 - \xi - \eta) \vb{x}_1 + \eta \, \vb{x}_2 + \xi \, \vb{x}_3 \, , && \vb{x}:\Gamma \to \surf_e \, ,  && \Gamma = \{ (\xi, \eta) \in [0,1]^2 \; | \; \xi + \eta \leq 1 \} \, ,
	\end{align}
	where $\vb{x}_i$ represent the coordinates of the vertices of one triangle in the physical domain, see \cref{fig:trimap}.
	The corresponding Jacobi matrix reads
	\begin{align}
		\bm{J} = \D \vb{x} = \begin{bmatrix}
		    \vb{x}_3 - \vb{x}_1, & \vb{x}_2 - \vb{x}_1
		\end{bmatrix}  \in \R^{2 \times 2} \, .
	\end{align}

\subsection{The Bernstein-B\'ezier basis for triangles}
\begin{figure}
    	\centering
    	\begin{subfigure}{0.3\linewidth}
    		\centering
    		\includegraphics[width=0.7\linewidth]{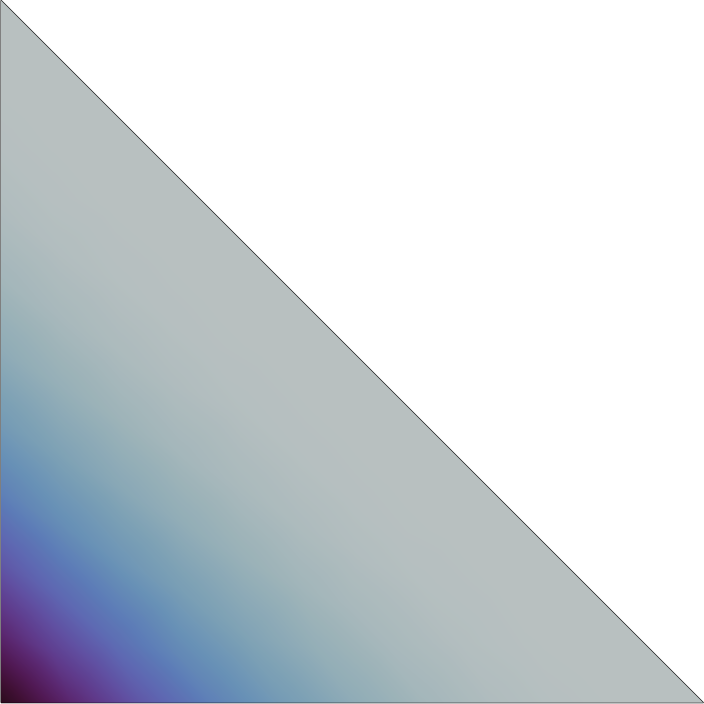}
    		\caption{}
    	\end{subfigure}
    	\begin{subfigure}{0.3\linewidth}
    		\centering
    		\includegraphics[width=0.7\linewidth]{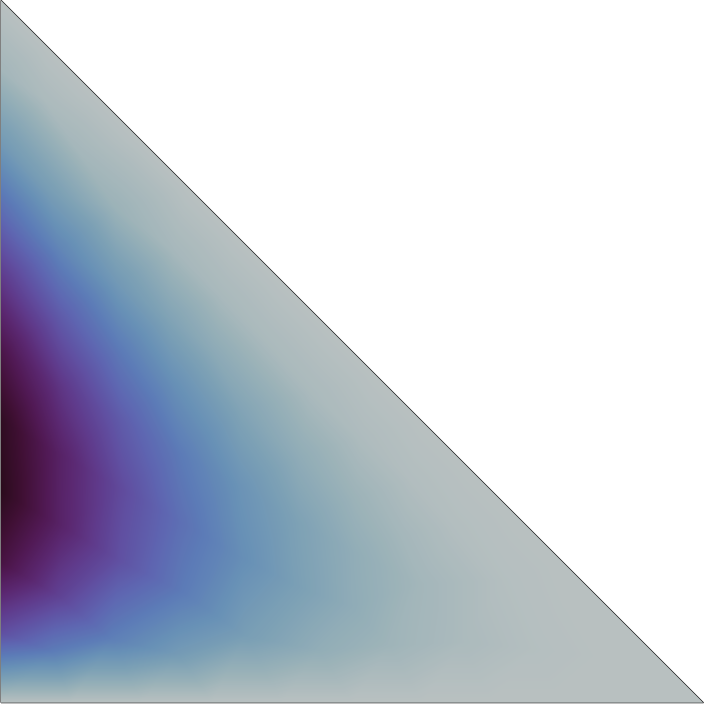}
    		\caption{}
    	\end{subfigure}
    	\begin{subfigure}{0.3\linewidth}
    		\centering
    		\includegraphics[width=0.7\linewidth]{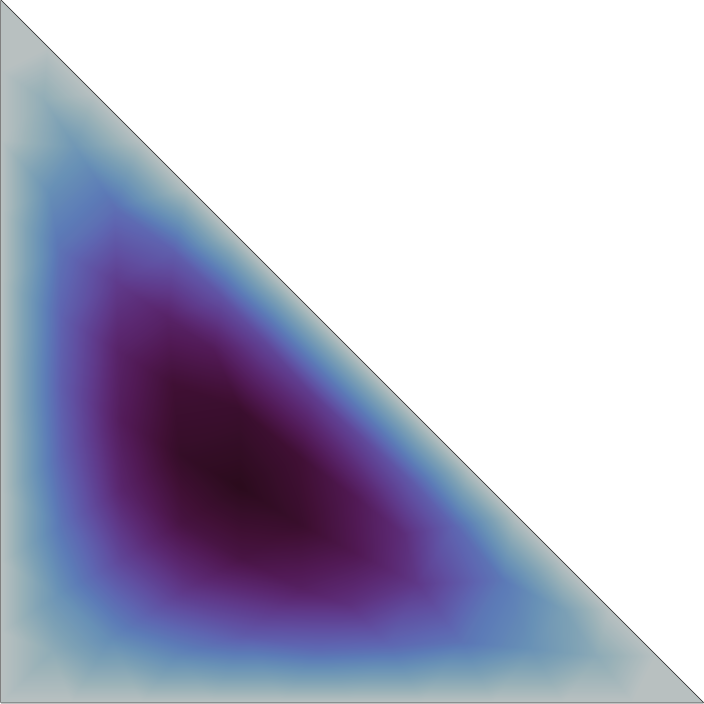}
    		\caption{}
    	\end{subfigure}
    	\caption{Cubic vertex (a), edge (b) and cell (c) B\'ezier base functions on the reference triangle.}
    	\label{fig:bezier}
    \end{figure}

\begin{figure}
		\centering
		\definecolor{asl}{rgb}{0.4980392156862745,0.,1.}
		\definecolor{asb}{rgb}{0.,0.4,0.6}
		\begin{tikzpicture}[line cap=round,line join=round,>=triangle 45,x=1.0cm,y=1.0cm]
			\clip(-0.5,-0.5) rectangle (11.5,3.5);
			\draw (0,0) node[circle,fill=asb,inner sep=1.5pt] {};
			\draw (0,2) node[circle,fill=asb,inner sep=1.5pt] {};
			\draw (2,0) node[circle,fill=asb,inner sep=1.5pt] {};
			\draw (2,2) node[circle,fill=asb,inner sep=1.5pt] {};
			\draw [color=asb,line width=.6pt] (0,0) -- (0,2) -- (2,2) -- (2,0) -- (0,0);
			\fill[opacity=0.1, asb] (0,0) -- (0,2) -- (2,2) -- (2,0) -- cycle;
			\draw (0,0) node[color=asb,anchor=north west] {$_{(0,0)}$};
			\draw (2,0) node[color=asb,anchor=north west] {$_{(1,0)}$};
			\draw (2,2) node[color=asb,anchor=south west] {$_{(1,1)}$};
			\draw (0,2) node[color=asb,anchor=south west] {$_{(0,1)}$};
			
			\draw [-to,color=black,line width=1.pt] (0,0) -- (3,0);
			\draw [-to,color=black,line width=1.pt] (0,0) -- (0,3);
			\draw (3,0) node[color=black,anchor=west] {$\alpha$};
			\draw (0,3) node[color=black,anchor=south] {$\beta$};
			
			\draw (8,0) node[circle,fill=asb,inner sep=1.5pt] {};
			\draw (8,2) node[circle,fill=asb,inner sep=1.5pt] {};
			\draw (10,0) node[circle,fill=asb,inner sep=1.5pt] {};
			\draw [color=asb,line width=.6pt] (8,0) -- (8,2) -- (10,0) -- (8,0);
			\fill[opacity=0.1, asb] (8,0) -- (8,2) -- (10,0) -- cycle;
			\draw (8.6,0.6) node[color=asb] {$\Gamma$};
			\draw (8,0) node[color=asb,anchor=north west] {$_{(0,0)}$};
			\draw (10,0) node[color=asb,anchor=north west] {$_{(1,0)}$};
			\draw (8,2) node[color=asb,anchor=south west] {$_{(0,1)}$};
			
			\draw [-to,color=black,line width=1.pt] (8,0) -- (11,0);
			\draw [-to,color=black,line width=1.pt] (8,0) -- (8,3);
			\draw (11,0) node[color=black,anchor=west] {$\xi$};
			\draw (8,3) node[color=black,anchor=south] {$\eta$};
			
			\draw (5.3,1.5) node[color=black,anchor=south] {$\bm{\xi}:\bm{\alpha}\to\Gamma$};
			\draw [-Triangle,color=black,line width=1.pt] (4.5,1.5) -- (6.5,1.5);
		\end{tikzpicture}
		\caption{Duffy transformation from a quadrilateral to a triangle by collapse of the coordinate system.}
		\label{fig:duffy2d}
	\end{figure}
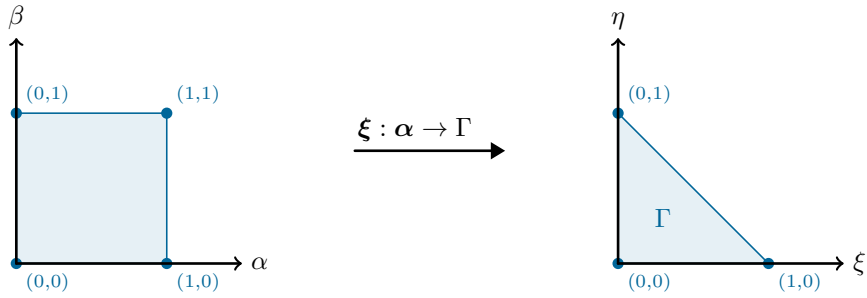

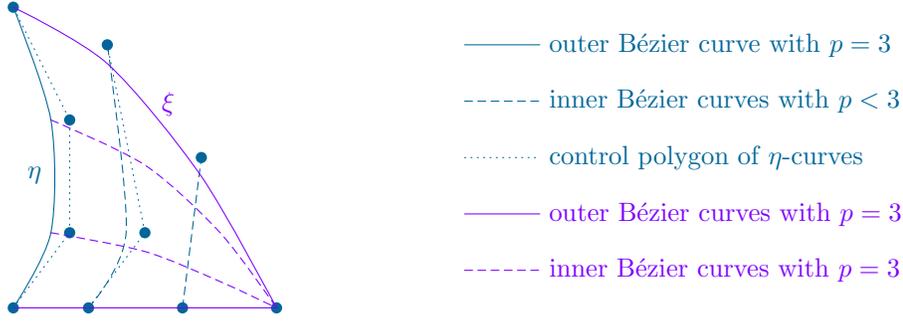
\begin{figure}
		\centering
		\definecolor{asl}{rgb}{0.4980392156862745,0.,1.}
		\definecolor{asb}{rgb}{0.,0.4,0.6}
		\begin{tikzpicture}[line cap=round,line join=round,>=triangle 45,x=1.0cm,y=1.0cm]
			\clip(-0.5,-0.5) rectangle (14,5);
			\draw (0,0) node[circle,fill=asb,inner sep=1.5pt] {};
			\draw (0.75,1) node[circle,fill=asb,inner sep=1.5pt] {};
			\draw (0.75,2.5) node[circle,fill=asb,inner sep=1.5pt] {};
			\draw (0,4) node[circle,fill=asb,inner sep=1.5pt] {};
			\draw (1,0) node[circle,fill=asb,inner sep=1.5pt] {};
			\draw (1.75,1) node[circle,fill=asb,inner sep=1.5pt] {};
			\draw (1.25,3.5) node[circle,fill=asb,inner sep=1.5pt] {};
			\draw (2.25,0) node[circle,fill=asb,inner sep=1.5pt] {};
			\draw (2.5,2) node[circle,fill=asb,inner sep=1.5pt] {};
			\draw (3.5,0) node[circle,fill=asb,inner sep=1.5pt] {};
			
			\draw [asb] plot [smooth] coordinates {(0,0) (0.5,1) (0.5,2.5) (0,4)};
			\draw [asb, densely dashed] plot [smooth] coordinates {(1,0) (1.5,1) (1.25,3.5)};
			\draw [color=asb, densely dashed] (2.25,0) -- (2.5,2);
			
			\draw [asl] plot [smooth] coordinates {(0,4) (1.25,3.25) (2.5,1.75) (3.5,0)};
			\draw [asl, densely dashed] plot [smooth] coordinates {(0.5,1) (1.75,0.75) (2.75,0.35) (3.5,0)};
			\draw [asl, densely dashed] plot [smooth] coordinates {(0.5,2.5) (1.75,1.9) (2.75,1) (3.5,0)};
			\draw [color=asl] (0,0) -- (3.5,0);
			
			\draw [color=asb,dotted] (0,0) -- (0.75,1) -- (0.75,2.5) -- (0,4);
			\draw [color=asb,dotted] (1,0) -- (1.75,1) -- (1.25,3.5);
			
			\draw (1.85,3) node[color=asl,anchor=north west] {$\xi$};
			\draw (0.5,2) node[color=asb,anchor=north east] {$\eta$};
			
			\draw [color=asb] (6,3.5) -- (7,3.5);
			\draw (7,3.5) node[color=asb,anchor=west] {outer B\'ezier curve with $p=3$};
			\draw [color=asb, densely dashed] (6,2.75) -- (7,2.75);
			\draw (7,2.75) node[color=asb,anchor=west] {inner B\'ezier curves with $p<3$};
			\draw [color=asb, dotted] (6,2) -- (7,2);
			\draw (7,2) node[color=asb,anchor=west] {control polygon of $\eta$-curves};
			\draw [color=asl] (6,1.25) -- (7,1.25);
			\draw (7,1.25) node[color=asl,anchor=west] {outer B\'ezier curves with $p=3$};
			\draw [color=asl, densely dashed] (6,0.5) -- (7,0.5);
			\draw (7,0.5) node[color=asl,anchor=west] {inner B\'ezier curves with $p=3$};
		\end{tikzpicture}
		\caption{B\'ezier triangle built by interpolating B\'ezier curves with an ever decreasing polynomial degree.}
		\label{fig:inter}
	\end{figure}

The base functions on the triangle reference element are defined using the binomial expansion of the barycentric coordinates on the domain $\Gamma$
	\begin{align}
		1 = (\lambda_1 + \lambda_2 + \lambda_3)^p = ([1- \xi - \eta] + \eta + \xi)^p \, .
	\end{align}
	As such, the B\'ezier base functions read 
	\begin{align}
		b_{ij}^p(\lambda_1,\lambda_2,\lambda_3) = \begin{pmatrix} p \\ i \end{pmatrix}
		\begin{pmatrix} p-i \\ j \end{pmatrix} \lambda_1^{p-i-j} \lambda_2^j \lambda_3^i \, ,
	\end{align}
	with the equivalent bivariate form
	\begin{align}
		b^p_{ij}(\xi,\eta) = \begin{pmatrix} p \\ i \end{pmatrix}
		\begin{pmatrix} p-i \\ j \end{pmatrix} 
		(1- \xi - \eta)^{p-i-j} \eta^j \xi^i \, ,
	\end{align}
    of which some examples are depicted in \cref{fig:bezier}.
	The Duffy transformation 
	\begin{align}
		&\bm{\xi}:[0,1]^2  \to \Gamma \, , && \{\alpha, \beta\} \mapsto  \{\xi, \eta\} \, , 
	\end{align}
	given by the relations
	\begin{align}
		\xi = \alpha \, , &&
		\alpha = \xi \, , && \eta = (1- \alpha) \beta \, ,  &&\beta = \dfrac{\eta}{1-\xi} \, , 
	\end{align}
	allows to view the triangle as a collapsed quadrilateral, see \cref{fig:duffy2d}. 
	Inserting the Duffy map into the definition of the B\'ezier base function yields the split 
	\begin{align}
		b^p_{ij}(\xi,\eta) &= \begin{pmatrix} p \\ i \end{pmatrix}
		\begin{pmatrix} p-i \\ j \end{pmatrix} 
		(1- \xi - \eta)^{p-i-j} \eta^j \xi^i  \notag \\
		&= \begin{pmatrix} p \\ i \end{pmatrix}
		\begin{pmatrix} p-i \\ j \end{pmatrix} (1-\alpha - [1-\alpha]\beta)^{p-i-j}(1-\alpha)^j \beta^j \alpha^i \notag\\
		&= \begin{pmatrix} p \\ i \end{pmatrix}
		\begin{pmatrix} p-i \\ j \end{pmatrix} (1-\alpha)^{p-i-j}(1-\beta)^{p-i-j}(1-\alpha)^j\beta^j\alpha^i  \\
		&= \begin{pmatrix} p \\ i \end{pmatrix} (1-\alpha)^{p-i} \alpha^i
		\begin{pmatrix} p-i \\ j \end{pmatrix} (1-\beta)^{p-i-j} \beta^j \notag \\
		&= b_i^p(\alpha) \, b_j^{p-i}(\beta) \, .
        \notag 
	\end{align}
	In other words, the Duffy transformation results in a natural factorization of the B\'ezier triangle into Bernstein base functions \cite{AinsworthOpt}. The latter allows for fast evaluation using sum factorization. Further, it is now clear that B\'ezier triangles are given by the interpolation of B\'ezier curves, where the degree of the polynomial decreases between each curve, see \cref{fig:inter}.   
	In order to compute gradients on the reference domain one applies the chain rule 
	\begin{align}
		&\nabla_\xi b_{ij}^p = (\D_\alpha \bm{\xi})^{-T} \nabla_\alpha b_{ij}^p \, , && \D_\alpha \bm{\xi} = \begin{bmatrix}
			1 & 0 \\
			-\beta & 1 - \alpha
		\end{bmatrix} \, , && (\D_\alpha \bm{\xi})^{-T} = \dfrac{1}{1-\alpha} \begin{bmatrix}
			1- \alpha & \beta \\
			0 & 1
		\end{bmatrix} \, .
	\end{align}
    The factorization is naturally suited for the use of dual numbers since the $\bm{\alpha}$-gradient of a base function reads
    \begin{align}
    	\nabla_\alpha b_{ij}^p(\alpha, \beta) = \begin{bmatrix}
    		b_j^{p-i} \dfrac{\dd }{\dd \alpha} b_i^p \\[2ex]
    		b_i^p \dfrac{\dd }{\dd \beta} b_j^{p-i}
    	\end{bmatrix} \, ,
    \end{align} 
    such that only the derivatives of the Bernstein base functions with respect to their parameter are required.
	
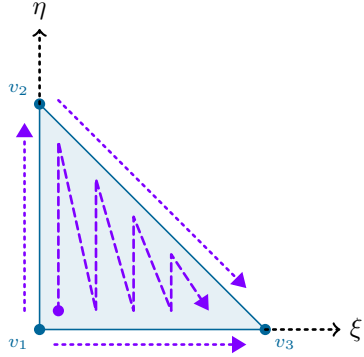
\begin{figure}
		\centering
		\definecolor{asl}{rgb}{0.4980392156862745,0.,1.}
		\definecolor{asb}{rgb}{0.,0.4,0.6}
		\begin{tikzpicture}[line cap=round,line join=round,>=triangle 45,x=1.0cm,y=1.0cm]
			\clip(-0.5,-0.5) rectangle (4.5,4.5);
			\draw (0,0) node[circle,fill=asb,inner sep=1.5pt] {};
			\draw (0,3) node[circle,fill=asb,inner sep=1.5pt] {};
			\draw (3,0) node[circle,fill=asb,inner sep=1.5pt] {};
			\draw [color=asb,line width=.6pt] (0,0) -- (0,3) -- (3,0) -- (0,0);
			\fill[opacity=0.1, asb] (0,0) -- (0,3) -- (3,0) -- cycle;
			\draw (0,0) node[color=asb,anchor=north east] {$_{v_1}$};
			\draw (3,0) node[color=asb,anchor=north west] {$_{v_3}$};
			\draw (0,3) node[color=asb,anchor=south east] {$_{v_2}$};
			
			\draw [-to,color=black,line width=1.pt, dotted] (3,0) -- (4,0);
			\draw [-to,color=black,line width=1.pt, dotted] (0,3) -- (0,4);
			\draw (4,0) node[color=black,anchor=west] {$\xi$};
			\draw (0,4) node[color=black,anchor=south] {$\eta$};
			
			\draw (0.25,0.25) node[circle,fill=asl,inner sep=1.5pt] {};
			\draw [-Triangle,color=asl,line width=1.pt, densely dashed] (0.25,0.25) -- (0.25,2.5) -- (0.75,0.25) -- (0.75,2) -- (1.25,0.25) -- (1.25,1.5) -- (1.75,0.25) -- (1.75,1.) -- (2.25,0.25);
			
			\draw [-Triangle,color=asl,line width=1.pt, dotted] (0.25,3.05) -- (2.75,0.55);
			\draw [-Triangle,color=asl,line width=1.pt, dotted] (-0.2,0.25) -- (-0.2,2.75);
			\draw [-Triangle,color=asl,line width=1.pt, dotted] (0.2,-0.2) -- (2.75,-0.2);
		\end{tikzpicture}
		\caption{Traversal order of base functions. The purple lines represent the order in which the base functions are constructed by the factorized evaluation. Note that the traversal order on each edge is intrinsically from the lower to the higher vertex index.}
		\label{fig:traversaltri}
	\end{figure}

	The Duffy transformation induces an intrinsic optimal order of traversal of the base functions, compare \cref{fig:traversaltri}, namely 
	\begin{align}
		(i,j)=(0,0) \to (0,1) \to...\to (2,2)\to...\to(i,p-i)\to...\to(p,0) \, ,
	\end{align}
	which respects a clockwise orientation of the element, compare \cite{SKY2022115298}. Thus, the order of the sequence of discrete values on common edges is determined by the global orientation.
	In order to relate a base function to a polytopal piece of the element, one observes the following result.
	\begin{observation} [Triangle base functions]
		
		The polytope of each base function $b_{ij}^p(\xi,\eta)$ can be determined as follows:
		\begin{itemize}
			\item The indices $(0,0)$, $(0,p)$ and $(p,0)$ represent the first, second and last vertex base functions, respectively.
			\item The indices $(0, j)$ with $0< j < p$ and $(i, 0)$ with $0<i<p$ represent the first and second edge base functions, respectively. Base functions of the slanted edge are given by $( i, p-i )$ with $0<i<p$.
			\item The remaining index combinations are cell base functions.
		\end{itemize}
	\end{observation}
	With the latter observation, the construction of vertex-, edge- and cell base functions follows the intrinsic traversal order induced by the Duffy transformation and relates to a specific polytope via index-pairs.

\subsection{N\'ed\'elec elements of the second type}
We construct the base functions for the N\'ed\'elec element of the second type using the polytopal template methodology introduced in \cite{skypoly}.
The template sets read
\begin{align}
		\tem_1 &= \{\vb{e}_2,\vb{e}_1\} \, , & \tem_2 &= \{\vb{e}_1 + \vb{e}_2,\vb{e}_1\} \, , & \tem_3 &= \{\vb{e}_1 + \vb{e}_2,-\vb{e}_2\} \, , \notag \\
		\tem_{12} &= \{\vb{e}_2,-\vb{e}_1\} \, , & \tem_{13} &= \{\vb{e}_1,\vb{e}_2\} \, , & \tem_{23} &= \{ (1/2) (\vb{e}_1 - \vb{e}_2),\vb{e}_1 + \vb{e}_2\} \, , \notag \\
		\tem_{123} &= \{\vb{e}_1,\vb{e}_2\} \, .
	\end{align}
The space of B\'ezier polynomials is split across the polytopes of the reference triangle into
\begin{align}
    \Ber^p(\Gamma) = \left \{ \bigoplus_{i=1}^3 \ver^p_i (\Gamma) \right \} \oplus \left \{ \bigoplus_{j \in \mathcal{J}} \edge^p_{j} (\Gamma) \right \} \oplus  \cell^p_{123}(\Gamma) \, , && \mathcal{J} = \{(1,2),(1,3),(2,3)\} \, , 
\end{align}
where $\ver^p_i$ are the sets of the vertex base functions, $\edge^p_{j}$ are the sets of edge base functions, $\cell^p_{123}$ is the set of cell base functions, and the $\oplus$ indicates summation over non-overlapping spaces.
Consequently, the N\'ed\'elec basis is given by
\begin{align}
			&\Nedtwo^p = \left \{ \bigoplus_{i=1}^3 \ver^p_i \otimes \tem_i \right \} \oplus \left \{ \bigoplus_{j \in \mathcal{J}} \edge^p_{j} \otimes \tem_{j} \right \} \oplus  \{\cell^p_{123} \otimes \tem_{123}\} \, , && \mathcal{J} = \{(1,2),(1,3),(2,3)\} \, . 
			\label{eq:ber-decomp}
		\end{align}
Using the B\'ezier basis one finds the following base functions, which inherit the optimal complexity of the underlying basis. 
	\begin{definition} [B\'ezier-N\'ed\'elec II triangle basis]
		
		The following base functions are defined on the reference triangle.
		\begin{itemize}
			\item On the edges the base function reads
			\begin{align}
				&e_{12}: & \bm{\vartheta}(\xi, \eta) &= b_{00}^p \vb{e}_2 \, , & \bm{\vartheta}(\xi, \eta) &= b_{0p}^p (\vb{e}_1 + \vb{e}_2)  \,, \notag \\ 
				&& \bm{\vartheta}(\xi, \eta) &= b_{0j}^p\vb{e}_2 \, , \quad 0<j<p \, , \notag \\ 
				&e_{13}: & \bm{\vartheta}(\xi, \eta) &= b_{00}^p \vb{e}_1 \, , & \bm{\vartheta}(\xi, \eta) &= b_{p0}^p (\vb{e}_1 + \vb{e}_2)  \,, \notag \\  
				&& \bm{\vartheta}(\xi, \eta) &= b_{i0}^p\vb{e}_1 \, , \quad 0<i<p \, , \notag \\ 
				&e_{23}: & \bm{\vartheta}(\xi, \eta) &= b_{0p}^p \vb{e}_1 \, , & \bm{\vartheta}(\xi, \eta) &= -b_{p0}^p \vb{e}_2  \,, \notag \\ && \bm{\vartheta}(\xi, \eta) &= (1/2)\, b_{i,p-i}^p(\vb{e}_1 - \vb{e}_2) \, , \quad 0<i<p \, , 
			\end{align}
			where the first two base functions for each edge are the vertex-edge base functions and the third equation generates pure edge base functions. 
			\item The cell base functions read
			\begin{align}
				&c_{123}: & \bm{\vartheta}(\xi,\eta) &= -b_{0j}^p\vb{e}_1 \, , && 0<j<p \, , \notag \\
				&& \bm{\vartheta}(\xi,\eta) &= b_{i0}^p \vb{e}_2 \, , && 0<i<p \, , \notag \\
				&& \bm{\vartheta}(\xi,\eta) &= b_{i,p-i}^p (\vb{e}_1 + \vb{e}_2) \, , && 0<i<p \, , \notag \\ &&\bm{\vartheta}(\xi,\eta) &= b_{ij}^p \vb{e}_2 \, , && 0<i<p \, , \quad 0<j<p-i \, , \notag \\
				&&\bm{\vartheta}(\xi,\eta) &= b_{ij}^p \vb{e}_1 \, , && 0<i<p \, , \quad 0<j<p-i \, ,
			\end{align}
			where the first three are the respective edge-cell base functions. The remaining two are pure cell base functions.
		\end{itemize}
	\end{definition}

\subsection{N\'ed\'elec elements of the first type}
In order to construct the N\'ed\'elec element of the first type we rely on the construction of the kernel introduced in \cite{Zaglmayr2006} via the exact de Rham sequence and the polytopal template for the non-kernel base functions following \cite{skypoly}.
The complete N\'ed\'elec space reads
	\begin{align}
		\Ned^p &= \Ned^0 \oplus  \left \{ \bigoplus_{j \in \mathcal{J} } \nabla \edge^{p+1}_j \right \} \oplus \nabla \cell^{p+1}_{123} \oplus \left \{ \bigoplus_{i = 1}^2 \ver^{p}_i \otimes \tem_i  \right \} \oplus \left \{ \bigoplus_{j \in \mathcal{J}} \edge^{p}_j \otimes \tem_j  \right \} \oplus \{ \cell^{p}_{123} \otimes \tem_{123} \} \, , \notag \\ \mathcal{J} &= \{(1,2),(1,3),(2,3)\} \, ,
	\end{align}
	where we relied on the decomposition \cref{eq:ber-decomp}.
	Applying the construction to the B\'ezier basis yields the following base functions.
    \begin{definition} [B\'ezier-N\'ed\'elec I triangle basis]
    	
    	We define the base functions on the reference triangle. 
    	\begin{itemize}
    		\item On the edges we employ the lowest order N\'ed\'elec base functions and the edge gradients
    		\begin{align}
    			&e_{12}: & \bm{\vartheta}(\xi, \eta) &=  \bm{\vartheta}_1^I \, ,
    			\notag \\
    			&& \bm{\vartheta}(\xi, \eta) &= \nabla_\xi b_{0j}^{p+1} \, , & & 0< j < p+1 \, , \notag \\
    			&e_{13}: &\bm{\vartheta}(\xi, \eta) &=  \bm{\vartheta}_2^I \, ,
    			\notag \\
    			&& \bm{\vartheta}(\xi, \eta) &= \nabla_\xi b_{i0}^{p+1} \, , & & 0< i < p+1 \, , \notag \\
    			&e_{23}: & \bm{\vartheta}(\xi, \eta) &=  \bm{\vartheta}_3^I \, ,
    			\notag \\
    			&& \bm{\vartheta}(\xi, \eta) &= \nabla_\xi b_{i,p+1-i}^{p+1} \, , & & 0< i < p+1 \, .
    		\end{align}
    	    \item The cell functions read
    	    \begin{align}
    	    	&c_{123}: & \bm{\vartheta}(\xi, \eta) &= b_{00}^{p} \bm{\vartheta}_3^I \, , \notag \\
    	    	& & \bm{\vartheta}(\xi, \eta) &= b_{0p}^{p} \bm{\vartheta}_2^I \, , \notag \\
    	    	& & \bm{\vartheta}(\xi, \eta) &=  b_{0j}^{p}(\bm{\vartheta}_3^I - \bm{\vartheta}_2^I) \, , & & 0 < j < p \, , \notag \\
    	    	& & \bm{\vartheta}(\xi, \eta) &=  b_{i0}^{p}(\bm{\vartheta}_1^I + \bm{\vartheta}_3^I) \, , & & 0 < i < p \, , \notag \\
    	    	& & \bm{\vartheta}(\xi, \eta) &=  b_{i,p-i}^{p}(\bm{\vartheta}_1^I - \bm{\vartheta}_2^I) \, , & & 0 < i < p \, , \notag \\
    	    	& & \bm{\vartheta}(\xi, \eta) &=  b_{ij}^{p}(\bm{\vartheta}_1^I - \bm{\vartheta}_2^I + \bm{\vartheta}_3^I) \, , & & 0 < i < p \, , \quad 0<j<p-i \, , \notag \\
    	    	& & \bm{\vartheta}(\xi, \eta) &= \nabla_\xi b_{ij}^{p+1} \, , & & 0 < i < p+1 \, , \quad 0<j<p+1-i \, ,
    	    \end{align}
            where the last formula gives the cell gradients and the remaining base functions are non-gradients.
    	\end{itemize}
    \end{definition}
The definition relies on the base functions of the lowest order N\'ed\'elec element of the first type \cite{skypoly,Anjam2015}
\begin{align}
		&\bm{\vartheta}^I_1(\xi,\eta) = \begin{bmatrix} 
			\eta \\ 1 - \xi
		\end{bmatrix} \, , &&
		\bm{\vartheta}^I_2(\xi,\eta) = \begin{bmatrix} 
			1 - \eta \\ \xi
		\end{bmatrix} \, , &&
		\bm{\vartheta}^I_3(\xi,\eta) = \begin{bmatrix} 
			\eta \\ - \xi
		\end{bmatrix} \, .
		\label{eq:nedtri}
	\end{align}

\section{Tetrahedral elements}
\begin{figure}
    	\centering
    	\definecolor{asl}{rgb}{0.4980392156862745,0.,1.}
    	\definecolor{asb}{rgb}{0.,0.4,0.6}
    	\begin{tikzpicture}
    		\begin{axis}
    			[
    			width=30cm,height=17cm,
    			view={50}{15},
    			enlargelimits=true,
    			xmin=-1,xmax=2,
    			ymin=-1,ymax=2,
    			zmin=-1,zmax=2,
    			domain=-10:10,
    			axis equal,
    			hide axis
    			]
    			\draw[-to, line width=1.pt, color=black](0., 0., 0.)--(1.5,0.,0.);
    			\draw[color=black] (1.6,0,0) node[] {$\xi$};
    			\draw[-to, line width=1.pt, color=black](0., 0., 0.)--(0.,1.5,0.);
    			\draw[color=black] (0.,1.6,0) node[] {$\eta$};
    			\draw[-to, line width=1.pt, color=black](0., 0., 0.)--(0.,0.,1.5);
    			\draw[color=black] (0.,0.,1.6) node[] {$\zeta$};
    			\addplot3[color=asb][line width=0.6pt,mark=*]
    			coordinates {(0,0,0)(1,0,0)(0,1,0)(0,0,0)};
    			\addplot3[color=asb][line width=1.pt,mark=*]
    			coordinates {(0,0,0)(0,0,1)};
    			\addplot3[color=asb][line width=0.6pt]coordinates {(1,0,0)(0,0,1)};
    			\addplot3[color=asb][line width=0.6pt]coordinates {(0,1,0)(0,0,1)};
    			\draw[color=asb] (0.2,0.2,0.15) node[anchor=south east] {$\Omega$};
    			\draw[color=asb] (0,0,0) node[anchor=south east] {$_{v_{1}}$};
    			\draw[color=asb] (1,0,0) node[anchor=north east] {$_{v_{4}}$};
    			\draw[color=asb] (0,1,0) node[anchor=south west] {$_{v_{3}}$};
    			\draw[color=asb] (0,0,1) node[anchor=north east] {$_{v_{2}}$};
    			\fill[opacity=0.1, asb] (axis cs: 0,0,0) -- (axis cs: 1,0,0) -- (axis cs: 0,1,0) -- (axis cs: 0,0,1) -- cycle;
    			\draw[to-, line width=1.pt, color=asl](0.8,0.2,0)--(0.2,0.8,0);
    			\draw[color=asl] (0.55,0.5,0.) node[anchor=west] {$\bm{\tau}$};
    			\draw[to-, line width=1.pt, color=asl](0.33,0.33,0.33)--(0.5,0.5,0.5);
    			\draw[color=asl] (0.5,0.5,0.55) node[anchor=west] {$\bm{\nu}$};
    			
    			\addplot3[color=asb][line width=0.6pt,mark=*]
    			coordinates {(3.6,3,0.5)(3.6,4,0.5)(4,3,0.2)(3,4,1.1)};
    			\addplot3[color=asb][line width=0.6pt]coordinates {(3.6,3,0.5)(4,3,0.2)};
    			\addplot3[color=asb][line width=0.6pt]coordinates {(3.6,3,0.5)(3,4,1.1)};
    			\addplot3[color=asb][line width=0.6pt]coordinates {(3.6,4,0.5)(3,4,1.1)};
    			\fill[opacity=0.1, asb] (axis cs: 3.6,3,0.5) -- (axis cs: 4,3,0.2) -- (axis cs: 3.6,4,0.5) -- (axis cs: 3,4,1.1) -- cycle;
    			\draw[color=asb] (axis cs: 3.25,4.25,0.47) node[anchor=east] {$\body_e$};
    			\draw[color=asb] (axis cs: 3.6,3,0.5) node[anchor=east] {$_{\vb{x}_{2}}$};
    			\draw[color=asb] (axis cs: 3.6,4,0.5) node[anchor=west] {$_{\vb{x}_{1}}$};
    			\draw[color=asb] (axis cs: 4,3,0.2) node[anchor=east] {$_{\vb{x}_{3}}$};
    			\draw[color=asb] (axis cs: 3,4,1.1) node[anchor=east] {$_{\vb{x}_{4}}$};
    			\draw[-to, line width=1.pt, color=black](2.5, 2.5, 0.)--(3.5, 2.5,0.);
    			\draw[color=black] (3.6, 2.5,0.) node[] {$x$};
    			\draw[-to, line width=1.pt, color=black](2.5, 2.5, 0.)--(2.5, 3.5,0.);
    			\draw[color=black] (2.5, 3.6,0.) node[] {$y$};
    			\draw[-to, line width=1.pt, color=black](2.5, 2.5, 0.)--(2.5, 2.5,1.);
    			\draw[color=black] (2.5, 2.5,1.1) node[] {$z$};
    			
    			\draw[-to, line width=1.pt, color=asl](3.8 , 3.2 , 0.38)--(3.2 , 3.8 , 0.92);
    			\draw[color=asl] (3.16 , 3.76 , 0.5) node[anchor=west] {$\vb{t}$};
    			\draw[-to, line width=1.pt, color=asl](3.508, 3.31 , 0.594)--(3.3, 3. , 0.7);
    			\draw[color=asl] (3.3, 3. , 0.71) node[anchor=east] {$\vb{n}$};
    			
    			\draw[-Triangle, line width=0.6pt,color=black] (1.3,1.3,0.51) -- (2,2,0.536);
    			\draw[color=black] (1.65,1.65,0.55) node[anchor=south] {$\vb{x}:\Omega \to \body_e$};
    		\end{axis}
    	\end{tikzpicture}
    	\caption{Barycentric mapping of the reference tetrahedron to an element in the physical domain.}
    	\label{fig:tetmap}
    \end{figure}
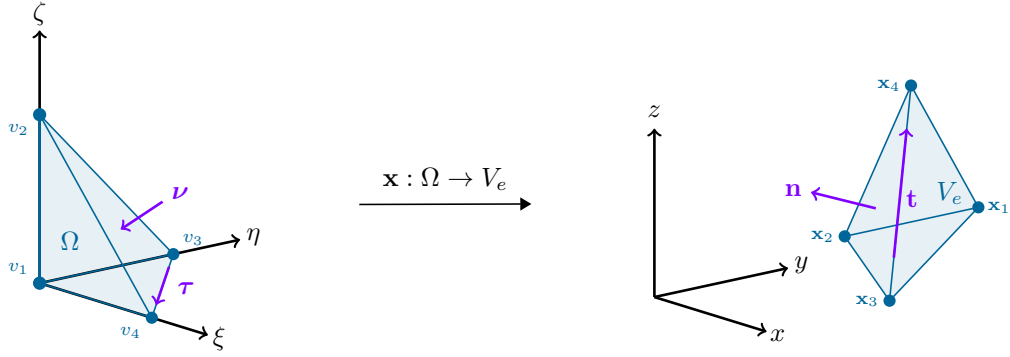

The tetrahedral elements are mapped from the reference tetrahedron $\Omega$ by the three-dimensional barycentric coordinates onto the physical domain $\body_e$, see \cref{fig:tetmap}
	\begin{align}
		\vb{x}(\xi, \eta, \zeta) &= (1 - \xi - \eta - \zeta) \vb{x}_1 + \zeta \, \vb{x}_2 + \eta \, \vb{x}_3 + \xi \, \vb{x}_4 \, , 
		&&\vb{x}:\Omega \to \body_e \, , \notag \\
		 \Omega &= \{(\xi,\eta,\zeta) \in [0,1]^3 \; | \; \xi + \eta + \zeta \leq 1\} \, .
	\end{align}
	The corresponding Jacobi matrix reads
	\begin{align}
		\bm{J} = \D \vb{x} = \begin{bmatrix}
			\vb{x}_4 - \vb{x}_1, & \vb{x}_3 - \vb{x}_1, & \vb{x}_2 - \vb{x}_1
		\end{bmatrix} \in \R^{3 \times 3} \, .
	\end{align}

\subsection{The Bernstein-B\'ezier basis for tetrahedra}
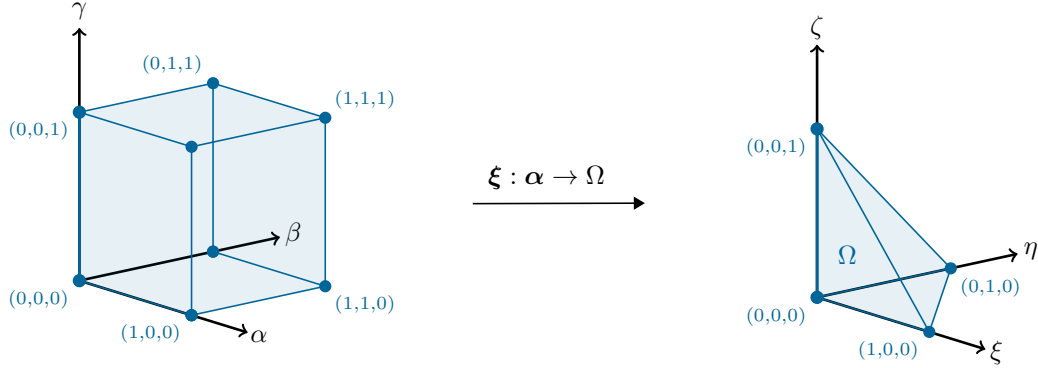
\begin{figure}
    	\centering
    	\definecolor{asl}{rgb}{0.4980392156862745,0.,1.}
    	\definecolor{asb}{rgb}{0.,0.4,0.6}
    	\begin{tikzpicture}
    		\begin{axis}
    			[
    			width=30cm,height=17cm,
    			view={50}{15},
    			enlargelimits=true,
    			xmin=-1,xmax=2,
    			ymin=-1,ymax=2,
    			zmin=-1,zmax=2,
    			domain=-10:10,
    			axis equal,
    			hide axis
    			]
    			\draw[-to, line width=1.pt, color=black](0., 0., 0.)--(1.5,0.,0.);
    			\draw[color=black] (1.6,0,0) node[] {$\alpha$};
    			\draw[-to, line width=1.pt, color=black](0., 0., 0.)--(0.,1.5,0.);
    			\draw[color=black] (0.,1.6,0) node[] {$\beta$};
    			\draw[-to, line width=1.pt, color=black](0., 0., 0.)--(0.,0.,1.5);
    			\draw[color=black] (0.,0.,1.6) node[] {$\gamma$};
    			\addplot3[color=asb][line width=0.6pt,mark=*]
    			coordinates {(0,0,0)(1,0,0)};
    			\addplot3[color=asb][line width=1.pt,mark=*]
    			coordinates {(0,0,0)(0,0,1)};
    			\addplot3[color=asb][line width=0.6pt,mark=*]coordinates {(1,0,0)(1,1,0)(0,1,0)};
    			\addplot3[color=asb][line width=0.6pt,mark=*]coordinates {(1,0,0)(1,0,1)(1,1,1)(0,1,1)(0,1,0)};
    			\addplot3[color=asb][line width=0.6pt]coordinates {(1,1,0)(1,1,1)};
    			\addplot3[color=asb][line width=0.6pt]coordinates {(1,0,1)(0,0,1)(0,1,1)};
    			\draw[color=asb] (0,0,0) node[anchor=north east] {$_{(0,0,0)}$};
    			\draw[color=asb] (1,0,0) node[anchor=north east] {$_{(1,0,0)}$};
    			\draw[color=asb] (0,0,1) node[anchor=north east] {$_{(0,0,1)}$};
    			\draw[color=asb] (1,1,0) node[anchor=north west] {$_{(1,1,0)}$};
    			\draw[color=asb] (1,1,1) node[anchor=south west] {$_{(1,1,1)}$};
    			\draw[color=asb] (0,1,1) node[anchor=south east] {$_{(0,1,1)}$};
    			\fill[opacity=0.1, asb] (axis cs: 0,0,0) -- (axis cs: 1,0,0) -- (axis cs: 1,0,1) -- (axis cs: 0,0,1) -- cycle;
    			\fill[opacity=0.1, asb] (axis cs: 1,0,0) -- (axis cs: 1,1,0) -- (axis cs: 1,1,1) -- (axis cs: 1,0,1) -- cycle;
    			\fill[opacity=0.1, asb] (axis cs: 0,0,1) -- (axis cs: 1,0,1) -- (axis cs: 1,1,1) -- (axis cs: 0,1,1) -- cycle;
    			
    			\draw[-to, line width=1.pt, color=black](3., 3., 0.)--(4.5,3.,0.);
    			\draw[color=black] (4.6,3,0) node[] {$\xi$};
    			\draw[-to, line width=1.pt, color=black](3., 3., 0.)--(3.,4.5,0.);
    			\draw[color=black] (3.,4.6,0) node[] {$\eta$};
    			\draw[-to, line width=1.pt, color=black](3., 3., 0.)--(3.,3.,1.5);
    			\draw[color=black] (3.,3.,1.6) node[] {$\zeta$};
    			\addplot3[color=asb][line width=0.6pt,mark=*]
    			coordinates {(3,3,0)(4,3,0)(3,4,0)(3,3,0)};
    			\addplot3[color=asb][line width=1.pt,mark=*]
    			coordinates {(3,3,0)(3,3,1)};
    			\addplot3[color=asb][line width=0.6pt]coordinates {(4,3,0)(3,3,1)};
    			\addplot3[color=asb][line width=0.6pt]coordinates {(3,4,0)(3,3,1)};
    			\draw[color=asb] (3.2,3.2,0.15) node[anchor=south east] {$\Omega$};
    			\draw[color=asb] (3,3,0) node[anchor=north east] {$_{(0,0,0)}$};
    			\draw[color=asb] (4,3,0) node[anchor=north east] {$_{(1,0,0)}$};
    			\draw[color=asb] (3,4,0) node[anchor=north west] {$_{(0,1,0)}$};
    			\draw[color=asb] (3,3,1) node[anchor=north east] {$_{(0,0,1)}$};
    			\fill[opacity=0.1, asb] (axis cs: 3,3,0) -- (axis cs: 4,3,0) -- (axis cs: 3,4,0) -- (axis cs: 3,3,1) -- cycle;
    			
    			\draw[-Triangle, line width=0.6pt,color=black] (1.6,1.6,0.51) -- (2.3,2.3,0.536);
    			\draw[color=black] (1.9,1.9,0.55) node[anchor=south] {$\bm{\xi}:\bm{\alpha} \to \Omega$};
    		\end{axis}
    	\end{tikzpicture}
    	\caption{Duffy mapping of the unit hexahedron to the unit tetrahedron.}
    	\label{fig:tetduffy}
    \end{figure}

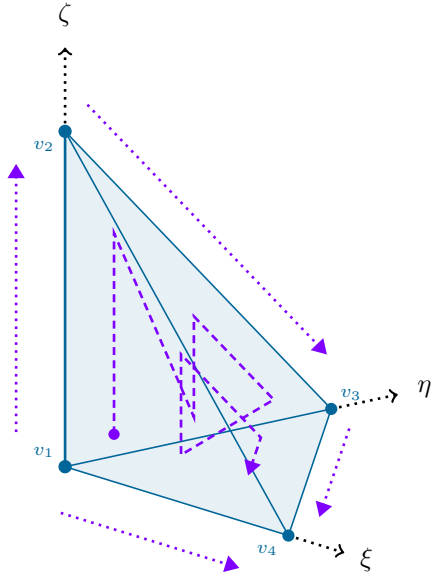
\begin{figure}
    	\centering
    	\definecolor{asl}{rgb}{0.4980392156862745,0.,1.}
    	\definecolor{asb}{rgb}{0.,0.4,0.6}
    	\begin{tikzpicture}
    		\begin{axis}
    			[
    			width=30cm,height=25cm,
    			view={50}{15},
    			enlargelimits=true,
    			xmin=-1,xmax=2,
    			ymin=-1,ymax=2,
    			zmin=-1,zmax=2,
    			domain=-10:10,
    			axis equal,
    			hide axis
    			]
    			\draw[-to, line width=1.pt, color=black, dotted](1., 0., 0.)--(1.25,0.,0.);
    			\draw[color=black] (1.35,0,0) node[] {$\xi$};
    			\draw[-to, line width=1.pt, color=black, dotted](0., 1., 0.)--(0.,1.25,0.);
    			\draw[color=black] (0.,1.35,0) node[] {$\eta$};
    			\draw[-to, line width=1.pt, color=black, dotted](0., 0., 1.)--(0.,0.,1.25);
    			\draw[color=black] (0.,0.,1.35) node[] {$\zeta$};
    			\addplot3[color=asb][line width=0.6pt,mark=*]
    			coordinates {(0,0,0)(1,0,0)(0,1,0)(0,0,0)};
    			\addplot3[color=asb][line width=1.pt,mark=*]
    			coordinates {(0,0,0)(0,0,1)};
    			\addplot3[color=asb][line width=0.6pt]coordinates {(1,0,0)(0,0,1)};
    			\addplot3[color=asb][line width=0.6pt]coordinates {(0,1,0)(0,0,1)};
    			\draw[color=asb] (0,0,0) node[anchor=south east] {$_{v_{1}}$};
    			\draw[color=asb] (1,0,0) node[anchor=north east] {$_{v_{4}}$};
    			\draw[color=asb] (0,1,0) node[anchor=south west] {$_{v_{3}}$};
    			\draw[color=asb] (0,0,1) node[anchor=north east] {$_{v_{2}}$};
    			\fill[opacity=0.1, asb] (axis cs: 0,0,0) -- (axis cs: 1,0,0) -- (axis cs: 0,1,0) -- (axis cs: 0,0,1) -- cycle;
    			
    			\draw (0.1, 0.1, 0.1) node[circle,fill=asl,inner sep=1.5pt] {};
    			\draw[-Triangle, line width=1.pt, color=asl, densely dashed](0.1, 0.1, 0.1)--(0.1,0.1,0.7)--(0.1,0.4,0.1)
    			--(0.1,0.4,0.4)--(0.1,0.7,0.1)--(0.4,0.1,0.1)
    			--(0.4,0.1,0.4)--(0.4,0.4,0.1)--(0.7,0.1,0.1);
    			\draw[-Triangle, line width=1.pt, color=asl, dotted](-0.1, -0.1, 0.1)--(-0.1,-0.1,0.9);
    			\draw[-Triangle, line width=1.pt, color=asl, dotted](0.1, -0.1, -0.1)--(0.9,-0.1,-0.1);
    			\draw[-Triangle, line width=1.pt, color=asl, dotted](0.1, 0, 1.1)--(0.1,0.9,0.2);
    			\draw[-Triangle, line width=1.pt, color=asl, dotted](0.2, 0.9, 0)--(0.9,0.2,0);
    		\end{axis}
    	\end{tikzpicture}
    	\caption{Order of traversal of tetrahedral B\'ezier base functions on the unit tetrahedron. The traversal order on each face agrees with an orientation of the vertices $f_{ijk} =\{v_i,v_j,v_k\}$ such that $i<j<k$. The traversal order on each edge is from the lower index vertex to the higher index vertex.}
    	\label{fig:tettrav}
    \end{figure}

\begin{figure}
		\centering
		\begin{subfigure}{0.24\linewidth}
			\centering
			\includegraphics[width=0.7\linewidth]{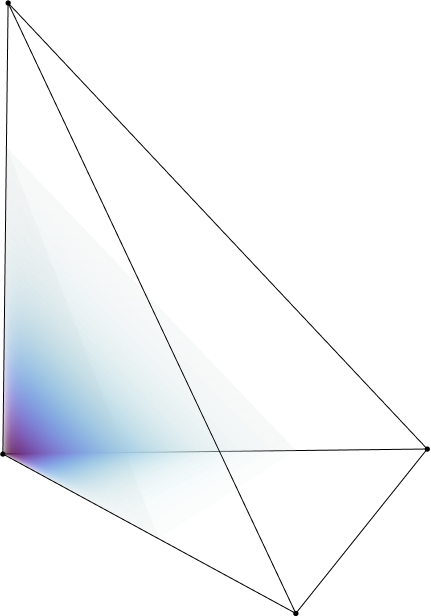}
			\caption{}
		\end{subfigure}
	    \begin{subfigure}{0.24\linewidth}
	    	\centering
	    	\includegraphics[width=0.7\linewidth]{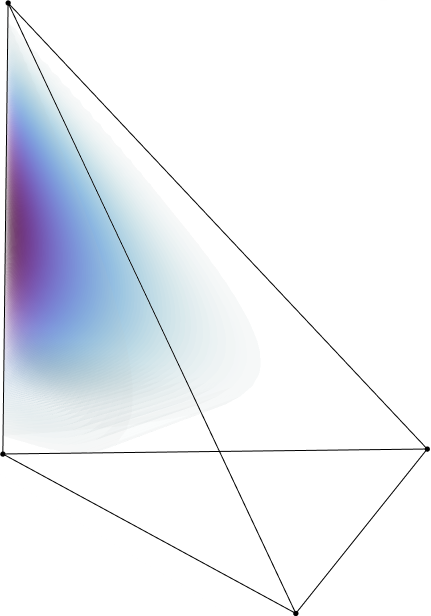}
	    	\caption{}
	    \end{subfigure}
        \begin{subfigure}{0.24\linewidth}
        	\centering
        	\includegraphics[width=0.7\linewidth]{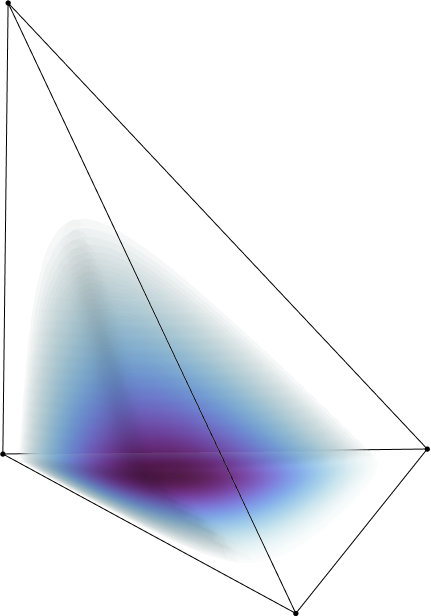}
        	\caption{}
        \end{subfigure}
        \begin{subfigure}{0.24\linewidth}
        	\centering
        	\includegraphics[width=0.7\linewidth]{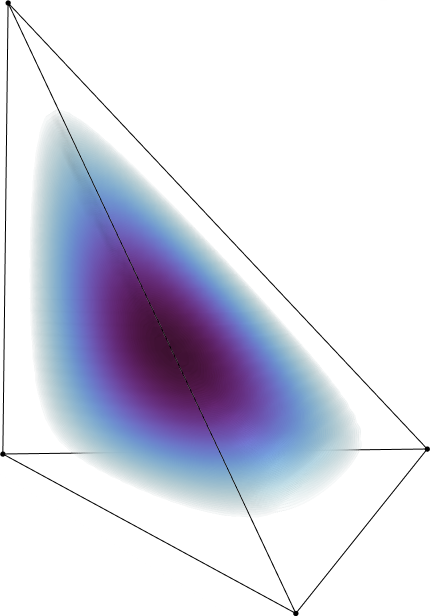}
        	\caption{}
        \end{subfigure}
		\caption{Quartic B\'ezier vertex (a), edge (b), face (c), and cell (c) base functions  on the reference tetrahedron.}
		\label{fig:beziertet}
	\end{figure}

Analogously to triangle elements, the B\'ezier tetrahedra on the unit tetrahedron $\Omega$ are defined using the barycentric coordinates by expanding the coefficients of
	\begin{align}
		(\lambda_1 + \lambda_2 + \lambda_3 + \lambda_4)^p = ([1- \xi - \eta - \zeta] + \zeta + \eta + \xi)^p = 1 \, ,
	\end{align}
    thus finding
	\begin{align}
		b_{ijk}^p(\lambda_1,\lambda_2,\lambda_3,\lambda_4) = \begin{pmatrix} p \\ i \end{pmatrix}
		\begin{pmatrix} p-i \\ j \end{pmatrix}
		\begin{pmatrix} p-i-j \\ k \end{pmatrix}
		\lambda_1^{p-i-j-k} \lambda_2^k \lambda_3^j \lambda_4^k \, ,
	\end{align}
	with the equivalent trivariate form
	\begin{align}
		b^p_{ijk}(\xi,\eta,\zeta) = \begin{pmatrix} p \\ i \end{pmatrix}
		\begin{pmatrix} p-i \\ j \end{pmatrix} 
		\begin{pmatrix} p-i-j \\ k \end{pmatrix}
		(1- \xi - \eta-\zeta)^{p-i-j-k} \zeta^k \eta^j \xi^i \, .
	\end{align}
    We construct the Duffy transformation by mapping the unit tetrahedron as a collapsed hexahedron
    \begin{align}
    	&\bm{\xi}: [0,1]^3 \to \Omega \, , && \{\alpha, \beta,\gamma\}\mapsto \{\xi,\eta,\zeta\} \, ,
    \end{align}
    using the relations
    \begin{align}
    	\xi &= \alpha \, , & \eta &= (1 - \alpha) \beta \, , & \zeta &= (1- \alpha)(1-\beta)\gamma \, , \notag \\
    	\alpha &= \xi \, , & \beta &=  \dfrac{\eta}{1-\xi} \, , & \gamma &= \dfrac{\zeta}{1- \xi - \eta} \, ,
    \end{align}
    as depicted in \cref{fig:tetduffy}.
    Applying the Duffy transformation to B\'ezier tetrahedra 
    \begin{align}
    	b^p_{ijk}(\xi,\eta,\zeta) &= \begin{pmatrix} p \\ i \end{pmatrix}
    	\begin{pmatrix} p-i \\ j \end{pmatrix} 
    	\begin{pmatrix} p-i-j \\ k \end{pmatrix}
    	(1- \xi - \eta-\zeta)^{p-i-j-k} \zeta^k \eta^j \xi^i \notag \\
    	 &= \begin{pmatrix} p \\ i \end{pmatrix}
    	\begin{pmatrix} p-i \\ j \end{pmatrix} 
    	\begin{pmatrix} p-i-j \\ k \end{pmatrix}
    	(1- \alpha - (1-\alpha)\beta-(1-\alpha)(1-\beta)\gamma)^{p-i-j-k} \notag \\
    	& \qquad \cdot (1-\alpha)^k(1-\beta)^k\gamma^k (1-\alpha)^j\beta^j \alpha^i \notag \\
    	&= \begin{pmatrix} p \\ i \end{pmatrix}
    	\begin{pmatrix} p-i \\ j \end{pmatrix} 
    	\begin{pmatrix} p-i-j \\ k \end{pmatrix}
    	(1-\alpha)^{p-i-j-k}(1-\beta)^{p-i-j-k}(1-\gamma)^{p-i-j-k} \\
    	& \qquad \cdot (1-\alpha)^k(1-\beta)^k\gamma^k (1-\alpha)^j\beta^j \alpha^i \notag \\
    	&= \begin{pmatrix} p \\ i \end{pmatrix}
    	(1-\alpha)^{p-i}\alpha^i
    	\begin{pmatrix} p-i \\ j \end{pmatrix} 
    	(1-\beta)^{p-i-j} \beta^j
    	\begin{pmatrix} p-i-j \\ k \end{pmatrix}
    	(1-\gamma)^{p-i-j-k}\gamma^k \notag \\
    	&=  b_i^p(\alpha) b_j^{p-i}(\beta)b_k^{p-i-j}(\gamma) \, , \notag
    \end{align}
    leads to an intrinsic factorization via univariate Bernstein base functions, which allow for fast evaluations using sum factorization \cite{AinsworthOpt}. Further, since the pair $b_j^{p-i}(\beta)b_k^{p-i-j}(\gamma)$ spans a B\'ezier triangle, it is clear that the multiplication with $b_i^p(\alpha)$ interpolates between that triangle and a point in space, effectively spanning a tetrahedron. In order to compute gradients the chain rule is employed with respect to the Duffy transformation
    \begin{align}
    	\nabla_\xi b_{ijk}^p &= (\D_\alpha \bm{\xi})^{-T} \nabla_\alpha b_{ijk}^p \, , \qquad \qquad \D_\alpha \bm{\xi} = \begin{bmatrix}
    		1 & 0 & 0 \\
    		-\beta & 1-\alpha & 0 \\
    		(\beta - 1)\gamma & (\alpha - 1)\gamma & (1-\alpha)(1-\beta)
    	\end{bmatrix} \, , \notag \\
        (\D_\alpha \bm{\xi})^{-T} &= \dfrac{1}{(1-\alpha)(1-\beta)} \begin{bmatrix}
        	(1-\alpha)(1-\beta) &  (1-\beta)\beta & \gamma\\
        0 & 1-\beta & \gamma\\
        0 & 0 & 1
    \end{bmatrix} \, .
    \end{align}
    We use dual numbers to compute the derivative of each Bernstein base function and construct the $\bm{\alpha}$-gradient
    \begin{align}
    	\nabla_\alpha b_{ijk}^p(\alpha, \beta, \gamma) = \begin{bmatrix}
    		b_j^{p-i} b_k^{p-i-j} \dfrac{\dd}{\dd \alpha} b_i^p \\[2ex]
    		b_i^p b_k^{p-i-j} \dfrac{\dd }{\dd \beta} b_j^{p-i} \\[2ex]
    		b_i^p b_j^{p-i} \dfrac{\dd }{\dd \gamma} b_k^{p-i-j}
    	\end{bmatrix} \, .
    \end{align} 
    The Duffy transformation results in the optimal order of traversal of the base functions depicted in \cref{fig:tettrav}. Note that the traversal order agrees with the oriental definitions introduced in \cite{SKY2022115298} and each oriented face has the same order of traversal as the triangle \cref{fig:traversaltri}.
    We relate the base functions to their respective polytopes using the index triplets.
    \begin{observation} [Tetrahedron base functions]
    	
    	The polytope of each base function $b_{ijk}^p(\xi,\eta,\zeta)$ is determined as follows.
    	\begin{itemize}
    		\item the indices $(0,0,0),(0,0,p),(0,p,0)$ and $(p,0,0)$ represent the respective vertex base functions;
    	    \item the first edge is associated with the triplet $(0,0,k)$ where $0<k<p$, the second with $(0,j,0)$ where $0<j<p$ and the third with $(i,0,0)$ where $0<i<p$. The slated edges are given by $(0,j,p-j)$ with $0<j<p$, $(i,0,p-i)$ with $0<i<p$ and $(i,p-i,0)$ with $0<i<p$, respectively;
    	    \item the base functions of the first face are given by 
    	    $(0,j,k)$ with $0<j<p$ and $0<k<p-j$. The second face is associated with the base functions given by the triplets 
    	    $(i,0,k)$ with $0<i<p$ and $0<k<p-i$. The base functions of the third face are related to the indices $(i,j,0)$ with $0<i<p$ and $0<j<p-i$. Lastly, the base functions of the slated face are given by $(i,j,p-i-j)$ with $0<i<p$ and $0<j<p-i$;
    	    \item the remaining indices correspond to the cell base functions.
    	\end{itemize}
    \end{observation}
    Examples of B\'ezier base functions on their respective polytopes are depicted in \cref{fig:beziertet}. 

\subsection{N\'ed\'elec elements of the second type}
The B\'ezier polynomial space is split according to the polytopes of the reference tetrahedron
\begin{align}
    &\Ber^p(\Omega) = \left\{ \bigoplus_{i=1}^4 \ver_i^p(\Omega)  \right\} \oplus \left\{ \bigoplus_{j \in \mathcal{J}  } \edge^p_j(\Omega)  \right\} \oplus \left\{ \bigoplus_{k \in \mathcal{K}} \face^p_k(\Omega) \right\} \oplus \cell^p_{1234}(\Omega) \, , \notag \\
    		&\mathcal{J} = \{ (1,2),(1,3),(1,4),(2,3),(2,4),(3,4) \} \, , \qquad \mathcal{K} = \{ (1,2,3),(1,2,4),(1,3,4),(2,3,4) \} \, ,
\end{align}
where $\ver_i^p$ are the sets of vertex base functions, $\edge^p_j$ are the sets of edge base functions, $\face^p_k$ are the sets of face base functions and $\cell^p_{1234}$ is the set of cell base functions.
We apply the template sets from \cite{skypoly}
\begin{align}
    	\tem_1 &= \{ \vb{e}_3,\vb{e}_2,\vb{e}_1 \} \, , & \tem_2 &= \{ \vb{e}_1 + \vb{e}_2 + \vb{e}_3 , \vb{e}_2 , \vb{e}_1 \} \, , & \tem_3 &= \{ \vb{e}_1 + \vb{e}_2 + \vb{e}_3, -\vb{e}_3 ,\vb{e}_1 \} \, , \notag \\
    	\tem_4 &= \{ \vb{e}_1 + \vb{e}_2 + \vb{e}_3, -\vb{e}_3, -\vb{e}_2\} \, , & \tem_{12} &= \{ \vb{e}_3, -\vb{e}_2, -\vb{e}_1 \} \, , & \tem_{13} &= \{ \vb{e}_2, \vb{e}_3, -\vb{e}_1 \} \, , \notag \\
    	\tem_{14} &= \{ \vb{e}_1, \vb{e}_3, \vb{e}_2 \} \, , & \tem_{23} &= \{ \vb{e}_2, \vb{e}_1 + \vb{e}_2 + \vb{e}_3, -\vb{e}_1 \} \, , & \tem_{24} &= \{ \vb{e}_1, \vb{e}_1 + \vb{e}_2 + \vb{e}_3, \vb{e}_2 \} \, , \notag \\
    	\tem_{34} &= \{ \vb{e}_1, \vb{e}_1 + \vb{e}_2 + \vb{e}_3, -\vb{e}_3 \} \, , & \tem_{123} &= \{ \vb{e}_3, \vb{e}_2, -\vb{e}_1 \} \, , & \tem_{124} &= \{ \vb{e}_3, \vb{e}_1, \vb{e}_2 \} \, , \notag \\
    	\tem_{134} &= \{\vb{e}_2, \vb{e}_1, -\vb{e}_3\} \, , & \tem_{234}  &= \{ \vb{e}_2, \vb{e}_1,\vb{e}_1 + \vb{e}_2 + \vb{e}_3\} \, , & \tem_{1234} &= \{ \vb{e}_3, \vb{e}_2, \vb{e}_1 \} \, ,
    \end{align}
to span the N\'ed\'elec element of the second type
    	\begin{align}
    		\Nedtwo^p &= \left\{ \bigoplus_{i=1}^4 \ver_i^p \otimes \tem_i \right\} \oplus \left\{ \bigoplus_{j \in \mathcal{J}  } \edge^p_j \otimes \tem_j \right\} \oplus \left\{ \bigoplus_{k \in \mathcal{K}} \face^p_k \otimes \tem_k  \right\} \oplus \{ \cell^p_{1234} \otimes \tem_{1234} \} \, , \notag \\
    		\mathcal{J} &= \{ (1,2),(1,3),(1,4),(2,3),(2,4),(3,4) \} \, , \qquad \mathcal{K} = \{ (1,2,3),(1,2,4),(1,3,4),(2,3,4) \} \, .
    	\end{align}
    We can now define the B\'ezier-N\'ed\'elec element of the second type for arbitrary powers while inheriting optimal complexity.
    \begin{definition} [B\'ezier-N\'ed\'elec II tetrahedral basis]
    	
    	We define the base functions on the reference tetrahedron:
    	\begin{itemize}
    		\item on the edges the base functions read
    		\begin{align}
    			&e_{12}: & \bm{\vartheta}(\xi,\eta,\zeta) &= b_{000}^p \vb{e}_3 \, , & \bm{\vartheta}(\xi,\eta,\zeta) &= b_{00p}^p (\vb{e}_1+\vb{e}_2+\vb{e}_3) \, , \notag \\
    			& & \bm{\vartheta}(\xi,\eta,\zeta) &= b_{00k}^p \vb{e}_3 \, ,  \quad  0< k <p \, , \notag \\
    			&e_{13}: & \bm{\vartheta}(\xi,\eta,\zeta) &= b_{000}^p \vb{e}_2 \, , & \bm{\vartheta}(\xi,\eta,\zeta) &= b_{0p0}^p (\vb{e}_1+\vb{e}_2+\vb{e}_3) \, , \notag \\
    			& & \bm{\vartheta}(\xi,\eta,\zeta) &= b_{0j0}^p \vb{e}_2 \, ,  \quad  0< j <p \, , \notag \\
    			&e_{14}: & \bm{\vartheta}(\xi,\eta,\zeta) &= b_{000}^p \vb{e}_1 \, , & \bm{\vartheta}(\xi,\eta,\zeta) &= b_{p00}^p (\vb{e}_1+\vb{e}_2+\vb{e}_3) \, , \notag \\
    			& & \bm{\vartheta}(\xi,\eta,\zeta) &= b_{i00}^p \vb{e}_1 \, ,  \quad  0< i <p \, , \notag \\
    			&e_{23}: & \bm{\vartheta}(\xi,\eta,\zeta) &= b_{00p}^p \vb{e}_2 \, , & \bm{\vartheta}(\xi,\eta,\zeta) &= -b_{0p0}^p \vb{e}_3 \, , \notag \\
    			& & \bm{\vartheta}(\xi,\eta,\zeta) &= b_{0j,p-j}^p \vb{e}_2 \, ,  \quad  0< j <p \, , \notag \\
    			&e_{24}: & \bm{\vartheta}(\xi,\eta,\zeta) &= b_{00p}^p \vb{e}_1 \, , & \bm{\vartheta}(\xi,\eta,\zeta) &= -b_{p00}^p \vb{e}_3 \, , \notag \\
    			& & \bm{\vartheta}(\xi,\eta,\zeta) &= b_{i0,p-i}^p \vb{e}_1 \, ,  \quad  0< i <p \, , \notag \\
    			&e_{34}: & \bm{\vartheta}(\xi,\eta,\zeta) &= b_{0p0}^p \vb{e}_1 \, , & \bm{\vartheta}(\xi,\eta,\zeta) &= -b_{p00}^p \vb{e}_2 \, , \notag \\
    			& & \bm{\vartheta}(\xi,\eta,\zeta) &= b_{i,p-i,0}^p \vb{e}_1 \, ,  \quad  0< i <p \, ,
    		\end{align}
    	    where the first two base functions on each edge are the vertex-edge base functions;
    	    \item the face base functions are given by
    	    \begin{align}
    	    	&f_{123}: & \bm{\vartheta}(\xi,\eta,\zeta) &= -b_{00k}^p\vb{e}_2 \, , && 0<k<p \, ,  \notag \\
    	    	&&\bm{\vartheta}(\xi,\eta,\zeta) &= b_{0j0}^p  \vb{e}_3 \, , && 0<j<p \, , \notag \\
    	    	&&\bm{\vartheta}(\xi,\eta,\zeta) &= b_{0j,p-j}^p  (\vb{e}_1+\vb{e}_2+\vb{e}_3) \, , && 0<j<p \, , \notag \\
    	    	&&\bm{\vartheta}(\xi,\eta,\zeta) &= b_{0jk}^p  \vb{e}_3 \, , && 0<j<p \, , \quad 0<k<p-j  \, ,  \notag \\
    	    	&&\bm{\vartheta}(\xi,\eta,\zeta) &= b_{0jk}^p  \vb{e}_2 \, , && 0<j<p \, , \quad 0<k<p-j  \, , \notag \\
    	    	&f_{124}: & \bm{\vartheta}(\xi,\eta,\zeta) &= -b_{00k}^p \vb{e}_1 \, , && 0<k<p \, ,  \notag \\
    	    	&&\bm{\vartheta}(\xi,\eta,\zeta) &= b_{i00}^p  \vb{e}_3 \, , && 0<i<p \, , \notag \\
    	    	&&\bm{\vartheta}(\xi,\eta,\zeta) &= b_{i0,p-i}^p  (\vb{e}_1+\vb{e}_2+\vb{e}_3) \, , && 0<i<p \, , \notag \\
    	    	&&\bm{\vartheta}(\xi,\eta,\zeta) &= b_{i0k}^p  \vb{e}_3 \, , && 0<i<p \, , \quad 0<k<p-i \, ,  \notag \\
    	    	&&\bm{\vartheta}(\xi,\eta,\zeta) &= b_{i0k}^p  \vb{e}_1 \, , && 0<i<p \, , \quad 0<k<p-i \, ,  \notag \\
    	    	&f_{134}: & \bm{\vartheta}(\xi,\eta,\zeta) &= -b_{0j0}^p \vb{e}_1 \, , && 0<j<p \, ,  \notag \\
    	    	&&\bm{\vartheta}(\xi,\eta,\zeta) &= b_{i00}^p  \vb{e}_2 \, , && 0<i<p \, , \notag \\
    	    	&&\bm{\vartheta}(\xi,\eta,\zeta) &= b_{i,p-i,0}^p  (\vb{e}_1+\vb{e}_2+\vb{e}_3) \, , && 0<i<p \, , \notag \\
    	    	&&\bm{\vartheta}(\xi,\eta,\zeta) &= b_{ij0}^p  \vb{e}_2 \, , && 0<i<p \, , \quad 0<j<p-i \, ,  \notag \\
    	    	&&\bm{\vartheta}(\xi,\eta,\zeta) &= b_{ij0}^p  \vb{e}_1 \, , && 0<i<p \, , \quad 0<j<p-i \, ,  \notag \\
    	    	&f_{234}: & \bm{\vartheta}(\xi,\eta,\zeta) &= -b_{0j,p-j}^p \vb{e}_1 \, , && 0<j<p \, ,  \notag \\
    	    	&&\bm{\vartheta}(\xi,\eta,\zeta) &= b_{i0,p-i}^p  \vb{e}_2 \, , && 0<i<p \, , \notag \\
    	    	&&\bm{\vartheta}(\xi,\eta,\zeta) &= -b_{i,p-i,0}^p \vb{e}_3 \, , && 0<i<p \, , \notag \\
    	    	&&\bm{\vartheta}(\xi,\eta,\zeta) &= b_{ij,p-i-j}^p  \vb{e}_2 \, , && 0<i<p \, , \quad 0<j<p-i \, ,  \notag \\
    	    	&&\bm{\vartheta}(\xi,\eta,\zeta) &= b_{ij,p-i-j}^p  \vb{e}_1 \, , && 0<i<p \, , \quad 0<j<p-i \, , 
    	    \end{align}
            where the first three formulas for each face are the edge-face base functions;
            \item finally, the cell base functions read
            \begin{align}
            	&c_{1234}: & \bm{\vartheta}(\xi,\eta,\zeta) &= -b_{0jk}^p \vb{e}_1 \, , && 0<j<p \, , \quad 0<k<p-j \, , \notag \\
            	&& \bm{\vartheta}(\xi,\eta,\zeta) &= b_{i0k}^p \vb{e}_2 \, , && 0<i<p \, , \quad 0<k<p-i \, , \notag \\
            	&& \bm{\vartheta} (\xi,\eta,\zeta) &= -b_{ij0}^p \vb{e}_3 \, , && 0<i<p \, , \quad 0<j<p-i \, , \notag \\
            	&& \bm{\vartheta}(\xi,\eta,\zeta) &= b_{ij,p-i-j}^p (\vb{e}_1 + \vb{e}_2 + \vb{e}_3) \, , && 0<i<p \, , \quad 0<j<p-i \, , \notag \\
            	&& \bm{\vartheta}(\xi,\eta,\zeta) &= b_{ijk}^p \vb{e}_3 \, , && 0<i<p \, , \quad 0<j<p-i \, ,  \quad 0<k<p-i-j \, , \notag \\
            && \bm{\vartheta}(\xi,\eta,\zeta) &= b_{ijk}^p\vb{e}_2 \, , && 0<i<p \, , \quad 0<j<p-i \, , \quad 0<k<p-i-j  \, , \notag \\
        && \bm{\vartheta}(\xi,\eta,\zeta) &= b_{ijk}^p\vb{e}_1 \, , && 0<i<p \, , \quad 0<j<p-i \, , \quad 0<k<p-i-j  \, , 
            \end{align}
        where the first four formulas are the face-cell base functions.
    	\end{itemize}
    \end{definition}

\subsection{N\'ed\'elec elements of the first type}
In order to construct the N\'ed\'elec element of first type on tetrahedra we introduce the template sets 
\begin{align}
    	\tem_{1} &= \{ \bm{\vartheta}^I_{4},\bm{\vartheta}^I_{5},\bm{\vartheta}^I_{6} \} \, , & 
    	\tem_{2} &= \{ -\bm{\vartheta}^I_{2},-\bm{\vartheta}^I_{3},\bm{\vartheta}^I_{6} \} \, , & 
    	\tem_{3} &= \{ -\bm{\vartheta}^I_{3},-\bm{\vartheta}^I_{5} \} \, , 
    	\notag \\
    	\tem_{12} &= \{ \bm{\vartheta}^I_{4}-\bm{\vartheta}^I_{2}, \bm{\vartheta}^I_{5}-\bm{\vartheta}^I_{3} \} \, , &
    	\tem_{13} &= \{ \bm{\vartheta}^I_{1}+\bm{\vartheta}^I_{4}, \bm{\vartheta}^I_{6}-\bm{\vartheta}^I_{3} \} \, , & 
    	\tem_{14} &= \{ \bm{\vartheta}^I_{1}+\bm{\vartheta}^I_{5}, \bm{\vartheta}^I_{2}+\bm{\vartheta}^I_{6} \} \, , 
    	\notag \\
    	\tem_{23} &= \{ \bm{\vartheta}^I_{1}-\bm{\vartheta}^I_{2}, \bm{\vartheta}^I_{6}-\bm{\vartheta}^I_{5} \} \, , &
    	\tem_{24} &= \{ \bm{\vartheta}^I_{1}-\bm{\vartheta}^I_{3}, \bm{\vartheta}^I_{4}+\bm{\vartheta}^I_{6} \} \, , &
    	\tem_{34} &= \{ \bm{\vartheta}^I_{2}-\bm{\vartheta}^I_{3}, \bm{\vartheta}^I_{4}-\bm{\vartheta}^I_{5} \} \, , \notag \\  
    	\tem_{123} &= \{ \bm{\vartheta}^I_{1}-\bm{\vartheta}^I_{2}+\bm{\vartheta}^I_{4}\} \, , &
    	\tem_{124} &= \{ \bm{\vartheta}^I_{1}-\bm{\vartheta}^I_{3}+\bm{\vartheta}^I_{5}\} \, , &
    	\tem_{134} &= \{ \bm{\vartheta}^I_{2}-\bm{\vartheta}^I_{3}+\bm{\vartheta}^I_{6}\} \, , 
    	\notag \\  
    	\tem_{234} &= \{ \bm{\vartheta}^I_{4}-\bm{\vartheta}^I_{5}+\bm{\vartheta}^I_{6}\} \, , 
    \end{align}
which are based on the lowest order N\'ed\'elec base functions on the unit tetrahedron
\begin{align}
    	\bm{\vartheta}_1(\xi, \eta ,\zeta) &= \begin{bmatrix}
    		\zeta \\ \zeta \\ 1 - \xi - \eta
    	\end{bmatrix} \, , & 
        \bm{\vartheta}_2(\xi, \eta ,\zeta) &= \begin{bmatrix}
        	\eta \\ 1 - \xi - \zeta \\  \eta
        \end{bmatrix} \, , &
        \bm{\vartheta}_3(\xi, \eta ,\zeta) &= \begin{bmatrix}
        	1 - \eta - \zeta \\ \xi  \\ \xi
        \end{bmatrix} \, , \notag \\
        \bm{\vartheta}_4(\xi, \eta ,\zeta) &= \begin{bmatrix}
        	0 \\ \zeta \\ - \eta
        \end{bmatrix} \, , & 
        \bm{\vartheta}_5(\xi, \eta ,\zeta) &= \begin{bmatrix}
        	\zeta \\ 0 \\ - \xi
        \end{bmatrix} \, , &
        \bm{\vartheta}_6(\xi, \eta ,\zeta) &= \begin{bmatrix}
        	\eta \\ -\xi  \\ 0
        \end{bmatrix} \, .
    \label{eq:nedtet}
    \end{align}
    For the non-gradient cell functions we use the construction introduced in \cite{AINSWORTH2018178}
    \begin{align}
    	&\mathcal{R}^p = \left \{ (p+1) b_{i-{e_j}}^p \nabla \lambda_j - \dfrac{i_j}{p+1} \nabla_\xi b_i^{p+1} \; | \;  i \in \mathcal{I}_o \right \} \, ,  
    \end{align}
    where $\mathcal{I}_o$ is the set of multi-indices of cell functions, $e_j$ is the unit multi-index with the value one at position $j$ and $i_j$ is the value of the $i$-multi-index at position $j$. 
    Note that only the first term in the cell functions is required to span the next space in the sequence due to
    \begin{align}
    	\curl \left ([p+1] b_{i-{e_j}}^p \nabla_\xi \lambda_j - \dfrac{i_j}{p+1} \nabla_\xi b_i^{p+1} \right ) = \curl([p+1] b_{i-{e_j}}^p \nabla_\xi \lambda_j) \, .
    \end{align}
    However, without the added gradient the function would not belong to $[\Po^p]^3 \oplus \bm{\xi}\times [\widetilde{\Po}]^3$ and consequently, would not be part of the N\'ed\'elec space.
    By limiting $\mathcal{R}^p$ to $\mathcal{R}^p_*$ such that $\mathcal{R}_*^p$ contains only the surface permutations with $\nabla \lambda_j = \vb{e}_j$ and the cell permutations with $j \in \{1,2\}$, one retrieves the necessary base functions.   
    The sum of the lowest order N\'ed\'elec base functions, the template base functions, gradient base functions, and the non-gradient cell base functions yields exactly $(p+4)(p+3)(p+1)/2$, thus satisfying the required dimensionality of the N\'ed\'elec space. 
    The complete space reads
    \begin{align}
    	\Ned^p = \, \Ned^0  &\oplus \left\{ \bigoplus_{i \in \mathcal{I}} \nabla \edge^{p+1}_{i}   \right\} \oplus \left\{ \bigoplus_{j \in \mathcal{J}} \nabla \face^{p+1}_j \right\} \oplus  \nabla \cell^{p+1}_{1234} \oplus \left\{ \bigoplus_{k = 1}^3 \ver^p_k \otimes \tem_k  \right\}  \oplus \left\{ \bigoplus_{i \in \mathcal{I}} \edge_i^{p} \otimes \tem_i \right\} 
    	\notag \\
    	& \oplus \left \{ \bigoplus_{j \in \mathcal{J}} \face^{p}_j \otimes \tem_j \right \} \oplus \mathcal{R}_*^{p+1} \, ,  \qquad \begin{aligned}
    		\mathcal{I} &= \{ (1,2),(1,3),(1,4),(2,3),(2,4),(3,4) \} \, , \\
    		\mathcal{J} &= \{ (1,2,3),(1,2,4),(1,3,4),(2,3,4) \} 
    	\end{aligned} \, .
    \end{align}
     Here, the B\'ezier basis is used to construct the higher order N\'ed\'elec base functions of the first type.
    \begin{definition} [B\'ezier-N\'ed\'elec I tetrahedral basis]
    	
    	The base functions are defined on the reference tetrahedron:
    	\begin{itemize}
    		\item for the edges we use the lowest order base functions from \cref{eq:nedtet}. The remaining edge base functions are given by the gradients
    		\begin{align}
    			&e_{12}: &  
    			\bm{\vartheta}(\xi,\eta,\zeta) &=  \bm{\vartheta}_1^I \, ,  \notag \\
    			&&
    			\bm{\vartheta}(\xi,\eta,\zeta) &= \nabla_\xi b_{00k}^{p+1} \, , & & 0 < k < p + 1 \, , \notag \\
    			&e_{13}: & 
    			\bm{\vartheta}(\xi,\eta,\zeta) &= \bm{\vartheta}_2^I \, , \notag \\ &&
    			\bm{\vartheta}(\xi,\eta,\zeta) &= \nabla_\xi b_{0j0}^{p+1} \, , & & 0 < j < p + 1 \, , \notag \\
    			&e_{14}: & 
    			\bm{\vartheta}(\xi,\eta,\zeta) &= \bm{\vartheta}_3^I \, , \notag \\ &&
    			\bm{\vartheta}(\xi,\eta,\zeta) &= \nabla_\xi b_{i00}^{p+1} \, , & & 0 < i < p + 1 \, , \notag \\
    			&e_{23}: &
    			\bm{\vartheta}(\xi,\eta,\zeta) &= \bm{\vartheta}_4^I \, , \notag \\ &&
    			\bm{\vartheta}(\xi,\eta,\zeta) &= \nabla_\xi b_{0j,p+1-j}^{p+1} \, , & & 0 < j < p + 1 \, , \notag \\
    			&e_{24}: &
    			\bm{\vartheta}(\xi,\eta,\zeta) &= \bm{\vartheta}_5^I \, , \notag \\ &&
    			\bm{\vartheta}(\xi,\eta,\zeta) &= \nabla_\xi b_{i0,p+1-i}^{p+1} \, , & & 0 < i < p + 1 \, , \notag \\
    			&e_{34}: &
    			\bm{\vartheta}(\xi,\eta,\zeta) &= \bm{\vartheta}_6^I \, , \notag \\ &&
    			\bm{\vartheta}(\xi,\eta,\zeta) &= \nabla_\xi b_{00k}^{p+1} \, , & & 0 < i < p + 1 \, ;
    		\end{align}
    	    \item on faces we employ both template base functions and gradients
    	    \begin{align}
    	    	&f_{123}: & 
    	    	\bm{\vartheta}(\xi,\eta,\zeta) &= b_{000}^{p} \bm{\vartheta}_4^I \, , \notag \\ &&
    	    	\bm{\vartheta}(\xi,\eta,\zeta) &= -b_{00p}^{p} \bm{\vartheta}_2^I \, , \notag \\ &&
    	    	\bm{\vartheta}(\xi,\eta,\zeta) &= b_{00k}^{p} (\bm{\vartheta}_4^I - \bm{\vartheta}_2^I) \, , & & 0 < k < p \, , \notag \\ &&
    	    	\bm{\vartheta}(\xi,\eta,\zeta) &= b_{0j0}^{p} (\bm{\vartheta}_1^I + \bm{\vartheta}_4^I) \, , & & 0 < j < p \, , \notag \\ &&
    	    	\bm{\vartheta}(\xi,\eta,\zeta) &= b_{0j,p-j}^{p} (\bm{\vartheta}_1^I - \bm{\vartheta}_2^I) \, , & & 0 < j < p \, , \notag \\ &&
    	    	\bm{\vartheta}(\xi,\eta,\zeta) &= b_{0jk}^{p} (\bm{\vartheta}_1^I - \bm{\vartheta}_2^I + \bm{\vartheta}_4^I) \, , & & 0 < j < p \, , \quad 0 < k < p - j \, , \notag \\ &&
    	    	\bm{\vartheta}(\xi,\eta,\zeta) &= \nabla_\xi b_{0jk}^{p+1} \, , & & 0 < j < p + 1 \, , \quad 0 < k < p + 1 - j \, , \notag \\
    	    	&f_{124}: & 
    	    	\bm{\vartheta}(\xi,\eta,\zeta) &= b_{000}^{p} \bm{\vartheta}_5^I \, , \notag \\ &&
    	    	\bm{\vartheta}(\xi,\eta,\zeta) &= -b_{00p}^{p} \bm{\vartheta}_3^I \, , \notag \\ &&
    	    	\bm{\vartheta}(\xi,\eta,\zeta) &= b_{00k}^{p} (\bm{\vartheta}_5^I - \bm{\vartheta}_3^I) \, , & & 0 < k < p \, , \notag \\ &&
    	    	\bm{\vartheta}(\xi,\eta,\zeta) &= b_{i00}^{p} (\bm{\vartheta}_1^I + \bm{\vartheta}_5^I) \, , & & 0 < i < p  \, , \notag \\ &&
    	    	\bm{\vartheta}(\xi,\eta,\zeta) &= b_{i0,p-i}^{p} (\bm{\vartheta}_1^I - \bm{\vartheta}_3^I) \, , & & 0 < i < p \, , \notag \\ &&
    	    	\bm{\vartheta}(\xi,\eta,\zeta) &= b_{i0k}^{p} (\bm{\vartheta}_1^I - \bm{\vartheta}_3^I + \bm{\vartheta}_5^I) \, , & & 0 < i < p \, , \quad 0 < k < p - i \, , \notag \\ &&
    	    	\bm{\vartheta}(\xi,\eta,\zeta) &= \nabla_\xi b_{i0k}^{p+1} \, , & & 0 < i < p + 1 \, , \quad 0 < k < p + 1 - i \, , \notag \\
    	    	&f_{134}: & 
    	    	\bm{\vartheta}(\xi,\eta,\zeta) &= b_{000}^{p} \bm{\vartheta}_6^I \, , \notag \\ &&
    	    	\bm{\vartheta}(\xi,\eta,\zeta) &= -b_{0p0}^{p} \bm{\vartheta}_3^I \, , \notag \\ &&
    	    	\bm{\vartheta}(\xi,\eta,\zeta) &= b_{0j0}^{p} (\bm{\vartheta}_6^I - \bm{\vartheta}_3^I) \, , & & 0 < j < p \, , \notag \\ &&
    	    	\bm{\vartheta}(\xi,\eta,\zeta) &= b_{i00}^{p} (\bm{\vartheta}_2^I + \bm{\vartheta}_6^I) \, , & & 0 < i < p  \, , \notag \\ &&
    	    	\bm{\vartheta}(\xi,\eta,\zeta) &= b_{i,p-i,0}^{p} (\bm{\vartheta}_2^I - \bm{\vartheta}_3^I) \, , & & 0 < i < p \, , \notag \\ &&
    	    	\bm{\vartheta}(\xi,\eta,\zeta) &= b_{ij0}^{p} (\bm{\vartheta}_2^I - \bm{\vartheta}_3^I + \bm{\vartheta}_6^I) \, , & & 0 < i < p \, , \quad 0 < j < p - i \, , \notag \\ &&
    	    	\bm{\vartheta}(\xi,\eta,\zeta) &= \nabla_\xi b_{ij0}^{p+1} \, , & & 0 < i < p + 1 \, , \quad 0 < j < p + 1 - i \, , \notag \\
    	    	&f_{234}: & 
    	    	\bm{\vartheta}(\xi,\eta,\zeta) &= b_{00p}^{p} \bm{\vartheta}_6^I \, , \notag \\ &&
    	    	\bm{\vartheta}(\xi,\eta,\zeta) &= -b_{0p0}^{p} \bm{\vartheta}_5^I \, , \notag \\ &&
    	    	\bm{\vartheta}(\xi,\eta,\zeta) &= b_{0j,p-j}^{p} (\bm{\vartheta}_6^I - \bm{\vartheta}_5^I) \, , & & 0 < j < p \, , \notag \\ &&
    	    	\bm{\vartheta}(\xi,\eta,\zeta) &= b_{i0,p-i}^{p} (\bm{\vartheta}_4^I + \bm{\vartheta}_6^I) \, , & & 0 < i < p  \, , \notag \\ &&
    	    	\bm{\vartheta}(\xi,\eta,\zeta) &= b_{i,p-i,0}^{p} (\bm{\vartheta}_4^I - \bm{\vartheta}_5^I)\, ,  & & 0 < i < p \, , \notag \\ &&
    	    	\bm{\vartheta}(\xi,\eta,\zeta) &= b_{ij,p-i-j}^{p} (\bm{\vartheta}_4^I - \bm{\vartheta}_5^I + \bm{\vartheta}_6^I) \, , & & 0 < i < p \, , \quad 0 < j < p - i \, , \notag \\ &&
    	    	\bm{\vartheta}(\xi,\eta,\zeta) &= \nabla_\xi b_{ij,p-i.j}^{p+1} \, , & & 0 < i < p + 1 \, , \quad 0 < j < p + 1 - i \, ;
    	    \end{align}
            \item the cell base functions read
            \begin{align}
            	&c_{1234}: & 
            	\bm{\vartheta}(\xi,\eta,\zeta) &= (p+2)b_{i-1,jk}^{p+1} \vb{e}_1  - \dfrac{i}{p+2} \nabla_\xi b_{ijk}^{p+2} \, , & & \begin{aligned}
            		&0<i<p + 2 \, , \\[-1ex] &0<j<p + 2-i \, , \\[-1ex] &0<k<p+2-i-j 
            	\end{aligned} \, ,  \notag \\ &&
            \bm{\vartheta}(\xi,\eta,\zeta) &= (p+2)b_{i,j-1,k}^{p+1} \vb{e}_2  - \dfrac{j}{p+2} \nabla_\xi b_{ijk}^{p+2} \, , & & \begin{aligned}
            	&0<i<p + 2 \, , \\[-1ex] &0<j<p + 2-i \, , \\[-1ex] &0<k<p+2-i-j 
            \end{aligned} \, ,  \notag \\ &&
            \bm{\vartheta}(\xi,\eta,\zeta) &= (p+2)b_{ij0}^{p+1} \vb{e}_3  - \dfrac{1}{p+2} \nabla_\xi b_{ij1}^{p+2} \, , & & \begin{aligned}
            	&0<i<p + 2 \, , \\[-1ex] &0<j<p + 2-i 
            \end{aligned} \, ,  \notag \\ &&
            	\bm{\vartheta}(\xi,\eta,\zeta) &=  \nabla_\xi b_{ijk}^{p+1}\, , & & \begin{aligned}
            		 	&0<i<p + 1 \, , \\[-1ex] &0<j<p + 1-i \, , \\[-1ex] &0<k<p+1-i-j 
            	\end{aligned} \, .
            \end{align}
    	\end{itemize}
    \end{definition}

\section{Numerical quadrature}
    Although the base functions are expressed using  $(\alpha,\beta,\gamma)$ the domain is either the reference triangle or the reference tetrahedron, which require fewer quadrature points than their counterparts given by the Duffy transformation (quad or hexahedron). As such, we employ a mixture of the efficient quadrature points introduced in \cite{Dunavant,XIAO2010663,WITHERDEN20151232,Papanicolopulos2016EfficientCO,Jas} for triangles and tetrahedra, where we avoid quadrature schemes with points on the edges or faces of the reference domain due to the recursion formula of the Bernstein polynomials \cref{eq:rec}. 
    The quadrature points are mapped to their equivalent expression in $(\alpha,\beta,\gamma)$. 
    Consequently, the integration over the reference triangle or tetrahedron reads
    \begin{align}
    	\int_{\surf_e} f(x,y) \, \dd \surf &= \int_{\Gamma} (f \circ (\xi, \eta))(\alpha,\beta) \, |\det \bm{J}| \, \dd \Gamma \, , \notag \\
    	\int_{\body_e} f(x,y,z) \, \dd \body &= \int_{\Omega} (f \circ (\xi, \eta, \zeta))(\alpha,\beta,\gamma) \, |\det \bm{J}| \, \dd \Omega \, .
    \end{align}

    For the lower order elements we use the Lagrangian-N\'ed\'elec construction from \cite{Sky2021,SKY2022115298}.
\section{Boundary conditions}
The degrees of freedom in \cite{Demkowicz2000} commute between the continuous and discrete spaces. As such, they allow to exactly satisfy the consistent coupling condition \cite{dagostino2021consistent}. We note that the functionals can be viewed as a hierarchical system of Dirichlet boundary problems. 
	In the case of hierarchical base functions \cite{Zaglmayr2006}, they can be solved independently. However, here the boundary value of each polytope is required in advance due to the non-hierarchical nature of Bernstein polynomials.  
	In other words, one must first solve the problem for vertices, then for edges, afterwards for faces, and finally for the cell. In our case the degrees of freedom for the cell are irrelevant since a cell is never part of the boundary.

	\subsection{Boundary vertices}
	The finite element mesh identifies each vertex with a tuple of coordinates. It suffices to evaluate the displacement field at the vertex 
	\begin{align}
		u_i^d =  \widetilde{u} \at_{\vb{x}_i} \, .
	\end{align}
    If the field is vectorial, each component is evaluated at the designated vertex.
    The boundary conditions of the microdistortion field are associated with tangential projections and as such do not have vertex-type degrees of freedom. This is the case since a vertex does not define a unique tangential plane. 
    
    \subsection{Boundary edges}
    The edge functionals from \cite{Demkowicz2000} for the $\Hone$-conforming subspace
    \begin{align}
        &l_{ij}(u) = \int_{\curv_i}  \dfrac{\partial q_j}{\partial s} \dfrac{\partial u}{\partial \curv} \, \dd \curv \, , \quad q \in \Po^p(\curv) \, ,
    \end{align}
    can be reformulated for a reference edge on a unit domain $\alpha \in [0,1]$.
    We parametrize the edge via
    \begin{align}
    	\vb{x}(\alpha) = (1- \alpha) \vb{x}_1 + \alpha \vb{x}_2 \, .
    \end{align}
    \begin{figure}
    	\centering
    	\definecolor{asl}{rgb}{0.4980392156862745,0.,1.}
    	\definecolor{asb}{rgb}{0.,0.4,0.6}
    	\begin{tikzpicture}[line cap=round,line join=round,>=triangle 45,x=1.0cm,y=1.0cm]
    		\clip(-7,-0.5) rectangle (12,3.5);
    		\draw (-6,1) node[circle,fill=asb,inner sep=1.5pt] {};
    		\draw (-4,1) node[circle,fill=asb,inner sep=1.5pt] {};
    		\draw [-to,color=black,line width=1.pt] (-6,1) -- (-3,1);
    		\draw [color=asb,line width=1.pt] (-6,1) -- (-4,1);
    		\draw (-3,1) node[color=black,anchor=west] {$\alpha$};
    		\draw (-6,1) node[color=asb,anchor=north west] {$_{0}$};
    		\draw (-4,1) node[color=asb,anchor=north west] {$_{1}$};
    		
    		\draw (-2.2,2) node[color=black,anchor=south] {$\bm{\xi}:\alpha \to \Gamma$};
    		\draw [-Triangle,color=black,line width=1.pt] (-3,2) -- (-1,2);
    		
    		\draw (0,2) node[circle,fill=asb,inner sep=1.5pt] {};
    		\draw (2,0) node[circle,fill=asb,inner sep=1.5pt] {};
    		\draw [color=asb,line width=.6pt] (0,0) -- (0,2) -- (2,0) -- (0,0);
    		\draw (2,0) node[color=asb,anchor=south west] {$_{\bm{\xi}_2}$};
    		\draw (0,2) node[color=asb,anchor=north east] {$_{\bm{\xi}_1}$};
    		\draw (0.6,0.6) node[color=black] {$\Gamma$};
    		
    		\draw [-to,color=asl,line width=1.pt] (0.75,1.25) -- (1.25,0.75);
    		\draw (1.25,0.75) node[color=asl,anchor=west] {$\bm{\tau}$};
    		
    		\draw [-to,color=black,line width=1.pt] (0,0) -- (3,0);
    		\draw [-to,color=black,line width=1.pt] (0,0) -- (0,3);
    		\draw (3,0) node[color=black,anchor=west] {$\xi$};
    		\draw (0,3) node[color=black,anchor=south] {$\eta$};
    		
    		\draw (8,3) node[circle,fill=asb,inner sep=1.5pt] {};
    		\draw (9,1) node[circle,fill=asb,inner sep=1.5pt] {};
    		\draw (8,3) node[color=asb,anchor=south west] {$_{\vb{x}_2}$};
    		\draw (9,1) node[color=asb,anchor=north east] {$_{\vb{x}_1}$};
    		\draw [color=asb,line width=.6pt] (8,3) -- (9,1);
    		\draw (6.6,0.6) node[color=black] {$\surf$};
    		\draw [-to,color=asl,line width=1.pt] (8.75,1.5) -- (8.25,2.5);
    		\draw (8.3,2.5) node[color=asl,anchor=west] {$\vb{t}$};
    		
    		\draw [-to,color=black,line width=1.pt] (6,0) -- (8,0);
    		\draw [-to,color=black,line width=1.pt] (6,0) -- (6,2);
    		\draw (8,0) node[color=black,anchor=west] {$x$};
    		\draw (6,2) node[color=black,anchor=south] {$y$};
    		
    		\draw (3.8,2) node[color=black,anchor=south] {$\vb{x}:\Gamma \to \surf$};
    		\draw [-Triangle,color=black,line width=1.pt] (3,2) -- (5,2);
    	\end{tikzpicture}
    	\caption{Barycentric mapping of edges from the unit domain to the reference triangle and onto the physical domain.}
    	\label{fig:edgemap}
    \end{figure}
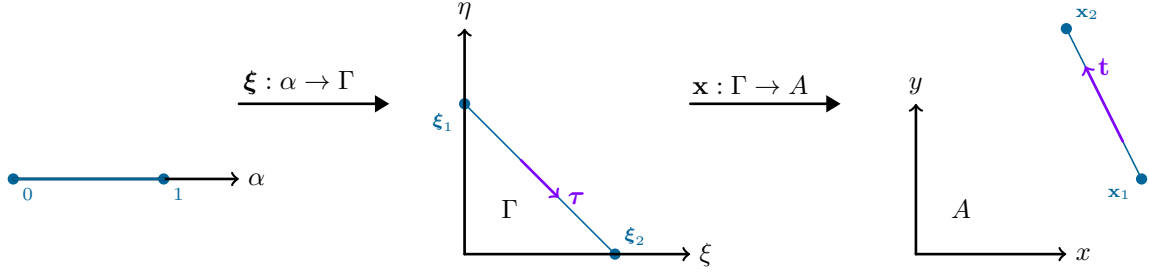
    As such, the following relation exists between the unit parameter and the arc-length parameter
    \begin{align}
    	 &\vb{t} = \dfrac{\dd}{\dd \alpha} \vb{x} = \vb{x}_2 - \vb{x}_1 \, , && \dd s =  \| \dd \vb{x} \| = \| \vb{x}_2 - \vb{x}_1 \| \dd \alpha = \|\vb{t} \| \dd \alpha \, .
    \end{align}
    By the chain rule we find
    \begin{align}
    	\dfrac{\dd u}{\dd s} = \dfrac{\dd u}{\dd \alpha} \dfrac{\dd \alpha}{\dd s} = \| \vb{t} \|^{-1} \dfrac{\dd u}{\dd \alpha} \, ,
    \end{align}
    for some function $u$. On edges, the test and trial functions are Bernstein polynomials parametrized by the unit domain. The function representing the boundary condition $\widetilde{u}(\vb{x})$ however, is parametrized by the Cartesian coordinates of the physical space. We find its derivative with respect to the arc-length parameter by observing
    \begin{align}
    	\dfrac{\dd}{\dd s} \widetilde{u} = \langle \dfrac{\dd}{\dd s} \vb{x} , \, \nabla_x \widetilde{u} \rangle \, .
    \end{align}
    The derivative of the coordinates with respect to the arc-length is simply the normed tangent vector
    \begin{align}
    	\dfrac{\dd}{\dd s} \vb{x} = \dfrac{\dd \vb{x}}{\dd \alpha} \dfrac{\dd \alpha}{\dd s} = \| \vb{t} \|^{-1} \vb{t} \, .
    \end{align}
    Consequently, the edge boundary condition is given by
    \begin{align}
    	\int_{\curv_i}  \dfrac{\partial q_j}{\partial s} \dfrac{\partial u}{\partial \curv} \, \dd \curv &= \int_0^1 \left( \|\vb{t}\|^{-1} \dfrac{\dd q_j}{\dd \alpha} \right) \left( \|\vb{t}\|^{-1} \dfrac{\dd u}{\dd \alpha} \right)  \|\vb{t} \| \, \dd \alpha \notag \\ &= \int_0^1 \left( \| \vb{t} \|^{-1} \dfrac{\dd q_j}{\dd \alpha} \right) \langle \| \vb{t} \|^{-1} \vb{t} , \,  \nabla_x \widetilde{u} \rangle  \| \vb{t} \| \, \dd \alpha
    	= \int_{\curv_i} \dfrac{\partial q_j}{\partial s} \dfrac{\partial \widetilde{u}}{\partial \curv}  \, \dd \curv  \qquad \forall \, q_j \in \Po^p(\alpha) \, ,
    \end{align}
    and can be solved by assembling the stiffness matrix of the edge and the load vector induced by the prescribed displacement field $\widetilde{u}$, representing volume forces
    \begin{align}
    	&k_{ij} = \int_0^1 \left( \|\vb{t}\|^{-1} \dfrac{\dd n_i}{\dd \alpha} \right) \left( \|\vb{t}\|^{-1} \dfrac{\dd n_j}{\dd \alpha} \right) \|\vb{t} \| \, \dd \alpha \, , && f_i = \int_0^1 \langle \| \vb{t} \|^{-1} \vb{t} , \,  \nabla_x \widetilde{u} \rangle \left( \| \vb{t} \|^{-1} \dfrac{\dd n_i}{\dd \alpha} \right) \| \vb{t} \| \, \dd \alpha \, .
    \end{align}
    
    Next we consider the Dirichlet boundary conditions for the microdistortion with the N\'ed\'elec space of the second type $\Nedtwo$. The problem reads
    \begin{align}
    	\int_{s_i} q_j \langle \vb{t} , \, \vb{p} \rangle \, \dd s = \int_{s_i} q_j \langle \vb{t} , \, \nabla_x \widetilde{u} \rangle \, \dd s \qquad \forall \, q_j \in \Po^p(s_i) \, .
    \end{align}
    Observe that on the edge the test functions $q_j$ are chosen to be the Bernstein polynomials. Further, by the polytopal template construction of the $\Nedtwo$-space there holds $\langle \vb{t} , \, \bm{\theta}_i \rangle |_s = n_i(\alpha)$. 
    Therefore, the components of the corresponding stiffness matrix and load vectors read
    \begin{align}
    	&k_{ij} = \int_0^1 n_i \, n_j \|\vb{t} \| \, \dd \alpha \, , && f_i = \int_0^1 n_i  \langle \vb{t} , \, \nabla_x \widetilde{u} \rangle \| \vb{t} \| \, \dd \alpha \, .
    \end{align}
    Note that in order to maintain the exactness property, the degree of the N\'ed\'elec spaces $\Ned^p , \Nedtwo^p$ is always one less than the degree of the subspace $\Ber^{p+1}$. 
    
    Lastly, we consider the N\'ed\'elec element of the first type. The problem is given by
    \begin{align}
    	\int_{s_i} q_j \langle \vb{t} , \, \vb{p} \rangle \, \dd s = \int_{s_i} q_j \langle \vb{t} , \, \nabla_x \widetilde{u} \rangle \, \dd s \qquad \forall \, q_j \in \Po^{p}(s_i) \, .
    \end{align}
   We define 
   \begin{align}
   	q_i = \dfrac{\dd }{\dd \alpha} n_i^{p+1} \, ,
   \end{align}
   and observe that on the edges the N\'ed\'elec base functions yield
   \begin{align}
   	   \langle \vb{t} , \, \bm{\theta}_j \rangle = \langle \vb{t} , \, \nabla_x n_j^{p+1} \rangle = \dfrac{\dd}{\dd \alpha} n_j^{p+1} \, .
   \end{align}
    Therefore, the components of the stiffness matrix and the load vector result in
    \begin{align}
    	&k_{ij} = \int_0^1 \dfrac{\dd n^{p+1}_i}{\dd \alpha} \, \dfrac{\dd n^{p+1}_j}{\dd \alpha} \|\vb{t} \| \, \dd \alpha \, , && f_i = \int_0^1 \dfrac{\dd n^{p+1}_i}{\dd \alpha}  \langle \vb{t} , \, \nabla_x \widetilde{u} \rangle \| \vb{t} \| \, \dd \alpha \, .
    \end{align}
    \subsection{Boundary faces}
    We start with the face boundary condition for the $\Hone$-conforming subspace. The problem reads
    \begin{align}
    	\int_{\surf_i} \langle \nabla_f q_j , \, \nabla_f u \rangle \, \dd \surf = \int_{\surf_i} \langle \nabla_f q_j , \, \nabla_f \widetilde{u} \rangle \, \dd \surf \qquad \forall \, q_j \in \Po^p(\surf_i) \, .
    \end{align}
    The surface is parameterized by the barycentric mapping from the unit triangle $\Gamma = \{ (\xi, \eta) \in [0,1]^2 \; | \; \xi + \eta \leq 1 \}$. The surface gradient is given by
    \begin{align}
    	\nabla_f \widetilde{u} = \nabla_x \widetilde{u} - \dfrac{1}{ \| \vb{n} \|^2} \langle \nabla_x \widetilde{u} , \, \vb{n} \rangle \vb{n} \, ,
    \end{align}
    where $\vb{n}$ is the surface normal. The surface gradient can also be expressed via
    \begin{align}
    	&\nabla_f u = \vb{e}^i \partial_i^x u = \vb{g}^\beta \partial_\beta^\xi u \, , && \beta \in \{1,2\} \, ,
    \end{align}
    where $\partial_\beta^x$ are partial derivates with respect to the physical coordinates, $\partial_\beta^\xi$ are partial derivatives with respect to the reference domain and $\vb{g}^\beta$ are the contravariant base vectors. The Einstein summation convention over corresponding indices is implied. 
    The covariant base vectors are given by
    \begin{align}
    	\vb{g}_\beta = \dfrac{\partial \vb{x}}{\partial \xi^\beta} \, .
    \end{align}
    One can find the contravariant vector orthogonal to the surface by
    \begin{align}
    	\vb{g}^3 = \vb{n} = \vb{g}_1 \times \vb{g}_2 \, .
    \end{align}
    We define the mixed transformation matrix
    \begin{align}
    	\bm{T} = \left[
    		\vb{g}_1 \, , \, \vb{g}_2 \, , \, \vb{g}^3
    	\right] \, .
    \end{align}
    Due to the orthogonality relation $\langle \vb{g}_i , \, \vb{g}^j \rangle = \delta_i^{\,j}$ the transposed inverse of $\bm{T}$ is clearly 
    \begin{align}
    	\bm{T}^{-T} = \left[
    		\vb{g}^1 \, , \, \vb{g}^2 \, , \, \vb{g}_3
    	\right] \, .
    \end{align}
    Thus, we can compute the surface gradient of functions parametrized by the reference triangle via
    \begin{align}
    	&\nabla_f u = \left[
    		\vb{g}^1 \, , \, \vb{g}^2 
    	\right] \nabla_\xi u = \bm{T}^{-T}_* \nabla_\xi u \, , && \bm{T}^{-T}_* = \left[
    	\vb{g}^1 \, , \, \vb{g}^2 
    \right] \, .
    \end{align}
    Further, there holds the following relation between the physical surface and the reference surface
    \begin{align}
    	\dd \surf = \| \vb{n} \| \dd \Gamma = \| \vb{g}^3 \| \dd \Gamma  = \sqrt{ \langle \vb{g}_1 \times \vb{g}_2 , \, \vb{g}^3 \rangle } \, \dd \Gamma = \sqrt{\det \bm{T}} \, \dd \Gamma \, .
    \end{align}
    Consequently, we can write the components of the stiffness matrix and load vector as
    \begin{align}
    	k_{ij} &= \int_\Gamma \langle \bm{T}_*^{-T} \nabla_\xi n_i , \, \bm{T}_*^{-T} \nabla_\xi n_j \rangle \sqrt{\det \bm{T}} \, \dd \Gamma \, , \notag \\ f_i &= \int_\Gamma \langle \bm{T}_*^{-T} \nabla_\xi n_i , \, \nabla_x \widetilde{u} - (\det \bm{T})^{-1}  \langle \nabla_x \widetilde{u} , \, \vb{n} \rangle \vb{n} \rangle \sqrt{\det \bm{T}} \, \dd \Gamma = \int_\Gamma \langle \bm{T}_*^{-T} \nabla_\xi n_i , \, \nabla_x \widetilde{u} \rangle \sqrt{\det \bm{T}} \, \dd \Gamma \, ,
    \end{align}
    with the orthogonality $\langle \vb{g}^\beta , \, \vb{n} \rangle = 0$ for $\beta \in \{1,2\}$.
    
    In order to embed the consistent coupling boundary condition to the microdistortion we deviate from the degrees of freedom defined in \cite{Demkowicz2000} and apply the simpler $\Hr{}$-projection
    \begin{align}
    	\langle \vb{q}_i , \, \vb{p} , \rangle_{\Hr{}} = \langle \vb{q}_i , \,  \nabla_f \widetilde{u} \rangle_{\Hr{}} \qquad \forall \, \vb{q}_i \in \Ned^p(\surf) \quad \text{or} \quad \forall \, \vb{q}_i \in \Nedtwo^p(\surf) \, .
    \end{align}
    Due to $\ker(\curl) = \nabla \Hone$ the problem reduces to
    \begin{align}
    	\int_{\surf_i} \langle \vb{q}_j , \, \vb{p} \rangle  + \langle \rot {\vb{q}_j} , \, \rot {\vb{p}} \rangle \, \dd \surf = \int_{\surf_i} \langle \vb{q}_j , \, \nabla_f \widetilde{u} \rangle \, \dd \surf \qquad \forall \, \vb{q}_j \in \Ned^p(\surf) \quad \text{or} \quad \forall \, \vb{q}_j \in \Nedtwo^p(\surf) \, .
    \end{align}
    We express the co- and contravariant Piola transformation from the two-dimensional reference domain to the three-dimensional physical domain using
    \begin{align}
    	&\bm{\theta}_i = \bm{T}_*^{-T} \bm{\vartheta}_i \, , &&
    	\di_x \bm{R} \, \bm{\theta}_i = \dfrac{1}{\sqrt{\det \bm{T}}} \di_\xi \bm{R} \, \bm{\vartheta}_i \, .
    \end{align}
    Thus, the stiffness matrix components and load vector components read
    \begin{align}
    	k_{ij} &= \int_\Gamma \langle \bm{T}_*^{-T} \bm{\vartheta}_i , \, \bm{T}_*^{-T} \bm{\vartheta}_j \rangle + \langle (\det \bm{T})^{-1/2} \di_\xi \bm{R} \, \bm{\vartheta}_i , \,  (\det \bm{T})^{-1/2} \di_\xi \bm{R} \, \bm{\vartheta}_j \rangle \sqrt{\det \bm{T}} \, \dd \Gamma \, , \notag \\
    	f_i &= \int_\Gamma \langle \bm{T}_*^{-T} \bm{\vartheta}_i , \, \nabla_x \widetilde{u} - (\det \bm{T})^{-1} \langle \nabla_x \widetilde{u} , \, \vb{n} \rangle \vb{n} \rangle \sqrt{\det \bm{T}} \, \dd \Gamma = \int_\Gamma \langle \bm{T}_*^{-T} \bm{\vartheta}_i , \, \nabla_x \widetilde{u} \rangle \sqrt{\det \bm{T}} \, \dd \Gamma \, ,
    \end{align} 
    where we again make use of the orthogonality between the surface tangent vectors and its normal vector.

\section{Numerical examples}
	In the following we test the finite element formulations with an artificial analytical solution in the antiplane shear model and with an analytical solution for an infinite plane under cylindrical bending in the three dimensional model. Finally, we benchmark the  ability of the finite element formulations to correctly interpolate between micro $\Cmic$ and macro $\Cmac$ stiffnesses as described by the characteristic length scale parameter $\Lc$.
	The majority of convergence results are presented by measuring the error in the Lebesgue norm over the domain
	\begin{align}
		&\| \widetilde{\ud} - \ud^h \|_{\Le} = \sqrt{\int_\body \| \widetilde{\ud} - \ud^h \|^2 \, \dd \body}  \, , && \| \widetilde{\Pm} - \Pm^h \|_{\Le} = \sqrt{\int_{\body} \| \widetilde{\Pm} - \Pm^h \|^2 \, \dd \body}  \, , 
	\end{align}
    in which context $\{\widetilde{\ud},\widetilde{\Pm}\}$ and $\{\ud^h,\Pm^h\}$ are the analytical and approximate subspace solutions, respectively.
    
\subsection{Compatible microdistortion} \label{ex:1}
\begin{figure}
    	\centering
    	\begin{subfigure}{0.3\linewidth}
    		\centering
    		\includegraphics[width=1\linewidth]{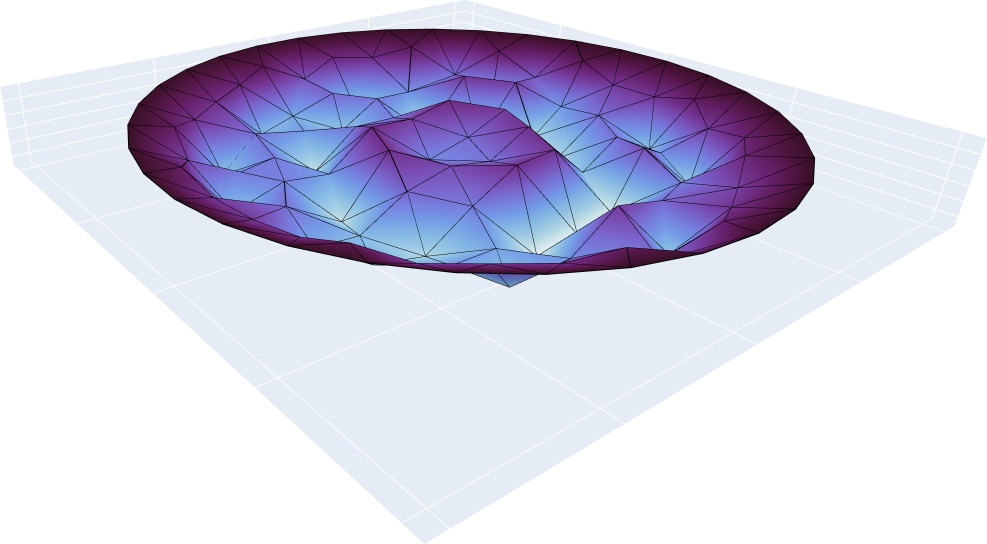}
    		\caption{}
    	\end{subfigure}
    	\begin{subfigure}{0.3\linewidth}
    		\centering
    		\includegraphics[width=1\linewidth]{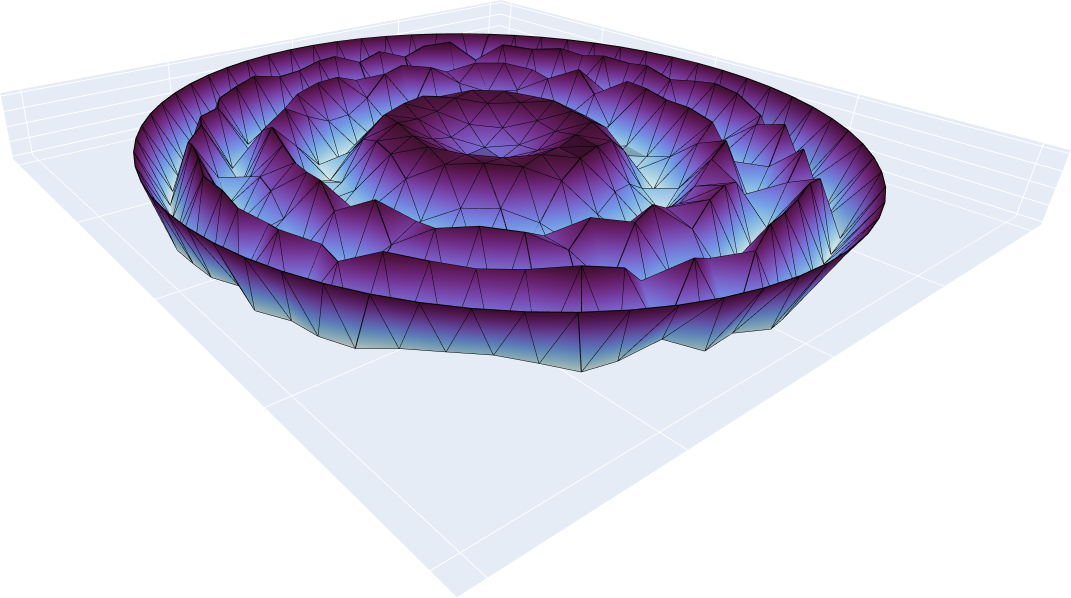}
    		\caption{}
    	\end{subfigure}
    	\begin{subfigure}{0.3\linewidth}
    		\centering
    		\includegraphics[width=1\linewidth]{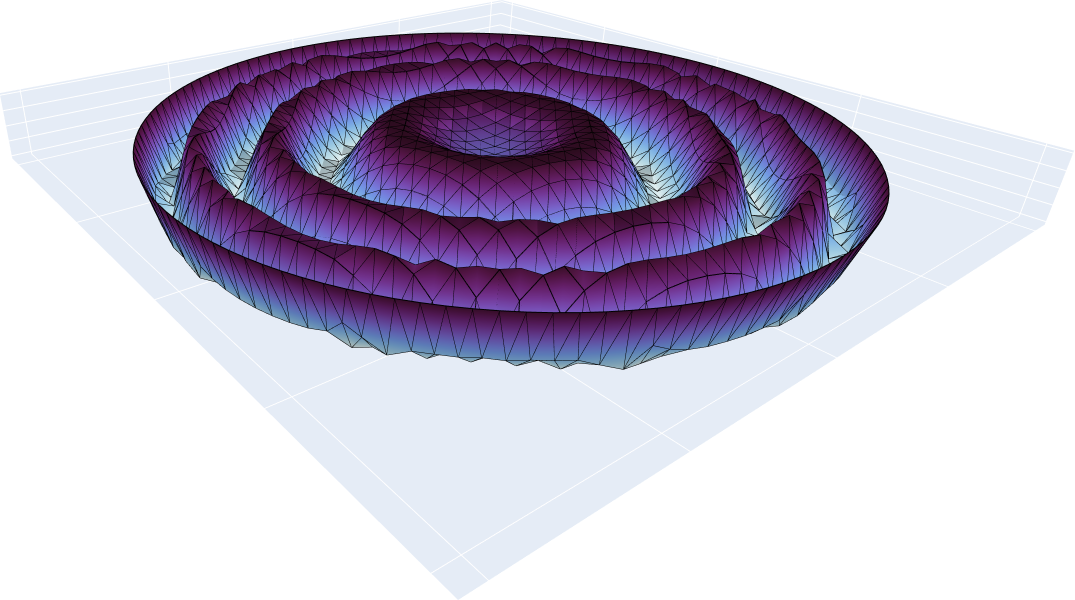}
    		\caption{}
    	\end{subfigure}
    	\begin{subfigure}{0.3\linewidth}
    		\centering
    		\includegraphics[width=1\linewidth]{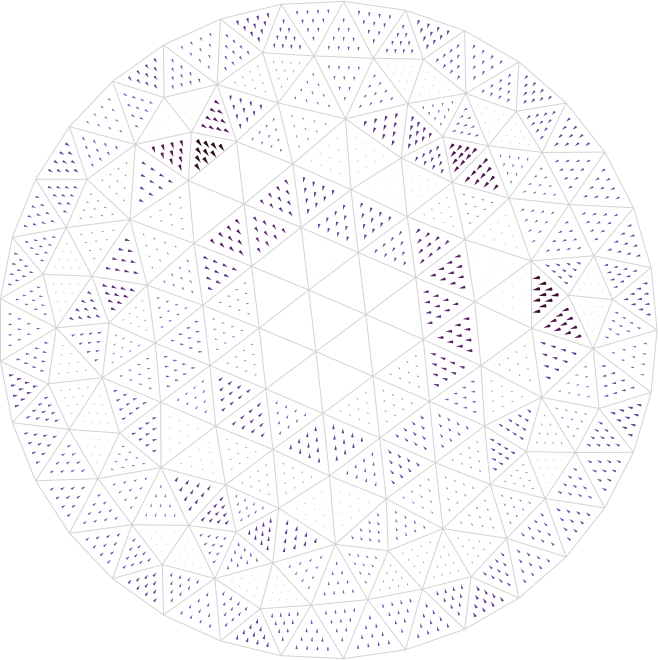}
    		\caption{}
    	\end{subfigure}
    	\begin{subfigure}{0.3\linewidth}
    		\centering
    		\includegraphics[width=1\linewidth]{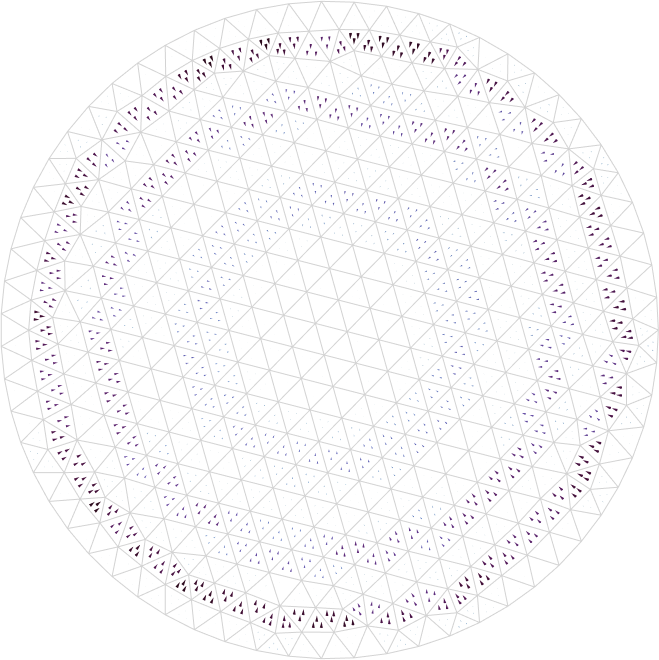}
    		\caption{}
    	\end{subfigure}
    	\begin{subfigure}{0.3\linewidth}
    		\centering
    		\includegraphics[width=1\linewidth]{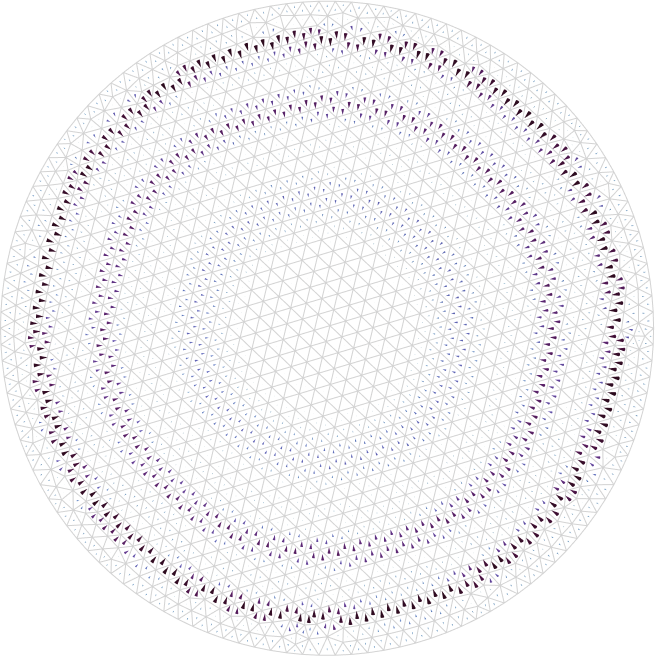}
    		\caption{}
    	\end{subfigure}
    	\begin{subfigure}{0.3\linewidth}
    		\centering
    		\includegraphics[width=1\linewidth]{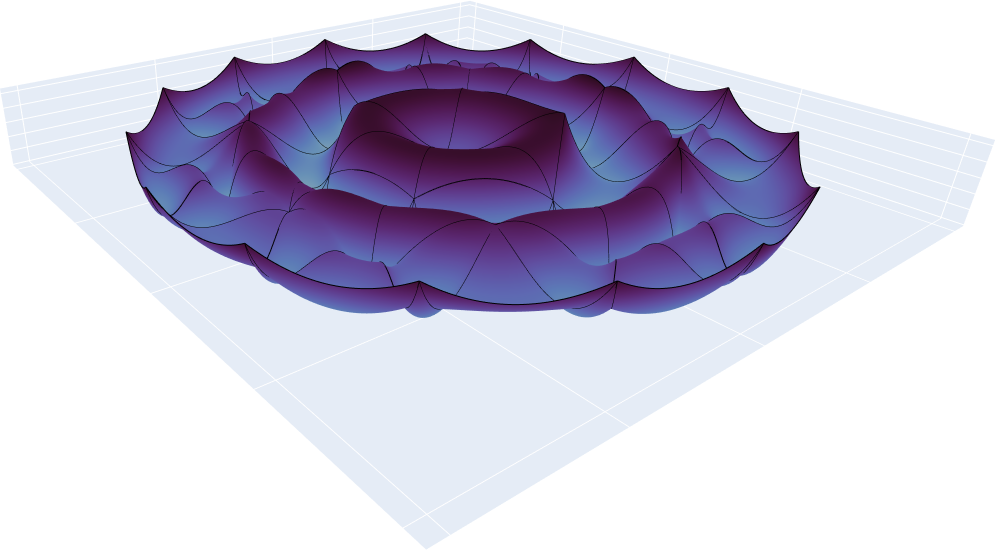}
    		\caption{}
    	\end{subfigure}
    	\begin{subfigure}{0.3\linewidth}
    		\centering
    		\includegraphics[width=1\linewidth]{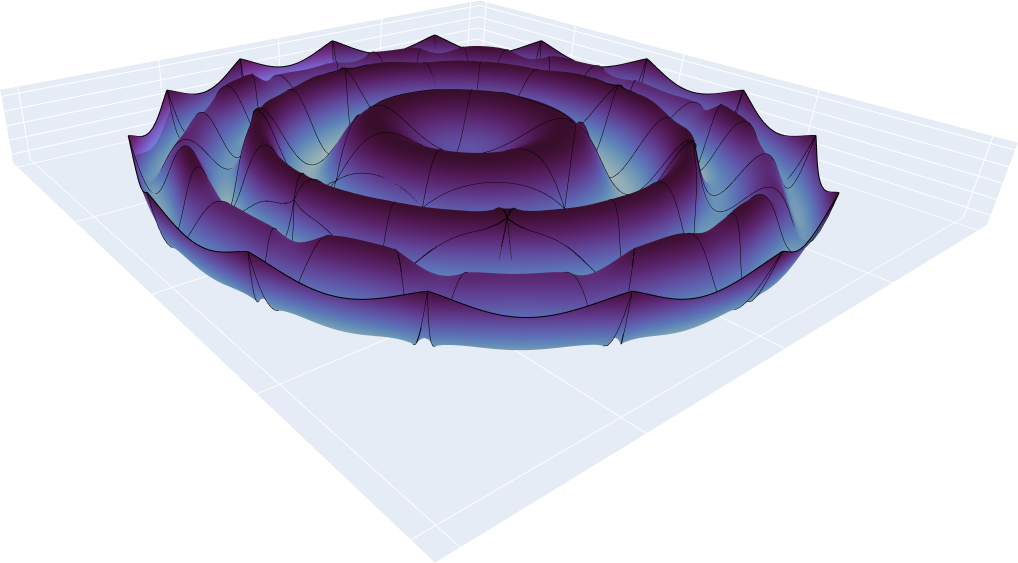}
    		\caption{}
    	\end{subfigure}
    	\begin{subfigure}{0.3\linewidth}
    		\centering
    		\includegraphics[width=1\linewidth]{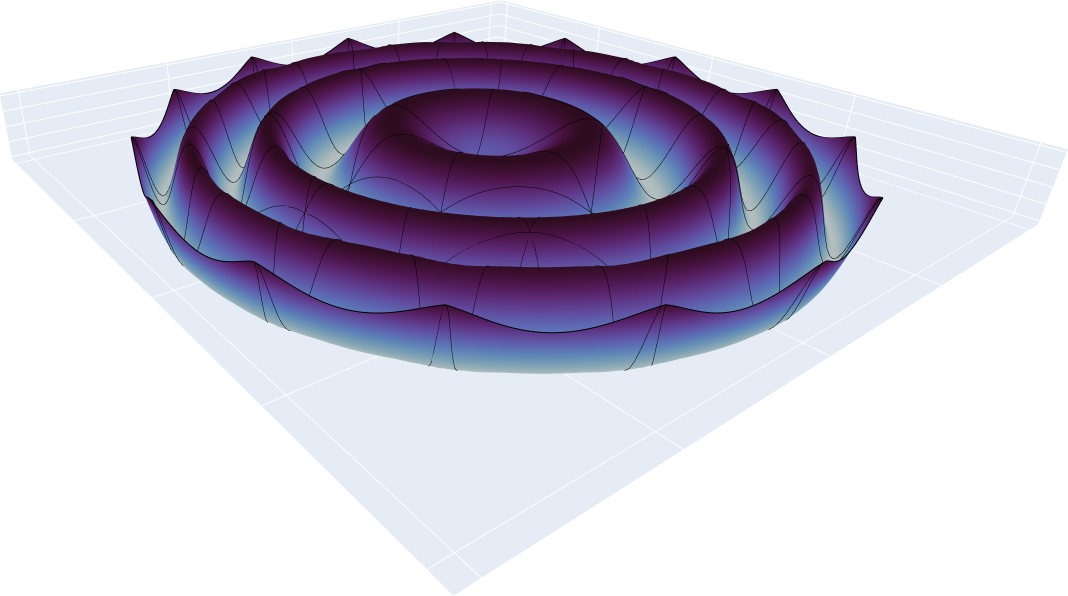}
    		\caption{}
    	\end{subfigure}
    	\begin{subfigure}{0.3\linewidth}
    		\centering
    		\includegraphics[width=1\linewidth]{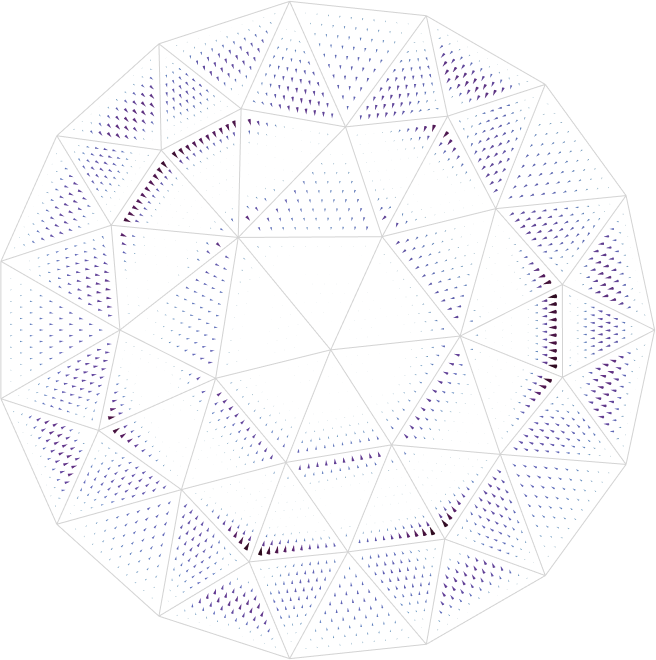}
    		\caption{}
    	\end{subfigure}
    	\begin{subfigure}{0.3\linewidth}
    		\centering
    		\includegraphics[width=1\linewidth]{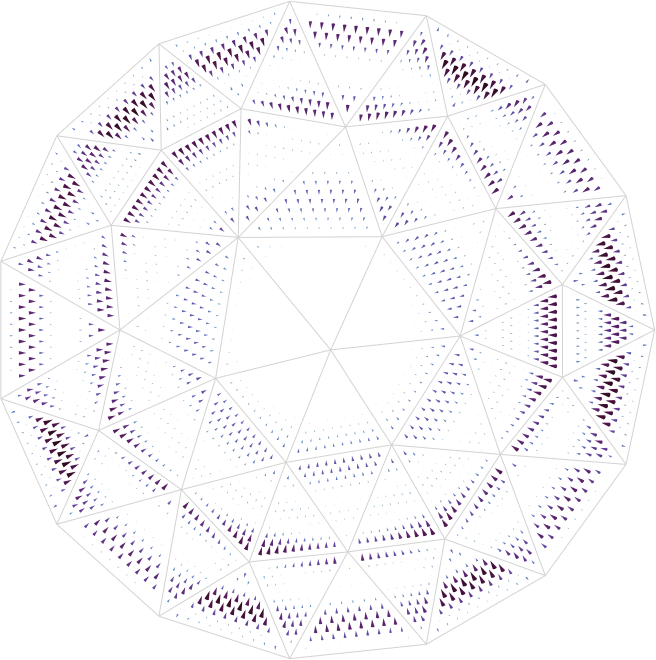}
    		\caption{}
    	\end{subfigure}
    	\begin{subfigure}{0.3\linewidth}
    		\centering
    		\includegraphics[width=1\linewidth]{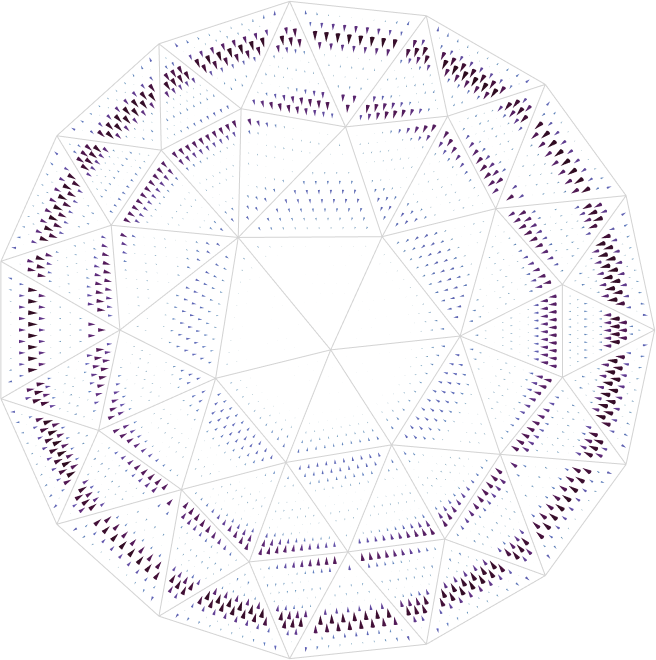}
    		\caption{}
    	\end{subfigure}
    	\caption{Depiction of the displacement field (a)-(c) and the microdistortion field (d)-(f) for the antiplane shear problem, for the linear element under h-refinement with $225$, $763$ and $2966$ elements, corresponding to $485$, $1591$ and $6060$ degrees of freedom. The p-refinement of the displacement field on the coarsest mesh of $57$ elements is visualized in (g)-(l) with $p \in \{3,5,7\}$, corresponding to $731$, $2072$ and $4097$ degrees of freedom.}
    	\label{fig:anti_p}
    \end{figure}

\begin{figure}
    	\centering
    	\begin{subfigure}{0.48\linewidth}
    		\centering
    		\begin{tikzpicture}
    			\definecolor{asl}{rgb}{0.4980392156862745,0.,1.}
    			\definecolor{asb}{rgb}{0.,0.4,0.6}
    			\begin{loglogaxis}[
    				/pgf/number format/1000 sep={},
    				axis lines = left,
    				xlabel={degrees of freedom},
    				ylabel={$\| \widetilde{u} - u^h \|_{\Le}$},
    				xmin=100, xmax=2e5,
    				ymin=1e-4, ymax=8000,
    				xtick={1e3,1e4,1e5},
    				ytick={1e-3,1e-1, 1e1},
    				legend style={at={(0.95,1)},anchor= north east},
    				ymajorgrids=true,
    				grid style=dotted,
    				]
    				\addplot[color=asl, mark=triangle] coordinates {
    					( 131 , 122.84271255729956 )
    					( 485 , 16.235821933959635 )
    					( 1591 , 5.033330206769331 )
    					( 6060 , 1.7411559824707534 )
    					( 37027 , 0.302615525389086 )
    				};
    				\addlegendentry{$\Lag^1 \times \Ned^0$}
    				
    				\addplot[color=violet, mark=o] coordinates {
    					( 317 , 21.990033994688236 )
    					( 1193 , 4.023797553797182 )
    					( 3943 , 1.268550265629483 )
    					( 15084 , 0.18274585290293446 )
    					( 92407 , 0.011528758508281563 )
    				};
    				\addlegendentry{$\Lag^2 \times \Nedtwo^1$}
    				
    				\addplot[color=asb, mark=square] coordinates {
    					( 731 , 9.15808287944714 )
    					( 2801 , 2.120283218033357 )
    					( 9347 , 0.18617381290244053 )
    					( 35972 , 0.012975143789274174 )
    					( 221207 , 0.00035331921471973813 )
    				};
    				\addlegendentry{$\Ber^3 \times \Nedtwo^2$}
    				
    				\addplot[color=blue, mark=pentagon] coordinates {
    					( 2072 , 3.2343114039220615 )
    					( 8042 , 0.13578857747499423 )
    					( 27022 , 0.004910428928331921 )
    					( 104442 , 6.966539793078105e-05 )
    					( 644002 , 3.0795786015737187e-07 )
    				};
    				\addlegendentry{$\Ber^5 \times \Nedtwo^{4}$}
    				
    				\addplot[color=cyan, mark=diamond] coordinates {
    					( 4097 , 0.4361035603837418 )
    					( 15983 , 0.0035940039219269438 )
    					( 53853 , 7.49190855210908e-05 )
    					( 208504 , 2.918199092694646e-07 )
    					( 1287057 , 4.3698070202838184e-10 )
    				};
    				\addlegendentry{$\Ber^{7} \times \Nedtwo^{6}$}
    				
    				\addplot[dashed,color=black, mark=none]
    				coordinates {
    					(0.5e3, 4e1)
    					(0.7e4, 2.8571428571428577)
    				};
    			
    			    \addplot[dashed,color=black, mark=none]
    			    coordinates {
    			    	(7e3, 0.2e-1)
    			    	(2e4, 0.0005073038414007921)
    			    };
    				
    			\end{loglogaxis}
    			\draw (2.7,3.7) 
    			node[anchor=south]{$\mathcal{O}(h^{2})$};
    			\draw (3.7,0.7) 
    			node[anchor=south]{$\mathcal{O}(h^{7})$};
    		\end{tikzpicture}
    		\caption{}
    	\end{subfigure}
    	\begin{subfigure}{0.48\linewidth}
    		\centering
    		\begin{tikzpicture}
    			\definecolor{asl}{rgb}{0.4980392156862745,0.,1.}
    			\definecolor{asb}{rgb}{0.,0.4,0.6}
    			\begin{loglogaxis}[
    				/pgf/number format/1000 sep={},
    				axis lines = left,
    				xlabel={degrees of freedom},
    				ylabel={$\|\widetilde{\vb{p}} - \vb{p}^h \|_{\Le}$},
    				xmin=100, xmax=2e5,
    				ymin=1e-4, ymax=1e+5,
    				xtick={1e3,1e4,1e5},
    				ytick={1e-3,1e-1, 1e1},
    				legend pos= north east,
    				ymajorgrids=true,
    				grid style=dotted,
    				]
    				\addplot[color=asl, mark=triangle] coordinates {
    					( 131 , 26.494279284901452 )
    					( 485 , 16.20102774183262 )
    					( 1591 , 10.395395643037588 )
    					( 6060 , 6.545273726959493 )
    					( 37027 , 2.7774965076073888 )
    				};
    				\addlegendentry{$\Lag^1 \times \Ned^0$}
    				
    				\addplot[color=violet, mark=o] coordinates {
    					( 317 , 18.367819092954853 )
    					( 1193 , 9.572242787607006 )
    					( 3943 , 5.033331448037817 )
    					( 15084 , 1.3963498316546792 )
    					( 92407 , 0.2239802649291176 )
    				};
    				\addlegendentry{$\Lag^2 \times \Nedtwo^1$}
    				
    				\addplot[color=asb, mark=square] coordinates {
    					( 731 , 14.725863637303336 )
    					( 2801 , 6.329300571977952 )
    					( 9347 , 1.0155148414436967 )
    					( 35972 , 0.1582418040304384 )
    					( 221207 , 0.017982726335169005 )
    				};
    				\addlegendentry{$\Ber^3 \times \Nedtwo^2$}

    				\addplot[color=blue, mark=pentagon] coordinates {
    					( 2072 , 7.688874882688885 )
    					( 8042 , 0.7757692365625642 )
    					( 27022 , 0.15143729234143524 )
    					( 104442 , 0.050847725824359515 )
    					( 644002 , 0.012879564773535995 )
    				};
    				\addlegendentry{$\Ber^5 \times \Nedtwo^{4}$}
    				
    				\addplot[color=cyan, mark=diamond] coordinates {
    					( 4097 , 2.0198836592302793 )
    					( 15983 , 0.3754862148775511 )
    					( 53853 , 0.1414468939416425 )
    					( 208504 , 0.050070402364277296 )
    					( 1287057 , 0.012686584197861417 )
    				};
    				\addlegendentry{$\Ber^{7} \times \Nedtwo^{6}$}
    				
    				\addplot[dashed,color=black, mark=none]
    				coordinates {
    					(5e2, 7e1)
    					(0.7e4, 18.708286933869708)
    				};
    			
    			    \addplot[dashed,color=black, mark=none]
    			    coordinates {
    			    	(2e3, 1e-1)
    			    	(5e4, 0.004)
    			    };
    				
    			\end{loglogaxis}
    			\draw (2.6,3.5) node[anchor=south]{$\mathcal{O}(h)$};
    			\draw (4.25, 1.35) node[anchor=north]{$\mathcal{O}(h^2)$};
    		\end{tikzpicture}
    		\caption{}
    	\end{subfigure}
    	\caption{Convergence of displacement (a) and the microdistortion (b) under h-refinement for multiple polynomial degrees for the antiplane shear problem.}
    	\label{fig:gradient_case}
    \end{figure}
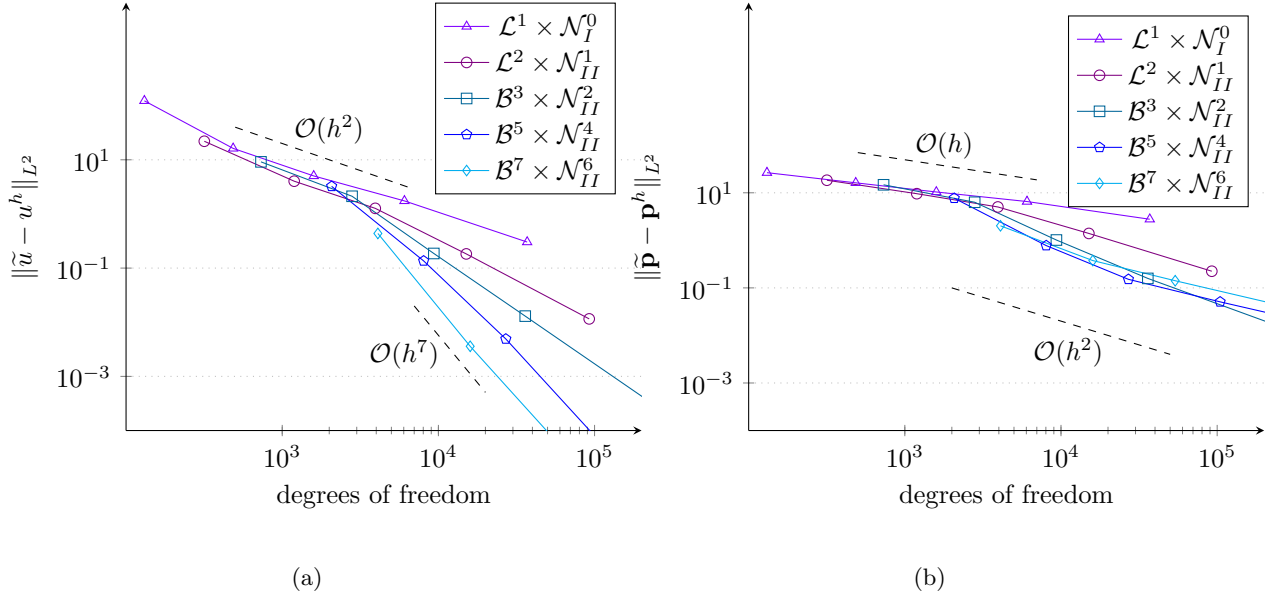

    In \cite{Sky2021} we explored the conditions for which the microdistortion $\vb{p}$ reduces to a gradient field, i.e. $\vb{p}$ is compatible. By defining the micro-moment with a scalar potential  
    \begin{align}
    	\vb{m} = \nabla \dfrac{100 - x^2 - y^2}{10} = -\dfrac{1}{5}\begin{bmatrix}
    		x \\ y 
    	\end{bmatrix} \, ,
    \end{align}
    and constructing an analytical solution for the displacement field
    \begin{align}
    	\widetilde{u} = \sin{\left(\dfrac{x^2 + y^2}{5} \right)} \, ,
    \end{align}
    we can recover the analytical solution of the microdistortion 
    \begin{align}
    	\vb{p} = \dfrac{1}{\mue + \mumi} (\vb{m} + \mue \nabla \widetilde{u}) &= \dfrac{1}{2}\left( -\dfrac{1}{5} \begin{bmatrix}
    		x \\ y
    	\end{bmatrix} + \dfrac{2}{5} \begin{bmatrix}  x \cos ([x^{2} + y^{2} ] / 5)\\ y \cos ([x^{2} + y^{2} ] / 5) \end{bmatrix}  \right)  = \dfrac{1}{5} \begin{bmatrix}  x \cos ([x^{2} + y^{2} ] / 5)\\ y \cos ([x^{2} + y^{2} ] / 5) \end{bmatrix} -\dfrac{1}{10} \begin{bmatrix}
    	x \\ y
    \end{bmatrix} \, ,
    \end{align}
    where for simplicity we set all material constants to one. Since $\vb{m}$ is a gradient field, the microdistortion $\vb{p}$ is also reduced to a gradient field and $\rot{\vb{p}} = 0$, see \cite{Sky2021}. Note that this result is specific to antiplane shear and does not generalize to the full three-dimensional model, compare \cite{SKY2022115298}. 
    We note that the microdistortion is not equal to the gradient of the displacement field
    and as such, their tangential projections on an arbitrary boundary are not automatically the same. However, for both the gradient of the displacement field and the micro-moment is the tangential projection on the boundary of the circular domain 
    $\overline{A} = \{\vb{x} \in \R^2 \; | \; \| \vb{x} \| \leq 10\}$ equal to zero 
    \begin{align}
    	\langle \nabla  \vb{t} , \, \widetilde{u}  \rangle \at_{\partial \surf} = \langle \vb{t} , \, \vb{m} \rangle \at_{\partial \surf} = 0 \, ,
    \end{align}
    and as such the microdistortion belongs to $\vb{p} \in \Hcz{,\surf}$. Consequently, we can set $\curv_D = \partial \surf$ and the consistent coupling condition remains compatible.
    With the displacement and the microdistortion fields at hand we derive the corresponding forces
    \begin{align}
    	f = \dfrac{1}{25} \left [2 x^{2} \sin{\left(\dfrac{x^2 + y^2}{5} \right)} + 2 y^{2} \sin{\left(\dfrac{x^2 + y^2}{5} \right)} - 10 \cos{\left(\dfrac{x^2 + y^2}{5} \right)} - 5 \right ] \, .
    \end{align}
    The approximation of the displacement and microdistortion fields using linear and higher order elements is shown in \cref{fig:anti_p}. We note that even with almost 3000 finite elements and 6000 degrees of freedom the linear formulation is incapable of finding an adequate approximation. On the other side of the spectrum, the higher order approximation (degree 7) with 57 elements and 4097 degrees of freedom yields very accurate results in the interior of the domain. However, the exterior of the domain is captured rather poorly. This is the case since the geometry of the circular domain is being approximated by linear triangles. Thus, in this setting, a finer mesh captures the geometry in a more precise manner.
    The effects of the geometry on the approximation of the solution are also clearly visible in the convergence graphs in \cref{fig:gradient_case}; only after a certain accuracy in the domain description is achieved do the finite elements retrieve their predicted convergence rates, compare \cite{SKY2022115298,Sky2021}. This is clearly observable when comparing the convergence curves of the linear and seventh order elements.
    The linear element generates quadratic convergence $p+1 = 1+1 =2$, whereas the seventh-order element yields the convergence slope $7$ (where $8$ is expected). Although the seventh-order formulation encompasses more degrees of freedom, it employs a coarser mesh and as such, generates higher errors at the boundary. The errors themselves can be traced back to the consistent coupling condition since, for a non-perfect circle the gradient of the displacement field induces tangential projections on the imperfect boundary. The influence of the latter effect is even more apparent in the convergence of the microdistortion, where the higher order formulations are unable to perform optimally on coarse meshes.
     
\subsection{Cylindrical bending} \label{ex:2}
	In order to test the capability of the finite element formulations to capture the intrinsic behaviour of the relaxed micromorphic model, we compare with analytical solutions of boundary-value problems.
	The first example considers the displacement and microdistortion fields under cylindrical bending \cite{Rizzi_bending} for infinitely extended plates. 
	Let the plates be defined as $\body = (-\infty,\infty)^2 \times [-1/2, 1/2]$, than the analytical solution for cylindrical bending reads
	\begin{align}
		&\vb{u} = \kappa\begin{bmatrix}
			- x z \\
			0 \\
			x^2 / 2
		\end{bmatrix} \, , &&
		\bm{P} = -\kappa\begin{bmatrix}
			[41 z + 20 \sqrt{82}\, \mathrm{sech}(\sqrt{41/2})\sinh(\sqrt{82}z)]/1681 & 0 & x \\
			0 & 0 & 0 \\
			-x & 0 & 0
		\end{bmatrix} \, , 
	\end{align}
    where $\mathrm{sech}(x)=1/\cosh(x)$, and for the following values of material constants
    \begin{align}
    	&\lame=\lammi = 0 \, , &&\mue=\muma = 1/2 \, , && \muc = 0 \, ,  && \Lc = 1  \, , &&\mumi = 20 \, .
    \end{align}
    The intensity of the curvature parameter $\kappa$ of the plate is chosen to be $\kappa = 14 / 200$.
    \begin{remark}
        The particular case of the cylindrical bending for which $\lame=\lammi=0$ (equivalent to a zero micro-Poisson's ratio) has been solved, along with its more general case ($\lame\neq\lammi\neq0$), in \cite{Rizzi_bending}.
        The advantage of considering this particular case is that a cut out finite plate of the infinite domain automatically exhibits the consistent coupling boundary conditions on its side surfaces.
    \end{remark}
    \begin{remark}
        Note that the general analytical solution for cylindrical bending does not depend on $\muc$, so we can set $\muc=0$ without loss of generality, compare \cite{Rizzi_bending}.
    \end{remark}
    We define the finite domain $\overline{\body} = [-10, 10]^2 \times [-1/2, 1/2]$ and the boundaries
    \begin{align}
    	\overline{\surf}_{D_1} &= \{-10\} \times [-10,10] \times  [-1/2, 1/2] \, , 
    && \overline{\surf}_{D_2} = \{10\} \times [-10,10] \times [-1/2, 1/2]  \, , \notag \\  \surf_{N} &= \partial \body \setminus \{ \overline{\surf}_{D_1} \oplus \overline{\surf}_{D_2} \} \, .
    \end{align}
    Additionally, on the Dirichlet boundary we impose the translated analytical solution 
    $\widetilde{\vb{u}} = \vb{u} - \begin{bmatrix} 0 & 0 & 3.5 \end{bmatrix}^T$. 
    
\begin{figure}
    	\centering
    	\begin{subfigure}{0.45\linewidth}
    		\centering
    		\includegraphics[width=0.7\linewidth]{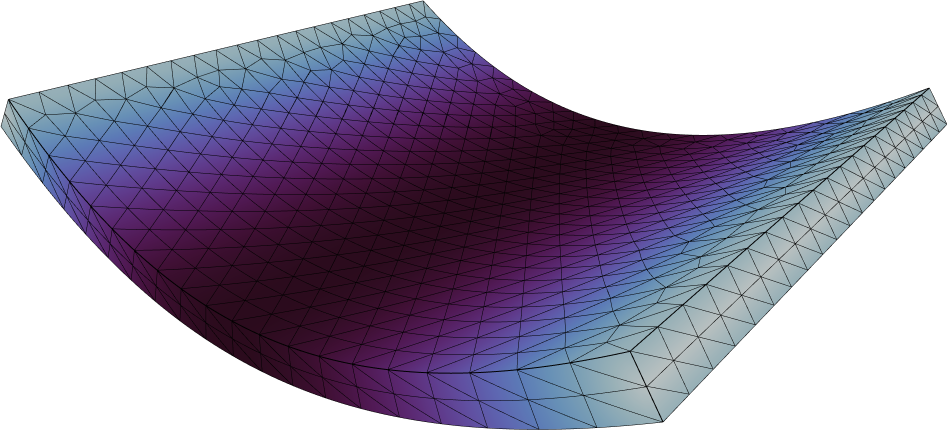}
    		\caption{}
    	\end{subfigure}
    	\begin{subfigure}{0.45\linewidth}
    		\centering
    		\includegraphics[width=0.7\linewidth]{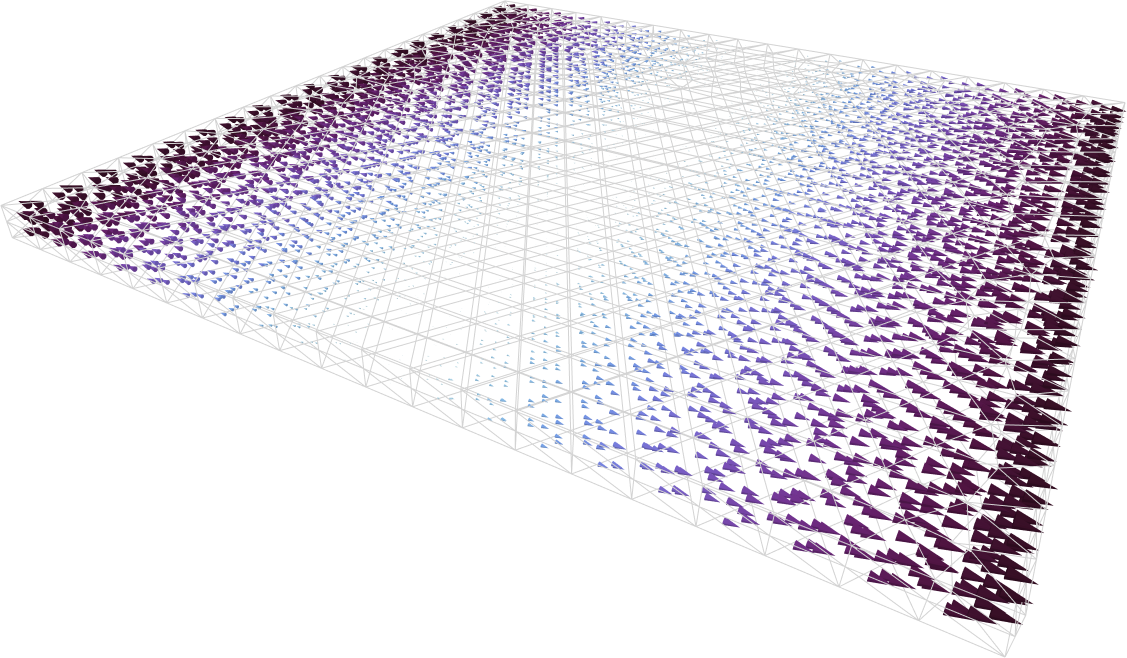}
    		\caption{}
    	\end{subfigure}
    	\caption{Displacement (a) and last row of the microdistortion (b) for the quadratic formulation using the N\'ed\'elec element of the first type.}
    	\label{fig:bending}
    \end{figure}

\begin{figure}
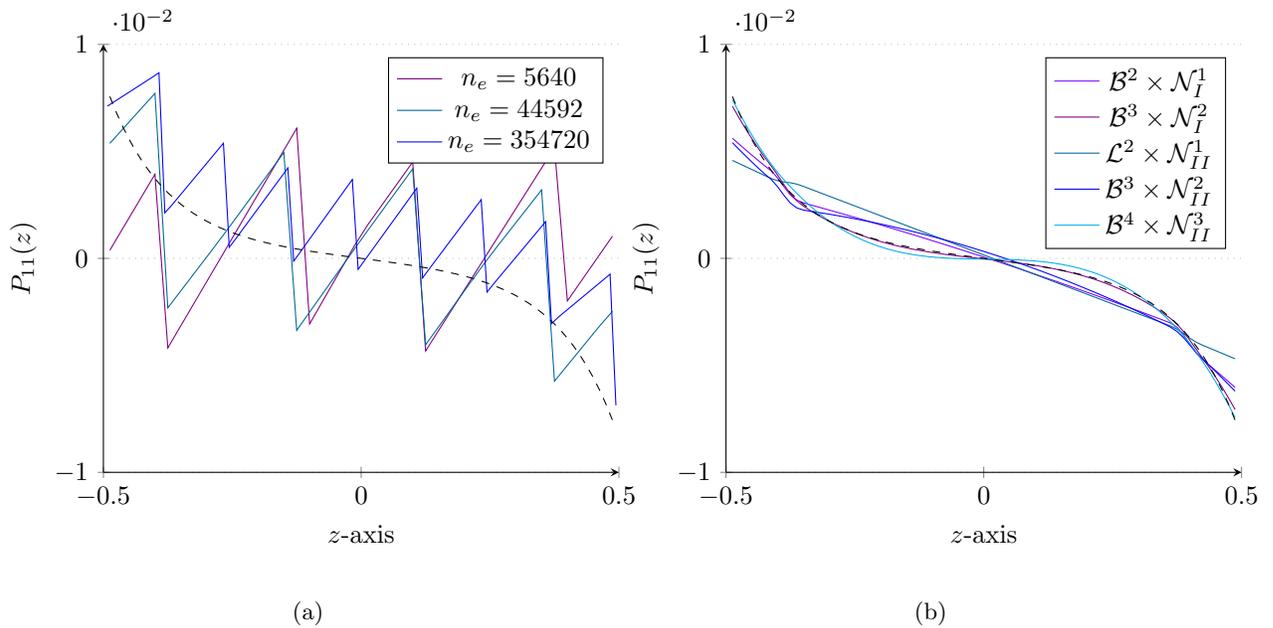

    	\centering
    	\begin{subfigure}{0.48\linewidth}
    		\centering
    		\input{graphs/bending}
    		\caption{}
    	\end{subfigure}
    	\begin{subfigure}{0.48\linewidth}
    		\centering
    		\input{graphs/bending2}
    		\caption{}
    	\end{subfigure}
    	\caption{Convergence of the lowest order formulation under h-refinement with $732$, $5640$ and $44592$ elements (a) and of the higher order formulations under p-refinement using $732$ elements(b) towards the analytical solution (dashed curve) of the $P_{11}(z)$ component at $x = y = 0$.}
    	\label{fig:bending2}
    \end{figure}

    The displacement field and the last row of the microdistortion are depicted in \cref{fig:bending}. The displacement field is dominated by its quadratic term and captured correctly. The last row of the microdistortion is a linear function and easily approximated even with linear elements.
    On the contrary, the $P_{11}$ component of the microdistortion is a hyperbolic function of the $z$-axis. The results of its approximation at $x = y = 0$ (the centre of the plane) are given in \cref{fig:bending2}.
    We observe that even increasing the number of linear finite elements to the extreme only results in better oscillations around the analytical solution. In comparison, higher order formulations converge towards the expected hyperbolic behaviour. The approximation of the quadratic N\'ed\'elec element of the first type is nearly perfect, whereas its second type counterpart clearly deviates from the analytical solution at $z \approx -0.25$. Taking the cubic second type element yields the desired result. This phenomenon is an evident indicator of the prominent role of the Curl of the microdistortion in this type of problems. Firstly, the microdistortion is a non-gradient field. Secondly, the Curl of the analytical solution induces an hyperbolic sine term. Such functions are often approximated using at least cubic terms in power series, thus explaining the necessity of such high order elements for correct computations.
    
\subsection{Bounded stiffness property} \label{ex:3}
    The characteristic length scale parameter $\Lc$ allows the relaxed micromorphic model to capture the transition from highly homogeneous materials to materials with a pronounced micro-structure by governing the influence of the micro-structure on the overall behaviour of the model. 
    We demonstrate this property of the model with an example, where we vary $\Lc$ and measure the resulting energy. 
	
	Let the domain be given by the axis-symmetric cube $\overline{\body} = [-1,1]^3$ with a total Dirichlet boundary 
	\begin{align}
		&\overline{\surf}_{D_1} = \{ (x,y,z) \in [-1,1]^3 \; | \; x = \pm 1 \} \, , && \overline{\surf}_{D_2} = \{ (x,y,z) \in [-1,1]^3 \; | \; y = \pm 1 \} \, ,
    \notag  \\ & \overline{\surf}_{D_3} = \{ (x,y,z) \in [-1,1]^3 \; | \; z = \pm 1 \} \, ,
	\end{align}
    we embed the periodic boundary conditions 
    \begin{align}
    	&\widetilde{u}\at_{\surf_{D_1}} = \begin{bmatrix}
    		(1-y^2) \sin(\pi[1-z^2]) / 10 \\ 0 \\ 0
    	\end{bmatrix} \, , && \widetilde{u}\at_{\surf_{D_2}} = \begin{bmatrix}
    	0 \\ (1-x^2) \sin(\pi[1-z^2]) / 10 \\ 0
    \end{bmatrix} \, , 
     \notag \\  & \widetilde{u}\at_{\surf_{D_3}} = \begin{bmatrix}
    	0 \\ 0  \\ (1-y^2) \sin(\pi[1-x^2]) / 10
    \end{bmatrix} \, .
    \end{align}
    The material parameters are chosen as
    \begin{align}
    	&\lamma = 2 \, , && \muma = 1 \, , && \lammi = 10 \, , && \mumi = 5 \, , && \muc = 1 \, ,
    \end{align}
    thus giving rise to the following meso-parameters via \cref{eq:muma}
    \begin{align}
    	&\lame = 2.5 \, , && \mue = 1.25 \, . 
    \end{align}
    
\begin{figure}
    	\centering
    	\begin{subfigure}{0.3\linewidth}
    		\centering
    		\includegraphics[width=1\linewidth]{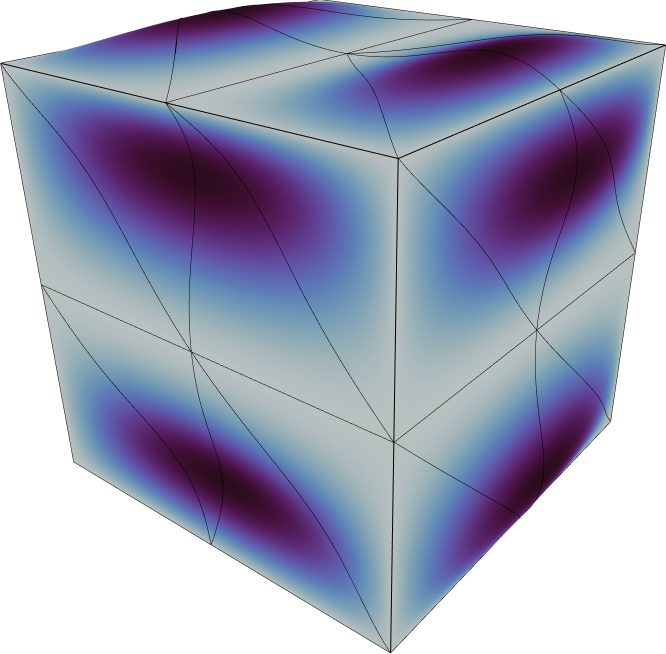}
    		\caption{}
    	\end{subfigure}
    	\begin{subfigure}{0.3\linewidth}
    		\centering
    		\includegraphics[width=1\linewidth]{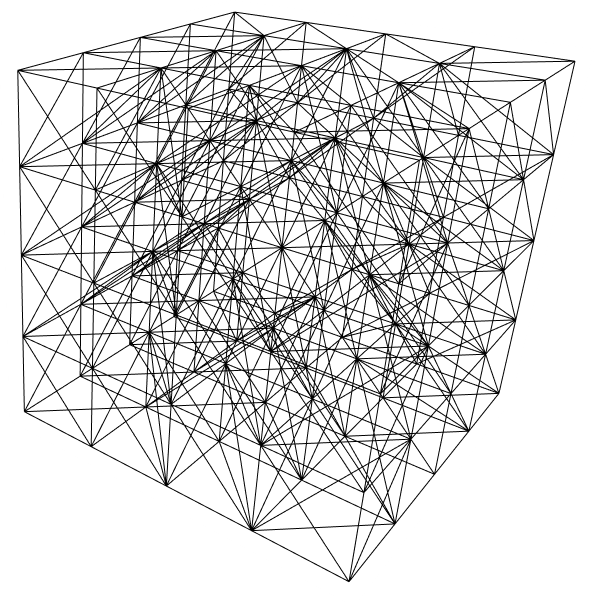}
    		\caption{}
    	\end{subfigure}
    	\begin{subfigure}{0.3\linewidth}
    		\centering
    		\includegraphics[width=1\linewidth]{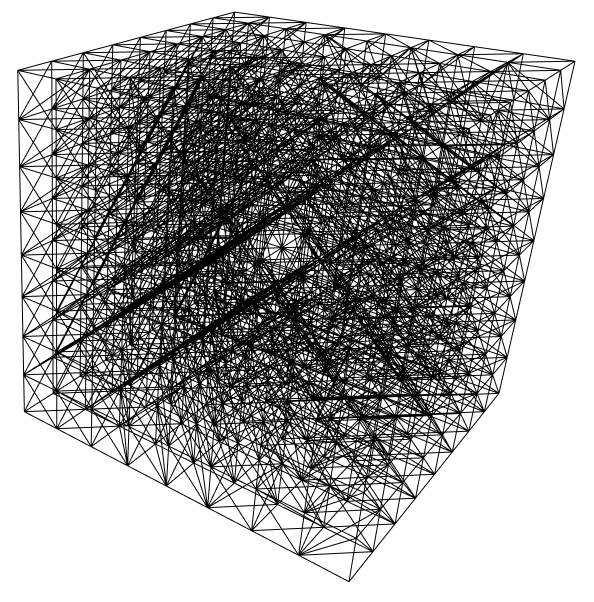}
    		\caption{}
    	\end{subfigure}
    	\caption{Displacement field of the Cauchy model on the coarsest mesh of 48 finite elements of the tenth order (a) and depictions of the meshes with 384 (b) and 3072 (c) elements, respectively.}
    	\label{fig:bounded2}
    \end{figure}

\begin{figure}
		\centering
		\begin{subfigure}{0.48\linewidth}
			\centering
			\definecolor{asl}{rgb}{0.4980392156862745,0.,1.}
\definecolor{asb}{rgb}{0.,0.4,0.6}
\begin{tikzpicture}
	\definecolor{npurple}{rgb}{0.4980392156862745,0.,1.}
	\begin{semilogxaxis}[
		/pgf/number format/1000 sep={},
		axis lines = left,
		xlabel={$\Lc$},
		ylabel={$I$},
		xmin=1e-4, xmax=1e+4,
		ymin=0, ymax=1,
		x dir=reverse,
		xtick={1e-3, 1e+0, 1e+3},
		ytick={0.1, 0.5, 0.9},
		legend style={at={(1,0.67)},anchor=east},
		ymajorgrids=true,
		grid style=dotted,
		]
		
		\addplot[color=asl, mark=triangle] coordinates {
			( 1000.0 ,  1.4327929268504584 )
			( 100.0 ,  0.6095494737008862 )
			( 31.622776601683793 ,  0.5414269049291108 )
			( 10.0 ,  0.5056151799627709 )
			( 5.623413251903491 ,  0.48630406079046395 )
			( 3.1622776601683795 ,  0.4621274043789848 )
			( 1.7782794100389228 ,  0.4321609929250506 )
			( 1.0 ,  0.39760034336986594 )
			( 0.5623413251903491 ,  0.3616034673158896 )
			( 0.31622776601683794 ,  0.3281135790197292 )
			( 0.1778279410038923 ,  0.3001818727076366 )
			( 0.1 ,  0.2790136877100162 )
			( 0.03162277660168379 ,  0.2544445937482407 )
			( 0.01 ,  0.24466938595392126 )
			( 0.001 ,  0.24014523927011808 )
			( 0.0001 ,  0.23967450966412682 )
		};
		\addlegendentry{$n_e = 384$}
		
		\addplot[color=violet, mark=o] coordinates {
			( 1000.0 ,  0.9067548116582305 )
			( 100.0 ,  0.8899121777631239 )
			( 31.622776601683793 ,  0.8797214426836334 )
			( 10.0 ,  0.8534465619317994 )
			( 5.623413251903491 ,  0.8282788525262215 )
			( 3.1622776601683795 ,  0.7915516138191729 )
			( 1.7782794100389228 ,  0.7419873407861056 )
			( 1.0 ,  0.681008725545968 )
			( 0.5623413251903491 ,  0.6134466947952979 )
			( 0.31622776601683794 ,  0.5464987836026256 )
			( 0.1778279410038923 ,  0.4870264804895582 )
			( 0.1 ,  0.4390910540685276 )
			( 0.03162277660168379 ,  0.3783123662058898 )
			( 0.01 ,  0.35060449080008393 )
			( 0.001 ,  0.3358418697980607 )
			( 0.0001 ,  0.33416539757789315 )
		};
		\addlegendentry{$n_e = 3072$}
		
		\addplot[color=asb, mark=square] coordinates {
			( 1000.0 ,  0.8723122918949631 )
			( 100.0 ,  0.868010686091502 )
			( 31.622776601683793 ,  0.8583413833154521 )
			( 10.0 ,  0.8311445605846346 )
			( 5.623413251903491 ,  0.8048745169487723 )
			( 3.1622776601683795 ,  0.7664708253698522 )
			( 1.7782794100389228 ,  0.714476585811559 )
			( 1.0 ,  0.6500273836027728 )
			( 0.5623413251903491 ,  0.577601222595007 )
			( 0.31622776601683794 ,  0.5042076540559882 )
			( 0.1778279410038923 ,  0.43700180519393117 )
			( 0.1 ,  0.3809166670275986 )
			( 0.03162277660168379 ,  0.30665803807216974 )
			( 0.01 ,  0.27137072027347037 )
			( 0.001 ,  0.2509494947562054 )
			( 0.0001 ,  0.24827875627326096 )
		};
		\addlegendentry{$n_e = 24576$}
		
		\addplot[color=blue, mark=pentagon] coordinates {
			( 1000.0 ,  0.8916200404165342 )
			( 100.0 ,  0.8606774959658097 )
			( 31.622776601683793 ,  0.8509668078206827 )
			( 10.0 ,  0.8236386716164382 )
			( 5.623413251903491 ,  0.7972378591170337 )
			( 3.1622776601683795 ,  0.7586317469199289 )
			( 1.7782794100389228 ,  0.7063264864679846 )
			( 1.0 ,  0.6413980456907228 )
			( 0.5623413251903491 ,  0.5682482280778567 )
			( 0.31622776601683794 ,  0.4938272135679687 )
			( 0.1778279410038923 ,  0.42529351172226915 )
			( 0.1 ,  0.3676730853203461 )
			( 0.03162277660168379 ,  0.2904250888490298 )
			( 0.01 ,  0.253232222175788 )
			( 0.001 ,  0.23147025022318554 )
			( 0.0001 ,  0.22847985074905167 )
		};
		\addlegendentry{$n_e = 48000$}
		
		\addplot[dashed,color=black, mark=none]
		coordinates {
			(1e+4, 0.16946198431835743)
			(1e-4, 0.16946198431835743)
		};
		
		\addplot[dashdotted,color=black, mark=none]
		coordinates {
			(1e+4, 0.8473099215917871)
			(1e-4, 0.8473099215917871)
		};
	\end{semilogxaxis}
	\draw (1.5,0.45) node[anchor=south]{$\mathbb{C}_\mathrm{macro}$};
	\draw (5.5,4.8) node[anchor=south]{$\mathbb{C}_\mathrm{micro}$};
\end{tikzpicture}
			\caption{}
		\end{subfigure}
	    \begin{subfigure}{0.48\linewidth}
	    	\centering
	    	\definecolor{asl}{rgb}{0.4980392156862745,0.,1.}
\definecolor{asb}{rgb}{0.,0.4,0.6}
\begin{tikzpicture}
	\definecolor{npurple}{rgb}{0.4980392156862745,0.,1.}
	\begin{semilogxaxis}[
		/pgf/number format/1000 sep={},
		axis lines = left,
		xlabel={$\Lc$},
		ylabel={$I$},
		xmin=1e-4, xmax=1e+4,
		ymin=0, ymax=1,
		x dir=reverse,
		xtick={1e-3, 1e+0, 1e+3},
		ytick={0.1, 0.5, 0.9},
		legend style={at={(1,0.71)},anchor=east},
		ymajorgrids=true,
		grid style=dotted,
		]
		
		\addplot[color=asl, mark=triangle] coordinates {
			( 1000.0 ,  0.9352021897201828 )
			( 100.0 ,  0.9322771153101306 )
			( 31.622776601683793 ,  0.9255311507001629 )
			( 10.0 ,  0.9063690103516739 )
			( 5.623413251903491 ,  0.887651454420386 )
			( 3.1622776601683795 ,  0.859971793085214 )
			( 1.7782794100389228 ,  0.8218983536469666 )
			( 1.0 ,  0.773599536427539 )
			( 0.5623413251903491 ,  0.7175785485133347 )
			( 0.31622776601683794 ,  0.6586983844926633 )
			( 0.1778279410038923 ,  0.6028418490398686 )
			( 0.1 ,  0.5548523573901727 )
			( 0.03162277660168379 ,  0.4900488603821017 )
			( 0.01 ,  0.46012663581923735 )
			( 0.001 ,  0.44503966750359325 )
			( 0.0001 ,  0.4434105752723326 )
		};
		\addlegendentry{$n_e = 384$}
		
		\addplot[color=violet, mark=o] coordinates {
			( 1000.0 ,  0.8560415931151273 )
			( 100.0 ,  0.8522333287897799 )
			( 31.622776601683793 ,  0.8434704222184762 )
			( 10.0 ,  0.8187442426015554 )
			( 5.623413251903491 ,  0.7947770735245078 )
			( 3.1622776601683795 ,  0.7595874255284596 )
			( 1.7782794100389228 ,  0.7116030386970928 )
			( 1.0 ,  0.6514432475461565 )
			( 0.5623413251903491 ,  0.5827525196611616 )
			( 0.31622776601683794 ,  0.5118047363280307 )
			( 0.1778279410038923 ,  0.4455082826911252 )
			( 0.1 ,  0.38906328562168563 )
			( 0.03162277660168379 ,  0.31240203566712466 )
			( 0.01 ,  0.2752857229631952 )
			( 0.001 ,  0.2548320046760102 )
			(0.0001, 0.2524529510681607)
		};
		\addlegendentry{$n_e = 3072$}
		
		\addplot[color=asb, mark=square] coordinates {
			( 1000.0 ,  0.847899086500748 )
			( 100.0 ,  0.8437635746268777 )
			( 31.622776601683793 ,  0.8342550963290448 )
			( 10.0 ,  0.8074737480557589 )
			( 5.623413251903491 ,  0.7815723344169285 )
			( 3.1622776601683795 ,  0.7436333371151003 )
			( 1.7782794100389228 ,  0.6920717742739054 )
			( 1.0 ,  0.6277179416425078 )
			( 0.5623413251903491 ,  0.5546013758640225 )
			( 0.31622776601683794 ,  0.47932859064829536 )
			( 0.1778279410038923 ,  0.40891685757509116 )
			( 0.1 ,  0.3485222782646808 )
			( 0.03162277660168379 ,  0.2645443462606974 )
			( 0.01 ,  0.22181345781335823 )
			( 0.001 ,  0.1962122074484653 )
			( 0.0001 ,  0.19287810761059862 )
		};
		\addlegendentry{$n_e = 24576$}
		
		\addplot[dashed,color=black, mark=none]
		coordinates {
			(1e+4, 0.16946198431835743)
			(1e-4, 0.16946198431835743)
		};
		
		\addplot[dashdotted,color=black, mark=none]
		coordinates {
			(1e+4, 0.8473099215917871)
			(1e-4, 0.8473099215917871)
		};
	\end{semilogxaxis}
	\draw (1.5,0.45) node[anchor=south]{$\mathbb{C}_\mathrm{macro}$};
	\draw (5.5,4.8) node[anchor=south]{$\mathbb{C}_\mathrm{micro}$};
\end{tikzpicture}
	    	\caption{}
	    \end{subfigure}
        \begin{subfigure}{0.48\linewidth}
        	\centering
        	\definecolor{asl}{rgb}{0.4980392156862745,0.,1.}
\definecolor{asb}{rgb}{0.,0.4,0.6}
\begin{tikzpicture}
	\definecolor{npurple}{rgb}{0.4980392156862745,0.,1.}
	\begin{semilogxaxis}[
		/pgf/number format/1000 sep={},
		axis lines = left,
		xlabel={$\Lc$},
		ylabel={$I$},
		xmin=1e-4, xmax=1e+4,
		ymin=0, ymax=1,
		x dir=reverse,
		xtick={1e-3, 1e+0, 1e+3},
		ytick={0.1, 0.5, 0.9},
		legend pos= south west,
		ymajorgrids=true,
		grid style=dotted,
		]
		
		\addplot[color=asl, mark=triangle] coordinates {
			( 1000.0 ,  0.843856326549687 )
			( 100.0 ,  0.8394962238517265 )
			( 31.622776601683793 ,  0.8297966781468089 )
			( 10.0 ,  0.8025436554273393 )
			( 5.623413251903491 ,  0.7762415308573348 )
			( 3.1622776601683795 ,  0.7378095517444951 )
			( 1.7782794100389228 ,  0.6857657311671885 )
			( 1.0 ,  0.6211468016218596 )
			( 0.5623413251903491 ,  0.5482555934163964 )
			( 0.31622776601683794 ,  0.47393951111784083 )
			( 0.1778279410038923 ,  0.40533455507262245 )
			( 0.1 ,  0.34755091934852583 )
			( 0.03162277660168379 ,  0.2701913341153719 )
			( 0.01 ,  0.2335135910913445 )
			( 0.001 ,  0.21354481897502536 )
			( 0.0001 ,  0.21119113904766879 )
		};
		\addlegendentry{$n_e = 384 \, , \; \Ber^3 \times \Ned^2$}
		
		\addplot[color=violet, mark=o] coordinates {
			( 1000.0 ,  0.8471521215569872 )
			( 100.0 ,  0.8428979575384361 )
			( 31.622776601683793 ,  0.8331210020091031 )
			( 10.0 ,  0.8056036508686608 )
			( 5.623413251903491 ,  0.779016082236531 )
			( 3.1622776601683795 ,  0.7401160894143561 )
			( 1.7782794100389228 ,  0.6873377106518075 )
			( 1.0 ,  0.6216247418367636 )
			( 0.5623413251903491 ,  0.5471916136380626 )
			( 0.31622776601683794 ,  0.4708121140453624 )
			( 0.1778279410038923 ,  0.3995807077884071 )
			( 0.1 ,  0.3386416301120596 )
			( 0.03162277660168379 ,  0.25414717698108663 )
			( 0.01 ,  0.21125079698219876 )
			( 0.001 ,  0.18556944119201188 )
			( 0.0001 ,  0.18217798515988448 )
		};
		\addlegendentry{$n_e = 3072 \, , \; \Ber^3 \times \Ned^2$}
		
		\addplot[color=asb, mark=square] coordinates {
			( 1000.0 ,  0.8438063693860487 )
			( 100.0 ,  0.8396531431933436 )
			( 31.622776601683793 ,  0.8303051721911059 )
			( 10.0 ,  0.8040502459663912 )
			( 5.623413251903491 ,  0.7787364395160762 )
			( 3.1622776601683795 ,  0.7417863871667447 )
			( 1.7782794100389228 ,  0.6918030905508379 )
			( 1.0 ,  0.6298110116431719 )
			( 0.5623413251903491 ,  0.5599998059117838 )
			( 0.31622776601683794 ,  0.489044354422035 )
			( 0.1778279410038923 ,  0.4238502757651924 )
			( 0.1 ,  0.3692319292642054 )
			( 0.03162277660168379 ,  0.29653601665361984 )
			( 0.01 ,  0.2623012903018341 )
			( 0.001 ,  0.2440473120699594 )
			( 0.0001 ,  0.24197069259582288 )
		};
		\addlegendentry{$n_e = 384 \, , \; \Ber^3 \times \Nedtwo^2$}
		
		\addplot[color=blue, mark=pentagon] coordinates {
			( 1000.0 ,  0.8471534446114187 )
			( 100.0 ,  0.8429146564527137 )
			( 31.622776601683793 ,  0.833172722324057 )
			( 10.0 ,  0.8057567305628952 )
			( 5.623413251903491 ,  0.7792706785201959 )
			( 3.1622776601683795 ,  0.7405243731577834 )
			( 1.7782794100389228 ,  0.6879618941776198 )
			( 1.0 ,  0.6225277878560909 )
			( 0.5623413251903491 ,  0.548430859533681 )
			( 0.31622776601683794 ,  0.4724474696822749 )
			( 0.1778279410038923 ,  0.40169578332449585 )
			( 0.1 ,  0.3413552957927539 )
			( 0.03162277660168379 ,  0.2583873842149924 )
			( 0.01 ,  0.21693671377670837 )
			( 0.001 ,  0.1926064623983449 )
			( 0.0001 ,  0.189489117630905 )
		};
		\addlegendentry{$n_e = 3072 \, , \; \Ber^3 \times \Nedtwo^2$}
		
		\addplot[dashed,color=black, mark=none]
		coordinates {
			(1e+4, 0.16946198431835743)
			(1e-4, 0.16946198431835743)
		};
		
		\addplot[dashdotted,color=black, mark=none]
		coordinates {
			(1e+4, 0.8473099215917871)
			(1e-4, 0.8473099215917871)
		};
	\end{semilogxaxis}
	\draw (5.5,0.45) node[anchor=south]{$\mathbb{C}_\mathrm{macro}$};
	\draw (5.5,4.8) node[anchor=south]{$\mathbb{C}_\mathrm{micro}$};
\end{tikzpicture}
        	\caption{}
        \end{subfigure}
        \begin{subfigure}{0.48\linewidth}
        	\centering
        	\definecolor{asl}{rgb}{0.4980392156862745,0.,1.}
\definecolor{asb}{rgb}{0.,0.4,0.6}
\begin{tikzpicture}
	\definecolor{npurple}{rgb}{0.4980392156862745,0.,1.}
	\begin{semilogxaxis}[
		/pgf/number format/1000 sep={},
		axis lines = left,
		xlabel={$\Lc$},
		ylabel={$I$},
		xmin=1e-4, xmax=1e+4,
		ymin=0, ymax=1,
		x dir=reverse,
		xtick={1e-3, 1e+0, 1e+3},
		ytick={0.1, 0.5, 0.9},
		legend pos= south west,
		ymajorgrids=true,
		grid style=dotted,
		]
		
		\addplot[color=asl, mark=triangle] coordinates {
			( 1000.0 ,  0.8519378346001252 )
			( 100.0 ,  0.8475127310549011 )
			( 31.622776601683793 ,  0.8376882972924693 )
			( 10.0 ,  0.8100715398373293 )
			( 5.623413251903491 ,  0.7833984055087697 )
			( 3.1622776601683795 ,  0.7443901506114571 )
			( 1.7782794100389228 ,  0.6915036835020355 )
			( 1.0 ,  0.6257376073636693 )
			( 0.5623413251903491 ,  0.5513940842201228 )
			( 0.31622776601683794 ,  0.4753429472237495 )
			( 0.1778279410038923 ,  0.4047403736139149 )
			( 0.1 ,  0.344734772874174 )
			( 0.03162277660168379 ,  0.2628441432800722 )
			( 0.01 ,  0.2229374302264403 )
			( 0.001 ,  0.2010069954749759 )
			( 0.0001 ,  0.19844855227628497 )
		};
		\addlegendentry{$\Ber^5 \times \Ned^4$}
		
		\addplot[color=violet, mark=o] coordinates {
			( 1000.0 ,  0.8475217761603515 )
			( 100.0 ,  0.8432664542828612 )
			( 31.622776601683793 ,  0.833486234499216 )
			( 10.0 ,  0.8059588470843813 )
			( 5.623413251903491 ,  0.779360386590064 )
			( 3.1622776601683795 ,  0.7404421774177246 )
			( 1.7782794100389228 ,  0.6876342372424231 )
			( 1.0 ,  0.6218750694652713 )
			( 0.5623413251903491 ,  0.5473727018734106 )
			( 0.31622776601683794 ,  0.4708936006155322 )
			( 0.1778279410038923 ,  0.39952378584256043 )
			( 0.1 ,  0.33839591727144347 )
			( 0.03162277660168379 ,  0.25330530755570535 )
			( 0.01 ,  0.20956748943723125 )
			( 0.001 ,  0.18287559348472582 )
			( 0.0001 ,  0.1793736785836238 )
		};
		\addlegendentry{$\Ber^7 \times \Ned^6$}
		
		\addplot[color=asb, mark=square] coordinates {
			( 1000.0 ,  0.8518816133605127 )
			( 100.0 ,  0.8475427427454422 )
			( 31.622776601683793 ,  0.8378083455355364 )
			( 10.0 ,  0.8104419965481727 )
			( 5.623413251903491 ,  0.78402196729697 )
			( 3.1622776601683795 ,  0.7454010729365357 )
			( 1.7782794100389228 ,  0.6930625423768626 )
			( 1.0 ,  0.6279973164019256 )
			( 0.5623413251903491 ,  0.5544556911313185 )
			( 0.31622776601683794 ,  0.479239496050169 )
			( 0.1778279410038923 ,  0.40948452112498934 )
			( 0.1 ,  0.3503980474322168 )
			( 0.03162277660168379 ,  0.2707197181682267 )
			( 0.01 ,  0.23290432814604872 )
			( 0.001 ,  0.21278579166927247 )
			( 0.0001 ,  0.21050111208407224 )
		};
		\addlegendentry{$\Ber^5 \times \Nedtwo^4$}
		
		\addplot[color=blue, mark=pentagon] coordinates {
			( 1000.0 ,  0.847582877529106 )
			( 100.0 ,  0.8432674589829707 )
			( 31.622776601683793 ,  0.8334888203777339 )
			( 10.0 ,  0.8059667165973848 )
			( 5.623413251903491 ,  0.7793738742811299 )
			( 3.1622776601683795 ,  0.740464797942429 )
			( 1.7782794100389228 ,  0.6876710142535634 )
			( 1.0 ,  0.6219325717730819 )
			( 0.5623413251903491 ,  0.54745904923473 )
			( 0.31622776601683794 ,  0.4710190812639736 )
			( 0.1778279410038923 ,  0.39970146890343783 )
			( 0.1 ,  0.3386409734009932 )
			( 0.03162277660168379 ,  0.25378108182629144 )
			( 0.01 ,  0.21058440679182933 )
			( 0.001 ,  0.18502577283364147 )
			( 0.0001 ,  0.1818045165659748 )
		};
		\addlegendentry{$\Ber^7 \times \Nedtwo^6$}
		
		\addplot[dashed,color=black, mark=none]
		coordinates {
			(1e+4, 0.16946198431835743)
			(1e-4, 0.16946198431835743)
		};
		
		\addplot[dashdotted,color=black, mark=none]
		coordinates {
			(1e+4, 0.8473099215917871)
			(1e-4, 0.8473099215917871)
		};
	\end{semilogxaxis}
	\draw (5.5,0.45) node[anchor=south]{$\mathbb{C}_\mathrm{macro}$};
	\draw (5.5,4.8) node[anchor=south]{$\mathbb{C}_\mathrm{micro}$};
\end{tikzpicture}
        	\caption{}
        \end{subfigure}
		\caption{Energy progression of the relaxed micromorphic model with respect to $\Lc$ using the linear (a), quadratic (b) and cubic (c) finite element formulations. The energy computed with the coarsest mesh of $48$ elements is depicted in (d) for various polynomial powers.}
		\label{fig:bounded1}
	\end{figure}
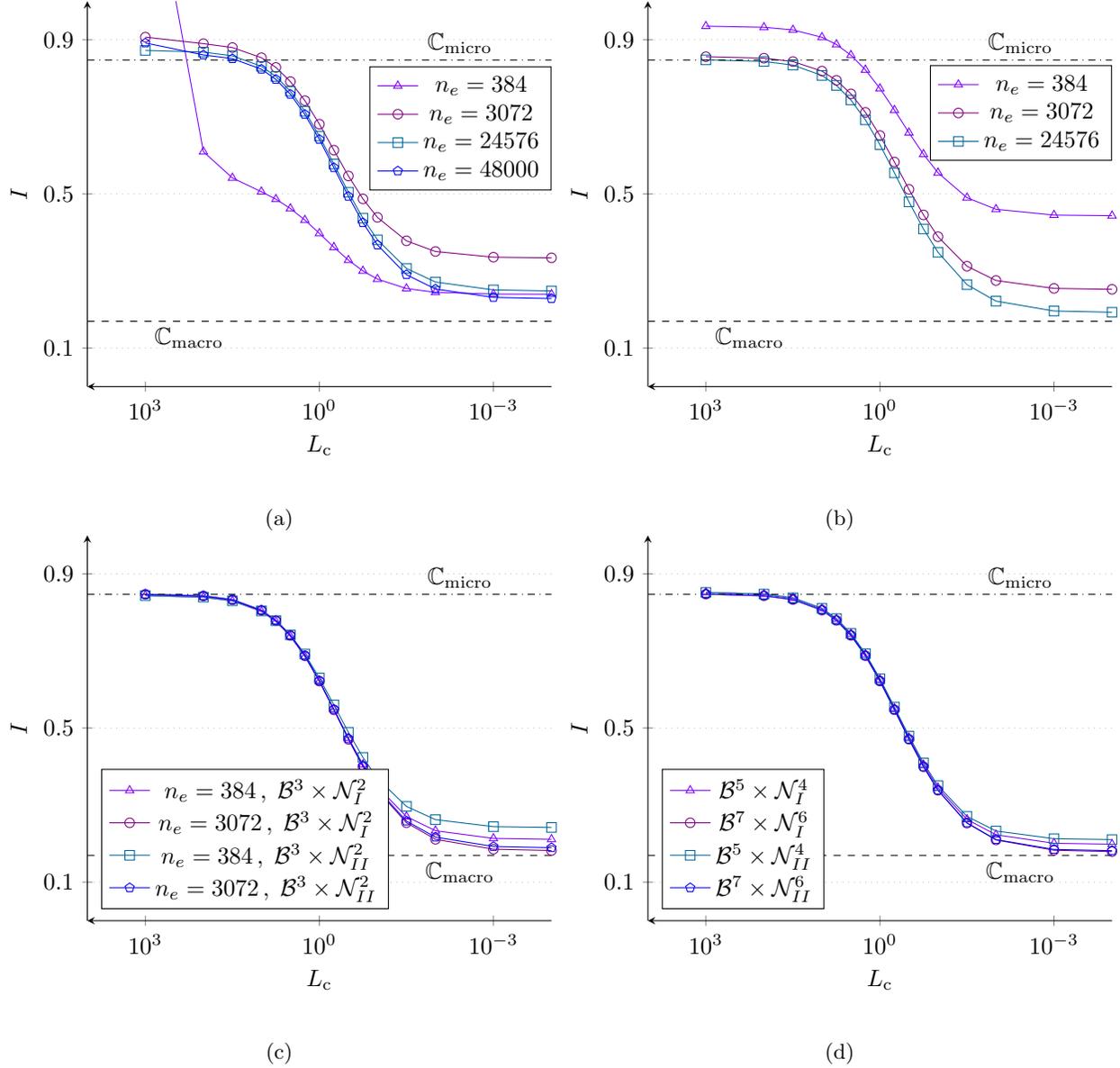

    The displacement field as well as some examples of the employed meshes are shown in \cref{fig:bounded2}. In order to compute the upper and lower bound on the energy we utilize the equivalent Cauchy model formulation with the micro- and macro elasticity parameters. In order to assert the high accuracy of the solution of the bounds we employ tenth order finite elements.
    The progression of the energy in dependence of the characteristic length parameter $\Lc$ is given in \cref{fig:bounded1}. We observe the high mesh dependency of the lower order formulations, where the energy is clearly overestimated. The higher order formulations all capture the upper bound correctly but diverge with respect to the result of the lower bound. Notably, the approximation using the N\'ed\'elec element of the first type is more accurate than the equivalent formulation with the N\'ed\'elec element of the second type, thus indicating the non-negligible involvement of the micro-dislocation in the energy. Using standard mesh coarseness the cubic element formulation with N\'ed\'elec elements of the first type yields satisfactory results. In order to achieve the same on highly coarse meshes, one needs to employ seventh order elements.

\section{Conclusions and outlook}
The intrinsic behaviour of the relaxed micromorphic model is revealed by the analytical solutions to boundary value problems. Clearly, the continuum exhibits hyperbolic and trigonometric solutions, which are not easily approximated by low order finite elements. The example provided in \cref{ex:2} demonstrates that cubic and higher order finite elements yield excellent results in approximate solutions of the model. 

The polytopal template methodology introduced in \cite{skypoly} allows to easily and flexibly construct $\Hc{}$-conforming vectorial finite elements that inherit many of the characteristics of an underlying $\Hone$-conforming basis, which can be chosen independently. 
In this work, we made use of Bernstein-B\'ezier polynomials. The latter boast optimal complexity properties manifesting in the form of sum factorization. The natural decomposition of their multi-variate versions into multiplications of univariate Bernstein base functions via the Duffy transformation allows to construct optimal iterators for their evaluation using recursion formulas. Further, this characteristic makes the use of dual numbers in the computation of their derivatives ideal. 
	Finally, the intrinsic order of traversal induced by the factorization is exploited optimally by the choice of clock-wise orientation of the reference element.  
	The consequence of these combined features is a high-performance hp-finite element program.
	
The ability of the relaxed micromorphic model to interpolate between the energies of homogeneous materials and materials with an underlying micro-structure using the characteristic length scale parameter $\Lc$ is demonstrated in \cref{ex:3}. It is also shown that in order to correctly capture the span of energies for the values of $\Lc$ either fine-discretizations or higher order elements are required.

The excellent performance of the proposed higher order finite elements in the linear static case is a precursor for their application in the dynamic setting, which is important since the relaxed micromorphic model is often employed in the computation of elastic waves (e.g., for acoustic metamaterials), where solutions for high frequency ranges are commonly needed.

The proposed computational scheme is lacking in its description of curved geometries. Due to the consistent coupling condition, this can easily lead to errors emanating from the boundary. Consequently, a topic for future works would be the investigation of curved finite elements \cite{Johnen2014,JOHNEN2013359} and their behaviour with respect to the model.

\section*{Acknowledgements}
Angela Madeo and Gianluca Rizzi acknowledge support from the European Commission through the funding of the ERC Consolidator Grant META-LEGO, N$^\circ$ 101001759.00

Patrizio Neff acknowledges support in the framework of the DFG-Priority Programme
2256 “Variational Methods for Predicting Complex Phenomena in Engineering Structures and Materials”, Neff 902/10-1, Project-No. 440935806.

\bibliographystyle{spmpsci}   

\footnotesize{
\bibliography{Ref}  
} 


\end{document}